\documentclass[12pt,a4paper]{article}
\usepackage[latin1]{inputenc}
\usepackage[french]{babel}
\usepackage{amsfonts}
\usepackage{amssymb}
\usepackage[all]{xypic}
\usepackage{bbm}

\newtheorem{lemme}{Lemme}[section]
\newtheorem{proposition}[lemme]{Proposition}
\newtheorem{definition}[lemme]{D\'efinition}
\newtheorem{theoreme}[lemme]{Th\'eor\`eme}
\newtheorem{corollaire}[lemme]{Corollaire}

\newcommand\dem{\noindent{\it D\'emonstration.}\ }
\newcommand\findem{\hfill$\square$}
\newcommand\hfld[2]{\ \smash{\mathop{\hbox to 7mm{\rightarrowfill}}
     \limits^{\scriptstyle#1}_{\scriptstyle#2}}\ }
\newcommand\hflg[2]{\smash{\mathop{\hbox to 7mm{\leftarrowfill}}
     \limits^{\scriptstyle#1}_{\scriptstyle#2}}}

\newcommand\ogg{\leavevmode\raise.3ex\hbox{$\scriptscriptstyle\langle\!\langle$}\,}
\newcommand\fgg{\leavevmode\raise.3ex\hbox{$\scriptscriptstyle\,\rangle\!\rangle$}}
\newcommand\rta{\rightarrow}

\newcommand\sv{{\scriptscriptstyle\vee}}
\newcommand\la{\langle}
\newcommand\ra{\rangle}

\newcommand\calA{{\mathcal A}}

\newcommand\calL{{\mathcal L}}
\newcommand\calM{{\mathcal M}}
\newcommand\calO{{\mathcal O}}

\newcommand\calJ{{\mathcal J}}

\newcommand\calT{{\mathcal T}}

\newcommand\frakg{{\mathfrak g}}
\newcommand\frakt{{\mathfrak t}}

\newcommand\FrakB{{\mathfrak B}}

\newcommand\build[3]{\mathrel{\mathop{\kern
0pt#1}\limits_{\textstyle #2}^{\textstyle #3}}}

\newcommand\bfB{{\mathbf B}}
\newcommand\CC{{\mathbbm C}}
\newcommand\GG{{\mathbbm G}}
\newcommand\FF{{\mathbbm F}}

\newcommand\ZZ{{\mathbbm Z}}
\newcommand\NN{{\mathbbm N}}
\newcommand\QQ{{\mathbbm Q}}

\newcommand\bbA{{\mathbbm A}}

\newcommand\rmH{{\rm H}}
\newcommand\rmR{{\rm R}}
\newcommand\bbX{{\mathbbm X}}
\newcommand\frakS{\mathfrak S}

\newcommand\GL{{\rm GL}}

\newcommand\Gal{{\rm Gal}}

\newcommand\pr{{\rm pr}}

\newcommand\car{{\bf car}}
\newcommand\Spec{{\rm Spec}}

\newcommand\SL{{\rm SL}}
\newcommand\PGL{{\rm PGL}}

\newcommand\Hom{{\rm Hom}}
\newcommand\Aut{{\rm Aut}}

\newcommand\Lie{{\rm Lie}}

\newcommand\ad{{\rm ad}}

\newcommand\act{{\rm act}}
\newcommand\cl{{\rm cl}}

\newcommand\Out{{\rm Out}}
\newcommand\ovl[1]{{\overline #1}}

\newcommand\hookr\hookrightarrow
\newcommand\hookl\hookleftarrow

\newcommand\inv{{\rm inv}}
\newcommand\isom{\,\smash{\mathop{\hbox to 5mm{\rightarrowfill}}
    \limits^\sim}\,}
\newcommand\cf{{\it cf.\ }}
\newcommand\wt[1]{{\widetilde #1}}
\newcommand\rmgssr{{\heartsuit}}
\newcommand\geom{{\rm geo}}

\begin{document}
\title{Fibration de Hitchin et endoscopie}
\author{Ng\^o Bao Ch\^au \footnote{Université Paris Sud, UMR 8628}}
\date{}
\maketitle

{\scriptsize{\it
\begin{flushright}
Mermoz, décidément, s'était retranché derrière son\\ ouvrage,
pareil au moissonneur qui, ayant
bien lié \\ sa gerbe, se couche  dans son champ.\\
\ \\
\rm Terre des hommes. {\bf Saint Exupéry}.
\end{flushright}}}

\begin{abstract}
We propose a geometric interpretation of the theory of elliptic
endoscopy, due to Langlands and Kottwitz, in terms of the Hitchin
fibration. As applications, we prove a global analog of a purity
conjecture, due to Goresky, Kottwitz and MacPherson. For unitary
groups, this global purity statement has been used, in a joint
work with G. Laumon, to prove the fundamental lemma over a local
fields of equal characteristics.
\end{abstract}

\bigskip
Nous proposons dans cet article une interprétation géométrique de
la théorie d'endoscopie elliptique, due à Langlands et Kottwitz, à
l'aide de la fibration de Hitchin.

L'espace de module $\calM$ des paires de Hitchin (\cf \cite{H})
dépend de la donnée d'une courbe projective lisse $X/\FF_q$ et
d'un fibré inversible $\calL$ sur $X$ et d'un schéma en groupes
réductifs $G$ au-dessus de $X$. Il paramètre les paires
$(E,\varphi)$ où $E$ est un $G$-torseur sur $X$ et où $\varphi$
est une section globale du fibré vectoriel $\ad(E)\otimes\calL$.
Le nombre de $\FF_q$-points de $\calM$ s'exprime formellement
comme le côté géométrique de la formule des traces sur l'algèbre
de Lie pour une certaine fonction test très simple. Le découpage
de la formule des traces en classes de conjugaison stables
correspond de façon évidente au morphisme de Hitchin
$f:\calM\rta\bbA$ associant à une paire $(E,\varphi)$ le polynôme
caractéristique de l'endomorphisme infinitésimal $\varphi$.
L'espace affine de Hitchin $\bbA$ dont la dimension dépend du
triplet $(X,G,\calL)$, paramètre ces polynômes caractéristiques.
La somme des intégrales orbitales globales dans une classe de
conjugaison stable donnée $a\in\bbA(k)$, compte formellement le
nombre de points rationnels de la fibre de Hitchin $f^{-1}(a)$ ;
ces nombres pouvant être infinis pour les caractéristiques $a$ non
elliptiques.

Pour un groupe $G$ semi-simple, il y a un ouvert dense $\bbA^{\rm
ell}$ de $\bbA$ au-dessus duquel le morphisme $f$ est propre, en
particulier de type fini et où les intégrales orbitales consistent
en des sommes finies. Cet ouvert contient les classes de
conjugaison stable rationnelles qui sont géométriquement
elliptiques c'est-à-dire les classes qui restent elliptiques après
toute extension du corps des constants. Dans ce travail, nous nous
sommes limités à cet ouvert elliptique. Notre notion d'ellipticité
est un peu différente de la notion standard en ce qu'elle concerne
uniquement le tore maximal et non le centre de $G$.

Du côté automorphe, Langlands \cite{RL} et Kottwitz \cite{K-EST}
ont réécrit la partie elliptique de la formule des traces en une
somme où on voit apparaître une première sommation sur l'ensemble
des classes de conjugaison stables des tores maximaux $T$ de $G$
définis sur le corps global de $X$, puis une deuxième sommation
sur les éléments $\kappa\in \hat T^\Gamma$ où $\hat T$ est le tore
dual complexe de $T$, muni d'une action finie du groupe de Galois
$\Gamma$ du corps $F$ des fonctions rationnelles de $X$. L'objet
de cet article est de donner une interprétation géométrique de
cette formule, autrement dit de trouver une décomposition de
l'image directe $f_*\QQ_\ell$, du moins de sa restriction à
l'ouvert elliptique.

De notre point de vue, il est naturel d'intervertir l'ordre dans
la sommation ci-dessus. On trouvera d'abord une {\em décomposition
canonique de la cohomologie perverse de la fibration de Hitchin en
somme directe}
$$^p\rmH^i (f^{\rm ell}_* \QQ_\ell) =\bigoplus_{[\kappa]}
\,^p \rmH^i(f^{\rm ell}_* \QQ_\ell)_{[\kappa]}.
$$
qu'on appellera la $[\kappa]$-décomposition. Cette somme s'étend
sur l'ensemble des classes de $W'$-conjugaison $[\kappa]$ des
caractères d'ordre fini $\kappa:\bbX^\vee\rta\CC^\times$ où on a
noté $\bbX^\vee$ le groupe des cocaractères du tore maximal
abstrait de $G$ et où $W'=W\rtimes \Theta$, $W$ étant le groupe de
Weyl, $\Theta$ étant le groupe de Galois d'un revêtement fini
étale de $X$ sur lequel $G$ devient une forme intérieure d'un
groupe constant.

Les relations précises entre $T$ et $\kappa$, comprenant notamment
$\kappa\in \hat T^\Gamma$, sont inhérentes à la formule de
Langlands et Kottwitz. Elles se traduiront, dans notre approche,
en {\em une estimation du support de la partie
$[\kappa]$-isotypique $^p \rmH^i(f_* \QQ_\ell)_{[\kappa]}$}. On
démontre que la partie $\kappa$-isotypique est supportée par le
réunion des images dans $\bbA$ des espaces affines de Hitchin
$\bbA_H$ des {\em groupes endoscopiques non ramifiés} $H$ associés
à $\kappa$, voir le théorème 10.4 qui est le résultat principal de
cet article.

Dans la démonstration récente du lemme fondamental pour les
groupes unitaires \cite{L-N}, cette description du support a joué
le rôle d'un énoncé de pureté analogue à une conjecture de
Goresky, Kottwitz et MacPherson \cite{GKM}. On pourrait s'attendre
à ce qu'il en est de même en général, voir la discussion à la fin
de §10.

L'ingrédient de base pour construire la $[\kappa]$-décomposition
est l'existence d'une action d'un champ de Picard relatif $P/\bbA$
sur la fibration de Hitchin $\calM/\bbA$. Par exemple, si
$G=\GL(n)$, on peut associer à toute caractéristique $a\in \bbA$,
un revêtement spectral $Y_a\rta X$ fini de degré $n$ qui est
génériquement étales pour des caractéristiques $a\in \bbA$
génériques. D'après Hitchin \cite{H} et Beauville, Narasimhan et
Ramanan \cite{BNR}, la fibre $f^{-1}(a)=\calM_a$ s'identifie à la
jacobienne compactifiée de $Y_a$. Dans ce cas, le champ de Picard
$P_a$ n'est autre que la jacobienne de $Y_a$ et l'action de $P_a$
sur $\calM_a$ est l'action évidente de la jacobienne sur la
jacobienne compactifiée par produit tensoriel. Dans le cas des
groupes classiques, on peut également deviner $P$ et son action
sur $\calM$. La définition uniforme du champ de Picard $P/\bbA$
est fondée sur le schéma en groupes $J$ des centralisateurs
réguliers de $G$.

Le champ de Picard relatif $P/\bbA$ agit sur la fibration de
Hitchin et par conséquent, sur les faisceaux de cohomologie
$^p\rmH^i(f_*\QQ_\ell)$. Le lemme d'homotopie (\cf \cite{L-N} 3.2)
nous dit que cette action se factorise à travers une action du
faisceau $\pi_0(P/\bbA)$ dont la fibre en un point géométrique
$a\in \bbA$ est le groupe des composantes connexes $\pi_0(P_a)$.

Sur l'ouvert $\bbA^{\rm ell}$ de $\bbA$, le faisceau
$\pi_0(P/\bbA)$ est un faisceau en groupes abéliens finis. Dans un
voisinage d'un point $a$ elliptique, $^p\rmH^i(f_*\QQ_\ell)$ se
dé\-com\-pose selon les caractères du groupe fini $\pi_0(P_a)$ et
c'est ainsi qu'on peut faire apparaître les caractères
endoscopiques $\kappa$.

L'action du champ de Picard $P$ sur la fibration de Hitchin peut
être considérée comme l'interprétation géométrique de l'opération
d'extraction de l'intégrale orbitale globale un nombre de
Tamagawa. L'intégrale orbitale globale divisée par ce volume, se
décompose alors en un produit d'intégrales orbitales locales.
Cette formule de produit s'interprète géo\-mé\-tri\-quement : le
quotient d'une fibre de Hitchin $\calM_a$ par le champ de Picard
fibre $P_a$ se décompose en un produit. A partir de cette formule
produit, on retrouve l'obstruction cohomologique, construite par
Langlands et Kottwitz, pour qu'une collection des classes de
conjugaison locales, toutes dans une classe de conjugaison stable
donnée, se recollent en une classe de conjugaison rationnelle :
l'obstruction se trouve tautologiquement dans $\rmH^1(k,P_a)$.

\begin{small}
\medskip {\em Remerciements.} J'ai été très
influencé par les travaux de Kottwitz sur l'endoscopie, tout
particulièrement dans le calcul crucial de $\pi_0(P_a)$.

J'ai profité de nombreuses conversations avec G. Laumon sur la
fibration de Hitchin et sur le lemme fondamental. Quand je
commençais à étudier la fibration de Hitchin et le lien avec
l'endoscopie, il m'a suggéré de compter le nombre de points de
l'espace de module de Hitchin : ce comptage s'est relevé
particulièrement inspirant pour la suite. Sa vigilance m'a
également permis de rectifier et d'améliorer l'énoncé de la
$[\kappa]$-décomposition.

Je remercie L.~Breen et M.~Raynaud pour leur aide précieuse dans
la réalisation de ce travail et tout particulièrement A. Abbes qui
m'a indiqué une référence utile à EGA IV.3. Par ailleurs, la tenue
d'un groupe de travail sur l'endoscopie et le lemme fondamental au
printemps 2003, à l'IHES et à l'Université Paris 13, m'a été une
grande aide. Je saisis l'occasion pour remercier ses participants,
en particulier P.-H. Chaudouard, J.-F. Dat, L. Fargues, A.
Genestier et L. Lafforgue. Je remercie J.-F. Dat et D. Gaitsgory
et tout particulièrement le referee pour leurs lectures attentives
du manuscript qui ont épargné à cet article quelques erreurs et O.
Gabber pour m'avoir signalé une restriction nécessaire sur la
caractéristique. Je remercie V. Drinfeld qui m'a suggéré l'énoncé
de la proposition 6.8 qui aura rendu plus clair le fil de
l'article.

\end{small}

\section{Comptage des paires de Hitchin}

Soit $k=\FF_q$ un corps fini et $\ovl k$ sa clôture algébrique.
Dans tout l'article, on va supposer que la caractéristique de $k$
est suffisamment grande. On en donnera une borne précise dans §2.

Soit $X$ une courbe projective lisse géométriquement connexe sur
$k$. Notons $\eta=\Spec(F)$ son point générique, $F$ étant le
corps des fonctions de $X$ et notons $|X|$ l'ensemble des points
fermés de $X$. Pour toute $v\in|X|$, notons $F_v$ la complétion de
$F$ pour la topologie $v$-adique. On va fixer un diviseur effectif
$D=\sum_{v\in |X|} d_v [v]$ de degré $d=\sum_{v\in|X|}
d_v\deg(v)$. Soit $\calO_X(D)$ le fibré en droites associé à ce
diviseur et $L_D$ le $\GG_m$-torseur sur $X$ associé à ce fibré en
droites. On supposera que $\deg(D)>2g-2$.

Soit $\GG$ un groupe réductif connexe {\em déployé} \footnote{Tout
groupe réductif connexe sur un corps fini est quasi-déployé et
cela suffit pour nous. Néanmoins, comme $\GG$ sert seulement comme
modèle auquel $G$ est localement isomorphe à pour la topologie
étale, il n'y a pas de raison de ne pas supposer $\GG$ déployé.}
sur $k$ et $G$ un schéma en groupes sur $X$ qui localement pour la
topologie étale de $X$, est isomorphe à $X\times \GG$. Le groupe
des automorphismes de $\GG$ est un groupe algébrique sur $k$ qui
s'insère dans la suite exacte
$$1\rta \GG^{\rm ad}\rta {\rm Aut}(\GG)\rta {\rm Out}(\GG)\rta 1$$
où $\GG^{\rm ad}$ est le quotient de $\GG$ par son centre et où
${\rm Out}(\GG)$ est un groupe discret, \cf \cite{Spr}. Le groupe
$\GG$ a des épinglages définis sur $k$, fixons-en un. Le
sous-groupe de ${\rm Aut}(\GG)$ fixant ce scindage s'envoie
bijectivement sur ${\rm Out}(\GG)$ et induit donc un scindage
${\rm Out}(\GG)\rta {\rm Aut}(\GG)$ de la suite exacte ci-dessus.

Soit $x$ un point géométrique de $X$. La donnée de $G$ induit un
${\rm Aut}(\GG)$-torseur $\tau_G$ sur $X$ dont la classe
d'isomorphisme est un élément
$$[\tau_G]\in \rmH^1(\pi_1(X,x),{\rm Aut}(\GG)(\ovl k)).$$
En utilisant la projection ${\rm Aut}(\GG)\rta {\rm Out}(\GG)$ on
obtient un ${\rm Out}(G)$-torseur $\tau^{\rm Out}_G$. La classe
d'isomorphisme $[\tau^{\Out}_G]\in \rmH^1(\pi_1(X,x),\Out(\GG))$
consiste en un homomorphisme $$\rho_G:\pi_1(X,x)\rta\Out(\GG)$$
puisque $\GG$ étant déployé, l'action de $\pi_1(X,x)$ sur
$\Out(\GG)$ via $\Gal(\ovl k/k)$ est triviale. Nous allons noter
$\Theta$ l'image de $\pi_1(X,x)$ dans $\Out(\GG)$. L'homomorphisme
surjectif $\rho_G:\pi_1(X,x)\rta\Theta$ définit un revêtement
$X_\Theta$ de $X$ étale galoisien de groupe de Galois $\Theta$ et
pointé par un point géométrique $x_\Theta$ de $X_\Theta$.

Le relèvement ${\rm Out}(\GG)\rta {\rm Aut}(\GG)$ construit à
partir du scindage, nous donne à nouveau un ${\rm
Aut}(\GG)$-torseur $\tau_G^{\rm ext}$. Ce torseur correspond à un
schéma en groupes quasi-déployé $G^{\rm ext}$ sur $X$ dont $G$ est
une forme intérieure. Le schéma en groupes $G$ est quasi-déployé
si les $\Aut(\GG)$-torseurs $\tau_G$ et $\tau_G^{\rm ext}$ sont
isomorphes.

On appelle {\em paires de Hitchin}, par rapport à $X,G$ et
$\calO_X(D)$, les paires $(E,\varphi)$ où $E$ est un $G$-torseur
sur $X$ et où $\varphi\in \rmH^0(X,\ad(E)\otimes \calO_X(D))$,
$\ad(E)$ étant le fibré vectoriel sur $X$ obtenu en tordant le
fibré vectoriel ${\rm Lie}(G)$ muni de l'action adjointe de $G$
par le $G$-torseur $E$. On considère le champ $\calM$ classifiant
ces paires. Il associe à tout $\FF_q$-schéma $S$ la catégorie en
groupoïdes des paires $(E,\varphi)$ où $E$ est un $G$-torseur sur
$X\times S$ et où $\varphi\in\rmH^0(\ad(E)\otimes \pr_X^*
\calO_X(D))$, $\pr_X:X\times S\rta X$ étant la projection évidente
sur $X$. C'est un champ algébrique.

On va décrire de façon adélique la catégorie $\calM(k)$ des paires
de Hitchin définies sur $k$. Soit $(E,\varphi)$ un objet de
$\calM(k)$. La classe d'isomorphisme de la fibre générique de $E$
définit un élément $\cl_\eta(E)$ de $\rmH^1(F,G)$.

\begin{lemme}
Soit $E$ un $G$-torseur sur $X$. L'élément $\cl_\eta(E)$ de
$\rmH^1(F,G)$ qui correspond à la classe d'isomorphisme du
$G$-torseur $E_\eta$, a l'image triviale dans $\rmH^1(F_v,G)$ pour
toutes les places $v$ de $X$. Inversement soit $c$ un élément de
$${\rm ker}^1(F,G)={\ker }[\rmH^1(F,G)\rta \prod_v \rmH^1(F_v,G)].$$
Alors, il existe un $G$-torseur $E$ sur $X$ dont la classe
d'isomorphisme de la fibre générique $\cl_\eta(E)$ est l'élément
$c$.
\end{lemme}

\dem Soit $\cl_v(E)$ l'image de $\cl_\eta(E)$ dans
$\rmH^1(F_v,G)$. Cet élément est aussi l'image de l'élément de
$\beta_v\in \rmH^1(\Spec(\calO_v),G)$, la classe d'isomorphisme du
$G$-torseur $E_{\calO_v}$ sur $\Spec(\calO_v)$. Ce $\rmH^1$ étant
trivial par un théorème de Lang, on déduit que l'élément
$\cl_v(E)\in\rmH^1(F_v,G)$ est trivial.

Inversement, soit $E_\eta$ un $G$-torseur sur le point générique
$\eta$ de $X$ ayant la classe d'isomorphisme $c\in\ker^1(F,G)$. On
peut prolonger $E_\eta$ en un $G$-torseur $E_U$ sur un ouvert $U$
de $X$. Sur un point $v$ du complémentaire de $U$, on peut
recoller $E_U$ avec le $G$-torseur triviale sur $\Spec(\calO_v)$,
parce que la restriction à $F_v$ de $E_U$ est un $G$-torseur
trivial. On obtient ainsi un $G$-torseur sur $X$ dont la classe
d'isomorphisme de la fibre générique est l'élément $c$. \hfill
$\square$

\medskip
Deux $G$-torseurs sur $X$ seront dits isogènes s'ils ont des
fibres génériques isomorphes. On a démontré que l'ensemble des
classes d'isogénie de $G$-torseur sur $X$ est en bijection
canonique avec le sous-ensemble $\ker^1(F,G)$ de $\rmH^1(F,G)$
constitué des éléments localement triviaux. Dans le cas où $G$ est
un groupe semi-simple adjoint ou simplement connexe, d'après un
théorème de Kneser, cet ensemble est réduit à l'élément neutre de
$\rmH^1(F,G)$. Ce n'est pas le cas en général. Notons au passage
que l'ensemble de ${\rm ker}^1(F,G)$ apparaît également dans le
comptage des points des variétés de Shimura (\cf \cite{K-JAMS}).

Pour tous $c\in {\rm ker}^1(F,G)$, on choisit un $G$-torseur
modèle ${\mathbbm E}_c$ sur $X$ dont la fibre générique a la
classe d'isomorphisme $c$ et qui est muni des trivialisations sur
tous les traits $\Spec(\calO_v)$ avec $v\in |X|$. Soit $G_c$ le
schéma en groupes des automorphismes de ${\mathbbm E}_c$ et
$\Lie(G_c)$ le schéma en algèbres de Lie associé. Les
trivialisations locales choisies induisent des isomorphismes entre
les changements de base à $\Spec(\calO_v)$ de $G$ et de $G_c$.

Considérons donc la catégorie $\mathcal R_c$ des triplets
$(E,\varphi,\iota)$ formée d'une paire de Hitchin
$(E,\varphi)\in\calM(k)$ et d'une rigidification de la fibre
générique $E_\eta$ $\iota:E_\eta \isom {\mathbbm E}_{c,\eta}$. Le
groupe des automorphismes d'un tel triplet est trivial puisqu'un
automorphisme de $E$ qui fixe la trivialisation générique $\iota$
est nécessaire\-ment l'identité si bien que la catégorie $\mathcal
R_c$ est discrète.

\begin{proposition}
On a une équivalence de catégories entre la catégorie $\mathcal
R_c$ des triplets $(E,\varphi,\iota)$ comme ci-dessus, et la
catégorie discrète des uplets $$(\gamma,(g_v)_{v\in|X|})$$
constitués
\begin{itemize}
\item d'un élément $\gamma\in\Lie(G_c)(F)$,
\item d'un $g_v\in G(F_v)/G(\calO_v)$ pour chaque
$v\in|X|$ tel que pour presque toutes $v$, $g_v$ est la classe
neutre,
\end{itemize}
qui vérifient $\ad(g_v)^{-1}(\gamma) \in
\varpi_v^{-d_v}\Lie(G)(\calO_v)$ pour tous $v\in|X|$. Ici,
$\varpi_v$ est une uniformisante de $F_v$.
\end{proposition}

\dem La donnée de la trivialisation $\iota$ induit un isomorphisme
de $F$-algèbres de Lie $\ad(E_\eta)\isom\Lie(G_c)_\eta$. La
section $\varphi\in \rmH^0(X,\ad(E)\otimes \calO_X(D))$ définit un
élément $\gamma\in \Lie(G_c)(F)$ puisque le fibré en droites
$\calO_X(D)$ est lui aussi trivialisé génériquement.

Puisqu'on s'est donné des trivialisations de ${\mathbbm E}_c$ sur
les traits $\Spec(\calO_v)$, le $G$-torseur $E$ muni d'un
isomorphisme avec ${\mathbbm E}_c$ au-dessus du point générique de
$X$, est équivalente à la donnée pour tous $v\in|X|$ d'un élément
$g_v\in G(F_v)/G(\calO_v)$ lesquels sont presque tous neutres.
Dans cette bijection
$$(E,\iota)\leftrightarrow (g_v)_{v\in |X|},$$
via l'identification des fibres génériques
$\ad(\iota):\ad(E)_\eta\isom \frakg_{c,\eta}$, les réseaux
$\ad(E)(\calO_v)$ et $\ad(g_v) \frakg(\calO_v)$ se correspondent.
Ainsi, pour que $\gamma\in\frakg_c(F)$ définit une section globale
de $\ad(E)\otimes \calO_X(D)$, il faut et il suffit que pour tous
$v\in |X|$, on a
$$\gamma\in \varpi_v^{-d_v} \ad(g_v) \frakg(\calO_v)$$
ce qui est équivalent à la condition $\ad(g_v)^{-1}\gamma\in
\varpi_v^{-d_v} \frakg(\calO_v)$. \findem

\begin{corollaire}
On a une équivalence de catégories entre $\calM(k)$ et la
catégorie dont les objets sont les triples
$(c;\gamma,(g_v)_{v\in|X|})$ avec
\begin{itemize}

\item $c\in \ker^1(F,G)$,

\item $\gamma\in\Lie(G_c)(F)$, \item pour tout $v\in|F|$ une
classe $g_v\in G(F_v)/G(\calO_v)$ qui est presque partout
triviale,

\item $\ad(g_v)^{-1}\gamma \in \varpi_v^{-d_v} \frakg(\calO_v)$
\end{itemize}
et dont l'ensemble des flèches d'un objet
$(c;\gamma,(g_v)_{v\in|X|})$ sur un autre objet
$(c';\gamma',(g'_v)_{v\in|X|})$ est
\begin{itemize}
\item vide si $c\not=c'$ ;
\item l'ensemble des éléments $\xi\in
G_c(F)$ tels que $\gamma'=\ad(\xi)\gamma$ et $g'_v=\xi g_v$ si
$c=c'$.
\end{itemize}
\end{corollaire}

Cet énoncé nous conduit à écrire une égalité qui pour l'instant
n'a pas de sens car les sommes en jeu sont infinies
$$|\calM(k)|=\sum_{\alpha\in \ker^1(F,G)}
\sum_{\gamma\in\frakg(\alpha)(F)/{\rm conj.}} O_\gamma(1_D)
$$
où $O_\gamma(1_D)$ est l'intégrale orbitale globale
$$
O_\gamma(1_D)=\int_{G_\gamma(F)\backslash G(\bbA_F)}
1_{D}(\ad(g)^{-1}\gamma) dg
$$
de la fonction
$$1_D=\bigotimes_{v\in |X|} 1_{\varpi_v^{-d_v}\frakg(\calO_v)},$$
$1_{\varpi_x^{-d_v}\frakg(\calO_v)}$ étant la fonction
caractéristique du compact ouvert $\varpi^{-d_x}\frakg(\calO_x)$
de $\frakg(F_v)$, le produit tensoriel sur toutes les places $v\in
|X|$ étant le produit tensoriel restreint. Cette égalité sera un
très bon guide pour la suite de l'article.

\section{Morphisme de Hitchin et section de Kostant}

On supposera désormais que la caractéristique de $k$ ne divise pas
l'ordre du groupe de Weyl $W$ de $\GG$.

Dans cette section, on rappelle la construction, due à Hitchin, du
morphisme qui associe à une paire de Hitchin, son {\ogg}polynôme
caractéristique\fgg. On commencera par rappeler le théorème de
restriction de Chevalley. La construction présentée ici qui est
bien entendu la même que la construction originale de Hitchin,
peut sembler singulièrement encombrante. Ceci est du en partie à
notre désir de traiter uniformément des schémas groupes réductifs
$G$ sur $X$. Par ailleurs, le langage de champs qui semble lourd
ici se montrera efficace dans la suite de l'article.

De règle générale, on commence par une construction sur le groupe
constant $\GG$, on vérifie que la construction est
$\Aut(\GG)$-équivariante, puis on le tord par le $\Aut(G)$-torseur
$\tau_G$ pour passer à $G$. On utilisera un indice $\GG$ ou la
police \ogg\ double-barre \fgg\ pour désigner les objets attachés
à $\GG$.

Dans §1, on a fixé un épinglage de $\GG$ défini sur $k$. Cet
épinglage consiste en un tore maximal $\mathbbm T$ de $\mathbbm
G$, un sous-groupe de Borel $\mathbbm B$ contenant $\mathbbm T$
tous les deux définis sur $k$, et enfin pour toute racine simple
$\alpha\in\Delta$, un vecteur non nul $\bf x_\alpha$ de l'espace
propre de $\mathbbm g_\alpha$ de ${\mathbbm g}=\Lie(\GG)$ où
$\mathbbm T$ agit par le caractère $\alpha$. L'hypothèse $\GG$
déployé implique que ces vecteurs sont définis sur $k$. On note
${\bf x}_+=\sum_{\alpha\in\Delta}{\bf x}_\alpha$. Pour toute
$\alpha\in\Delta$, on peut choisir de façon unique un vecteur non
nul ${\bf x}_{-\alpha}$ dans l'espace propre de $\mathbbm T$ dans
$\mathbbm g$ de caractère $-\alpha$ tel que $[{\bf x_\alpha,\bf
x_{-\alpha}}]=\alpha^\vee$. On note $\bf
x_-=\sum_{\alpha\in\Delta}{\bf x_{-\alpha}}$. Le sous-groupe de
$\Aut(\GG)$ fixant cet épinglage s'envoie bijectivement sur le
groupe $\Out(\GG)$ de sorte que $\Out(\GG)$ agit canoniquement sur
$\GG$ en fixant l'épinglage. On en déduit une action de
$\Aut(\GG)$ sur $\GG$ à travers $\Out(\GG)$ qui est différente de
l'action tautologique de $\Aut(\GG)$ sur $\GG$.

Notons $k[\mathbbm g]$ l'algèbre des fonctions polynômiales sur
$\mathbbm g$. L'action adjointe de $G$ sur $\mathbbm g$ induit une
action de $\GG$ sur cette algèbre $k[\mathbbm g]$. Soit
$k[\mathbbm g]^{\GG}$ la sous-algèbre de $k[\mathbbm g]$ des
fonctions invariantes sous l'action adjointe $\GG$. Notons
$\mathbbm t$ l'algèbre de Lie de $\mathbbm T$. Soit $k[\mathbbm
t]^W$ la sous-algèbre des $W$-invariants de $k[\mathbbm t]$. On
note $\car_\GG=\Spec(k[\mathbbm t]^W)$. Le groupe ${\rm Out}(\GG)$
agit sur $k[\mathbbm t]$ et sur $W$ donc sur ${\car_\GG}$. On en
déduit une action de $\Aut(\GG)$ sur $\car_\GG$ à travers son
quotient ${\rm Out}(\GG)$. L'énoncé suivant a été démontré par
Kostant en caractéristique nulle et étendu par Veldkamp aux
caractéristiques positives suffisamment grandes par rapport au
système de racines, voir \cite{Kos} et \cite{Vel}.

\begin{theoreme}
Supposons que $k$ est un corps de caractéristique ne divisant pas
l'ordre du groupe de Weyl $|W|$.
\begin{enumerate}
\item
La restriction de $\mathbbm g$ à $\mathbbm t$ induit un
isomorphisme $k[\mathbbm g]^{\GG}\isom k[\mathbbm t]^W$. De plus,
$k[\mathbbm t]^W$ est une algèbre de polynômes en les variables
$u_1,\ldots,u_r$ qui sont homogènes de degrés $m_1+1,\ldots,m_r+1$
qui ne dépendent pas des choix des $u_i$. On a donc un morphisme
$\GG_m$-équivariant $\chi:\mathbbm g\rta \car_\GG$ où
$\car_\GG=\Spec(k[\mathbbm t]^W)$ muni de l'action de $\GG_m$
défini par
$t.(u_1,\ldots,u_r)=(t^{m_1+1}u_1,\ldots,t^{m_r+1}u_r)$.

\item Soit $\mathbbm g^{\rm reg}$ l'ouvert de $\mathbbm g$ des éléments
dont la dimension du centralisateur est égale à la dimension de
$r$. Alors le morphisme restreint $\mathbbm g^{\rm
reg}\rta\car_\GG$ est un morphisme lisse dont les fibres sont des
orbites sous $G$.

\item Soit $\mathbbm g^{\bf x_+}$ l'algèbre de Lie du centralisateur
de ${\bf x_+}$. Alors, le sous-espace affine ${\bf x}_-+\mathbbm
g^{\bf x_+}$ est contenu dans $\mathbbm g^{\rm reg}$ et de plus,
la restriction de $\chi$ à cet espace définit un isomorphisme
${\bf x}_-+\mathbbm g^{\bf x_+} \rta \car_\GG$. Son inverse est
une section de Kostant $\epsilon:\car_\GG\rta \mathbbm g^{\rm
reg}$.
\end{enumerate}
\end{theoreme}

Nous avons besoin d'un énoncé sur la section de Kostant qui tient
compte de l'action de $\GG_m$. Considérons l'action $\rho(\GG_m)$
sur $\mathbbm g$ qui agit trivialement sur $\mathbbm t$ et pour
toute racine $\alpha$, $\rho(\GG_m)$ agit sur l'espace propre
$\mathbbm g_\alpha$ par $\rho(t)(x_\alpha)=t^{{\rm
ht}(\alpha)}x_\alpha$ pour tous $x_\alpha\in\mathbbm g_\alpha$. En
particulier $\rho(t){\bf x_+}=t {\bf x_+}$ de sorte que
$\rho(\GG_m)$ laisse stable le sous-espace vectoriel $\mathbbm
g^{\bf x_+}$ de $\mathbbm g$. Du fait que $\rho(t){\bf
x_-}=t^{-1}{\bf x_-}$, il ne laisse pas stable la section de
Kostant ${\bf x_-}+\mathbbm g^{\bf x_+}$. Considérons l'action
$\rho_+$ de $\GG_m$ sur $\mathbbm g$ définit par
$$\rho_+(t)(x)=t\rho(t)(x)$$
le composé de $\rho$ par l'homothétie. Il est clair que
$$\rho_+(t)({\bf x_-})={\bf x_-}$$
et que $\rho_+(\GG_m)$ laisse stable ${\bf x_-}+\mathbbm g^{\bf
x_+}$.

\begin{proposition}
L'isomorphisme $\epsilon^{-1}:{\bf x}_-+\mathbbm g^{\bf x_+} \rta
\car_\GG$ est $\GG_m$-équi\-va\-riante pour l'action de
$\rho_+(\GG_m)$ sur ${\bf x}_-+\mathbbm g^{\bf x_+}$ et pour
l'action de $\GG_m$ sur $\car_\GG$ mentionné dans le théorème.
\end{proposition}

\dem L'action de $\GG_m$ à travers $\rho$ sur $\mathbbm g$ laisse
stable $\mathbbm g^{\bf x_+}$ et est diagonalisable. Soit $\{{\bf
v}_1,\ldots,{\bf v}_r\}$ une base de $\mathbbm g^{\bf x_+}$ formée
de vecteurs propres sous l'action $\rho(\GG_m)$ de valeurs propres
$\rho(t){\bf v}_i=t^{n_i}{\bf v_i}$ avec $n_i\in\NN$. En ajoutant
l'homothétie, on a $\rho_+(t)=t^{n_i+1}(\bf v_i)$.

Soit $x\in \mathbbm g$ régulier. Il existe un unique $y= {\bf
x_-}+\sum_{i=1}^r v^*_i(y) {\bf v}_i$ conjugué à $x$. Les
fonctions $x\mapsto v^*_i(y)$ engendrent alors l'algèbre des
polynômes $\ad(\GG)$-invariantes sur $\mathbbm g$. Il reste à
démontrer qu'ils sont homogènes de bon degré.

Pour tout $t\in\GG_m$, $tx$ est conjugué à $ty=t{\bf
x_-}+\sum_{i=1}^r t v_i^*(y) {\bf v}_i$. Cet élément est conjugué
par $\ad\circ \rho (t)$ à
$$\ad\circ \rho (t)(ty)= {\bf x}_- +\sum_{i=1}^r t^{n_i+1}
v_i^*(y) {\bf v}_i$$ puisque pour toute racine simple
$\alpha\in\Delta$, $\la -\alpha,\rho^\vee\ra=-1$, si bien que les
polynômes $u_i$ sur $\mathbbm g$ qui correspondent à $v^*_i$ via
l'isomorphisme de Kostant $k[\mathbbm g]^{\ad(G)}\isom k[\bf
x_-+\mathbbm g^{x_+}]$, sont bien des polynômes homogènes de degré
$n_i+1$. En particulier, l'ensemble des $n_i$ et l'ensemble des
$m_i$ sont identiques,  ce qui est d'ailleurs un théorème de
Kostant. \findem

\begin{proposition} Le morphisme de Chevalley $\chi:\mathbbm g\rta\car_\GG$
est $\Aut(\GG)$-équivariant pour l'action de $\Aut(\GG)$ sur
$\mathbbm g$ induite de l'action tautologique sur $\GG$ et
l'action de $\Aut(\GG)$ sur $\car_G$ à travers $\Out(\GG)$. La
section de Kostant $\epsilon:\car_G\rta {\bf x_-}+\mathbbm g^{\bf
x_+}\subset \mathbbm g$ est $\Out(\GG)$-équivariante pour l'action
de $\Out(\GG)$ sur $\mathbbm g$ qui fixe l'épinglage.
\end{proposition}

Nous allons maintenant rappeler la construction du morphisme de
Hitchin dans un langage un peu différent de \cite{H}. Soit
${\mathbb S}$ un schéma de base quelconque. Pour tout ${\mathbb
S}$-groupe algébrique lisse $G$, on note ${\bfB}G$ le classifiant
de $G$ au-dessus de $\mathbb S$. Ce champ associe à tout $\mathbb
S$-schéma $S$ la catégorie en groupoïdes des $G$-torseurs sur $S$.
Pour tout $G$-torseur $E$ sur $S$, on note $h_E: S\rta {\bfB}G$ le
morphisme associé. Pour tout $\mathbb S$-schéma $V$ muni d'une
action relative de de $G$, on note $[V/G]$ le champ quotient de
$V$ par $G$. Par définition, il associe à tout $\mathbb S$-schéma
$S$ la catégorie dont les objets consistent en un $G$-torseur $E$
sur $S$ plus une section $\varphi:S\rta V_S\times^G E$ de l'image
réciproque de $V$ tordue par le torseur $E$. Cette construction
s'applique en particulier au $\mathbb S$-fibré vectoriel $V$ muni
d'une action $\nu:G\rta\GL(V)$. Pour tout $G$-torseur $E$ sur une
base $S$, on note $\nu(E)=V\times^G E$ le fibré vectoriel tordu
par $\nu$ sur $S$, de fibres isomorphes à $V$. Le quotient $[V/G]$
classifie les $G$-torseurs $E$ sur $X$ munis d'une section
$\varphi\in\rmH^0(S,\nu(E))$.

Maintenant $\mathbb S$ sera la courbe $X$ et $G$ sera le schéma en
groupes réductifs sur $X$ comme dans le paragraphe précédent. Le
groupe $G$ agit sur le fibré vectoriel $\Lie(G)$ par l'action
adjointe. Par ailleurs $\GG_m$ agit sur $\Lie(G)$ par homothétie
et ses deux actions commutent.

Soit $(E,\varphi)\in \calM(S)$ un point de $\calM$ à valeurs dans
un $k$-schéma $S$. La donnée de $E$ induit un morphisme $h_E:
X\times S\rta {\bfB}G$. Le fibré en droites $\calO_X(D)$ associé
au diviseur $D$ de $X$, induit un morphisme $h_D:X\rta
{\bfB}\GG_m$. On a donc un morphisme $h_E\times h_D :X\times S\rta
{\bfB}G\times {\bfB}\GG_m$. La donnée de la section $\varphi\in
\rmH^0(X\times S, \ad(E) \otimes \calO_X(D))$ est donc équivalente
à la donnée d'un morphisme
$$h_{E,\varphi}:X\times S\rta [\Lie(G)/(G\times \GG_m)]$$
qui relève $h_E\times h_D:X\times S\rta [\bfB G\times\bfB\GG_m]$.
La donnée des tels morphismes $h_{E,\varphi}$ détermine donc le
couple $(E,\varphi)$ et vice versa.

Soit $\tau_G$ le ${\rm Aut}(\GG)$-torseur associé à $G$ comme dans
§1. L'algèbre de Lie de $G$ s'obtient alors en tordant $\mathbbm
g$ par $\tau_G$
$$\Lie(G)=\mathbbm g\times^{\Aut(\GG)}\tau_G$$
de même que $G=\GG\times^{{\Aut(\GG)}}\tau_G$. On a une action
stricte de $\Aut(\GG)$ sur $[\mathbbm g/\GG\times\GG_m]$ de sorte
qu'on peut également tordre $[\mathbbm g/\GG\times\GG_m]$ par le
$\Aut(\GG)$-torseur $\tau_G$. On obtient ainsi
$$[\mathbbm g/\GG\times\GG_m]\times^{\Aut(\GG)}\tau_G
=[\Lie(G)/G\times\GG_m].$$

On a une action de ${\rm Aut}(\GG)$ à travers ${\rm Out}(\GG)$ sur
$\car_\GG$ de sorte qu'on peut tordre $\car_\GG$ par $\tau_G$
$$\car:=\car_\GG\times^{\Aut(\GG)}\tau_G.$$
Le morphisme de Chevalley $\chi:\mathbbm g\rta \car_\GG$ étant
$\Aut(\GG)\times\GG_m$-équivariant, il induit un morphisme de
$X$-schémas $\Lie(G)\rta \car$ qui est $G\times
\GG_m$-équivariant. On en déduit un morphisme de champs quotients
$$[\chi]:[\Lie(G)/(G\times\GG_m)]\rta \car/\GG_m.$$

Pour définir le morphisme de Hitchin, il suffit maintenant de
composer avec $[\chi]$ comme dans le diagramme
$$\xymatrix{
{X\times S} \ar[r]^{h_{E,\varphi}\hbox{\ \ \ \ \ \ }}
\ar@/^2pc/[rr]^{h_a} \ar[dr]_{h_E\times h_D} \ar[drr]^{\hbox{\ \ \
\ \ }h_D} & {[\Lie(G)/( G\times\GG_m)]} \ar[r]^{\hbox{\ \ \ \ \
}[\chi]} \ar[d]
&{[{\rm car}/\GG_m]} \ar[d]\\
&{{\bfB}G\times\bfB\GG_m} \ar[r] & {{\bfB}\GG_m}}
$$
Il associe donc une flèche $h_{E,\phi}$ au-dessus de $h_E\times
h_D$ une flèche $h_a$ au-dessus de $h_D$.

\begin{lemme}
Le foncteur qui associe à un $k$-schéma $S$ la catégorie des
flèches $S\times X\rta[\car/\GG_m]$ au-dessus de la flèche
$h_D:X\times S\rta\bfB\GG_m$ est représentable par un espace
affine $\bbA$ appelé l'espace affine de Hitchin.
\end{lemme}

\dem La donnée de $h_a$ est équivalente à la donnée d'une section
$$a:X\times S\rta \car\times^{\GG_m} L_D$$
de l'espace affine $\car$ tordu par le $\GG_m$-torseur $L_D$
au-dessus de $X\times S$. Le tordu $\car\times^{\GG_m}L_D$ est
l'espace total d'un fibré vectoriel sur $X$ car il provient après
torsion de l'espace affine $\car_\GG$. Par suite, le foncteur
ci-dessus mentionné est représentable par l'espace affine associé
à l'espace des sections globales de ce fibré vectoriel. \findem

\bigskip
Dans le cas où $G$ est isomorphe au-groupe constant $G\isom
X\times\GG$, $\car\times^{\GG_m}L_D$ est l'espace total du fibré
totalement décomposé $\bigoplus_{i=1}^r \calO_X((m_i+1)D)$.
L'espace de ses sections globales est
$$\bbA(k)=\bigoplus_{i=1}^r \rmH^0(X,\calO_X((m_i+1)D)).$$
C'est la description originale de Hitchin de cet espace affine.
Dans le cas où $G$ est le groupe unitaire associé à un revêtement
étale quadratique $X'\rta X$, le fibré inversible de carré trivial
sur $X$ associé à $X'\rta X$, intervient dans la description de
l'espace de Hitchin, voir \cite{L-N}.

Il sera bien commode de pouvoir choisir un $k$-point dans la fibre
$\calM_a$. Ceci nécessite des hypothèses supplémentaires sans
lesquelles $\calM_a(k)$ pourrait être vide.

\begin{proposition} On supposera que $G$ est
quasi-déployé c'est-à-dire que les ${\rm Aut}(\GG)$-torseurs
$\tau_G$ et $\tau_G^{\rm ext}$ construits dans §1 sont isomorphes.
On se donne de plus une racine carrée $L_D^{\otimes 1/2}$ de
$L_D$. Alors, il existe une section du morphisme de Hitchin.
\end{proposition}

\dem On va construire une section du morphisme de Hitchin en
modelant sur la section de Kostant.

On a vu que la section de Kostant $\car_\GG\isom {\bf
x_-}+\frakg^{\bf x_+}$ induit un isomorphisme
$$[\car_\GG/\GG_m]\isom [({\bf x_-}+\frakg^{\bf
x_+})/\rho_+(\GG_m)].$$ En considérant l'homomorphisme diagonal
$\GG_m\rta \GG_m\times\GG_m$, on obtient un morphisme
$$[({\bf x_-}+\frakg^{\bf
x_+})/\rho_+(\GG_m)]\rta  [\frakg/\GG_m\times\rho(\GG_m)]$$ où
dans champ but, le premier facteur de $\GG_m$ agit par
l'homothétie et où le second par $\rho$. La petite complication
vient du fait que l'homomorphisme $\rho:\GG_m\rta{\rm
Aut}(\frakg)$ ne factorise par un caractère $\GG_m\rta G$ et il
faut en prendre une racine carrée.

Notons donc $\GG_m^{[2]}\rta \GG_m$ l'homomorphisme d'une copie de
$\GG_m$ notée $\GG_m^{[2]}$ dans $\GG_m$ qui consiste à
l'élévation à la puissance $2$. Par le changement de base
$\bfB\GG_m^{[2]}\rta\bfB\GG_m$, on obtient les morphismes
$$[\car_\GG/\GG_m^{[2]}]\isom [({\bf x_-}+\frakg^{\bf
x_+})/\rho_+(\GG_m^{[2]})] \rta
[\frakg/\GG_m^{[2]}\times\rho(\GG_m^{[2]})]$$ où $\GG_m^{[2]}$
agit à travers l'homomorphisme $\GG_m^{[2]}\rta \GG_m$.

On a maintenant un cocaractère $2\rho:\GG_m^{[2]} \rta {\mathbbm
T}\rta \GG$ qui induit un morphisme
$$
[\frakg/\GG_m^{[2]}\times\rho(\GG_m^{[2]})] \rta
[\frakg/\GG_m^{[2]}\times \GG].
$$
En composant, on obtient donc un morphisme
$$[\car_\GG/\GG_m^{[2]}] \rta [\frakg/\GG_m^{[2]}\times \GG].
$$
qui relève le morphisme de Chevalley $[\frakg/\GG_m^{[2]}\times
\GG] \rta [\car_\GG/\GG_m^{[2]}]$. Le morphisme ci-dessus est
équivariant pour l'action de ${\rm Out}(\GG)$ de sorte qu'on peut
le tordre par le torseur $\tau^{\rm ext}(\GG)$. En utilisant
l'hypothèse $G$ quasi-déployé $\tau(G)=\tau^{\rm ext}(G)$, on en
déduit un morphisme
$$[\car/\GG_m^{[2]}] \rta [\Lie(G)/G\times \GG_m^{[2]}].
$$

Maintenant, la donnée d'une flèche $h_a:X\times S\rta
[\car/\GG_m]$ au-dessus de $h_D:X\rta \bfB\GG_m$ correspondant au
$\GG_m$-torseur $L_D$, plus une racine carrée $L_D^{\otimes 1/2}$
de $L_D$, correspond simplement à une flèche $h_a^{1/2} :X\times
S\rta [\car/\GG_m^{[2]}]$. En composant avec la section définie
ci-dessus on obtient une flèche
$$h_\phi^{1/2}:X\times S\rta [\Lie(G)/G\times \GG_m^{[2]}]$$
qui consiste en une flèche
$$h_\phi:X\times S\rta [\Lie(G)/G\times \GG_m]$$
au-dessus de $h_D:X\times S\rta \bfB\GG_m$, plus d'une racine
carrée de $L_D$. \findem

\section{Centralisateurs}

Dans cette section, on va construire un schéma en groupes lisse
$J$ au-dessus de $\car$ muni d'un homomorphisme $\chi^*J\rta I$
vers le schéma en groupes des centralisateurs $I$ au-dessus de
$\frak g$. Cet homomorphisme est un isomorphisme au-dessus de
l'ouvert $\frakg^{\rm reg}$ des éléments réguliers de $\frakg$.

 Considérons le schéma en
groupes $I$ sur $\mathbbm g$ des paires
$$I_\GG=\{(x,g)\in \mathbbm g \times \GG\mid \ad(g)x=x\}.$$
L'action adjointe de $G$ sur $\mathbbm g$ se relève à $I$ de façon
évidente $\ad(h)(x,g)=(\ad(h) x, hgh^{-1})$. L'action par
homothétie de $\GG_m$ sur $\mathbbm g$ se relève également
$t(x,g)=(tx,g)$. On en déduit un schéma en groupes
$$[I_\GG]:=[I_\GG/\GG\times \GG_m]$$
au-dessus de $[\mathbbm g/G\times \GG_m]$. De plus, l'action de
$\Aut(\GG)$ sur $[\mathbbm g/\GG\times \GG_m]$ relève à $[I_\GG]$
de sorte qu'on peut tordre $[I_\GG]$ par le $\Aut(\GG)$ torseur
$\tau_G$. On obtient ainsi un schéma en groupes $[I]$ au-dessus de
$[\frakg/G\times\GG_m]$. Notons un énoncé tautologique.

\begin{lemme}
Soit $(E,\varphi)\in\calM(S)$ une paire de Hitchin. Soit
$h_{E,\varphi}:X\times S\rta  [\frakg/G\times\GG_m]$ le morphisme
associé. Alors on a
$$I_{E,\varphi}=h_{E,\varphi}^* [I].$$
\end{lemme}

Le morphisme caractéristique de Chevalley $\chi:\mathbbm g\rta
\car_\GG$ admet une section, due à Kostant $\epsilon:\car_\GG\rta
\mathbbm g$ dont l'image est contenue dans l'ouvert dense
$\mathbbm g^{\rm reg}$ des éléments réguliers de $\mathbbm g$. On
note $J_\GG=\epsilon^* I_\GG$ et on l'appelle le {\em
centralisateur régulier}. La façon dont on construit ici $J_\GG$
dépend de la section de Kostant mais il n'en dépend pas en fait.

Par définition, $\mathbbm g^{\rm reg}$ est l'ouvert de $\mathbbm
g$ des éléments $x$ dont le centralisateur $I_x$ a dimension égale
à $r$ le rang de $\mathbbm g$. On sait de plus qu'au-dessus de
$\mathbbm g^{\rm reg}$, le schéma en groupes $I_\GG\rta \mathbbm
g$ est lisse. Puisque l'image de la section de Kostant $\epsilon$
est contenu dans $\mathbbm g^{\rm reg}$, le schéma en groupes
$J_\GG\rta \car_\GG$ est aussi lisse de dimension relative $r$.
Rappelons qu'ici on a utilisé l'hypothèse que la caractéristique
de $k$ ne divise pas l'ordre de groupe de Weyl de $\mathbbm G$. Il
sera intéressant de savoir si en petites caractéristiques, {\em la
restriction de $I_{\mathbbm G}$ à $\frakg^{\rm reg}$ est encore
plate}. Si c'est le cas, la plupart des résultats de ce travail
devrait se généraliser en petites caractéristiques.

\begin{proposition}
On a un isomorphisme canonique $\chi^* J_\GG|_{\mathbbm g^{\rm
reg}} \isom I_\GG|_{\mathbbm g^{\rm reg}}$ qui, de plus, se
prolonge en un homomorphisme de schéma en groupes $\chi^*
J_\GG\rta I_\GG$.
\end{proposition}

\dem Considérons le morphisme
$$\beta: \GG\times \car_\GG \rta \mathbbm g^{\rm reg}$$
défini par $(g,a)\mapsto \ad(g)\epsilon(a)$ où
$\epsilon:\car_\GG\rta\mathbbm g$ est la section de Kostant. C'est
un morphisme lisse et surjectif donc fidèlement plat. L'image
inverse de $I_\GG|_{\mathbbm g^{\rm reg}}$ par $\beta$ est
l'ensemble des triplets $(g,a;h)$ tel que
$\ad(hg)\epsilon(a)=\ad(g)\epsilon(a)$. L'image inverse de $\chi^*
J_\GG|_{\mathbbm g^{\rm reg}}$ sur $\GG\times\car_\GG$ est
l'ensemble des triplets $(g,a; h')$ tel que $\ad(h')\kappa(a)
=\kappa(a)$. On a clairement un isomorphisme
$$\beta^*  \chi^* J|_{\mathbbm g^{\rm reg}} \isom
\beta^* I|_{\mathbbm g^{\rm reg}}$$ défini par $(g,a,h')\mapsto
(g,a;h)$ avec $h=gh' g^{-1}$.

En vertu du théorème de descente fidèlement plate, pour démontrer
que cet isomorphisme descend en un isomorphisme $ \chi^*
J_\GG|_{\mathbbm g^{\rm reg}} \isom I_\GG|_{\mathbbm g^{\rm
reg}}$, il suffit de démontrer une égalité de cocycle sur
$$(G\times \car_\GG) \times_{\mathbbm g^{\rm reg}} (G\times \car_\GG).$$
Ce produit fibré est l'espace des $(g_1,g_2;a)$ tels que
$\ad(g_1)\epsilon(a)=\ad(g_2)\epsilon(a)$, ou autrement dit
$g_1^{-1}g_2 \in I_{\kappa(a)}$. L'égalité de cocycle à démontrer
consiste entre l'égalité entre les deux isomorphismes entre les
images inverses de $\chi^* J_\GG|_{\mathbbm g^{\rm reg}}$ et
$I_\GG|_{\mathbbm g^{\rm reg}}$ donné par $h'\mapsto g_1h'
g_1^{-1}$ et $h'\mapsto g_1h' g_1^{-1}$. Ils se diffèrent l'un de
l'autre de l'automorphisme intérieur ${\rm int}(g^{-1}g_2)$ de
$I_{\kappa(a)}$ lequel n'est autre l'identité parce que
$I_{\kappa(a)}$ est un groupe {\em commutatif}. On en déduit par
le théorème de descente fidèlement plate de Grothendieck, un
isomorphisme $\chi^* J_\GG|_{\mathbbm g^{\rm reg}} \isom
I_\GG|_{\mathbbm g^{\rm reg}}$.

Puisque $J_\GG$ est lisse sur $\car_\GG$, $\chi^* J_\GG$ est un
schéma en groupes lisse sur $\mathbbm g$. Puisque $\mathbbm
g-\mathbbm g^{\rm reg}$ est un fermé de codimension trois de
$\mathbbm g$, en particulier plus grande que deux,
$\chi^*J-\chi^*J|_{\mathbbm g^{\rm reg}}$ est un fermé de
codimension plus grand que 2 du schéma $\chi^* J_\GG$. Celui-ci
est par ailleurs un schéma lisse donc normal. Le morphisme $\chi^*
J_\GG|_{\mathbbm g^{\rm reg}} \rta I_\GG$ à valeur dans le schéma
affine $I_\GG$ se prolonge donc en un unique morphisme $\chi^*
J_\GG\rta I_\GG$ [EGA IV.4 20.4.12]. L'unicité du prolongement
implique que celui-ci est un homomorphisme de schémas en groupes
sur $\mathbbm g$. \hfill $\square$

\bigskip
Cette proposition généralise une opération bien connue dans le cas
$\GL(n)$. L'espace $\car_\GG$ classifie dans ce cas les polynômes
unitaires de degré $n$. Pour tout $a\in \car_\GG$ vu comme un
polynôme unitaire $a\in k[t]$, $J_a$ est le groupe des éléments
inversibles de l'algèbre $k[t]/(a)$. Pour toute application
linéaire $x:k^n \rta k^n$ de polynôme caractéristique $a$, on peut
faire agir l'algèbre $k[t]/(a)$ sur l'espace vectoriel $k^n$ en
faisant agir la variable $t$ comme l'endomorphisme $x$. Cette
action commute à $x$ si bien qu'elle définit un homomorphisme $J_a
\rta I_x$ où $I_x$ est le centralisateur de $x$.

\begin{proposition}
Il existe un schéma en groupes $[J]$ au-dessus du quotient
$[\car/\GG_m]$, uniquement déterminé à un unique isomorphisme
près, dont l'image inverse sur $\car$ est le schéma des
centralisateurs réguliers $J$. De plus, sur $[\frakg/G\times
\GG_m]$, on a un homomorphisme canonique $[\chi]^* [J] \rta [I]$
dont la restriction à $[\frakg^{\rm reg}/G\times\GG_m]$ est un
isomorphisme.
\end{proposition}

\dem L'action de $\GG_m$ sur $\mathbbm g$ par l'homothétie
préserve l'ouvert $\mathbbm g^{\rm reg}$ des éléments réguliers,
d'où une action de $\GG_m$ sur $I|_{\mathbbm g^{\rm reg}}$ et donc
une action de $\GG_m$ sur $\chi^* J|_{\mathbbm g^{\rm reg}}$,
compatible à l'action de $\GG_m$ sur $\mathbbm g^{\rm reg}$. Cette
action induit par la descente fidèlement plate une action de
$\GG_m$ sur $J$ compatible à l'action $\nu_\chi(\GG_m)$ sur
$\car$. En prenant le quotient au sens des champs de $J$ par cette
action de $\GG_m$, on obtient le schéma en groupes $[J_\GG]$ sur
$[\car/\nu_\chi(\GG_m]$ désiré. L'homomorphisme $[\chi]^* [J_\GG]
\rta [I_\GG]$ s'obtient également de $\chi^* J_\GG\rta I_\GG$ par
la descente fidèlement plate.

Pour passer du groupe constant $\GG$ au schéma en groupes $G$, on
utilise de nouveau le $\Aut(\GG)$-torseur $\tau_G$. En tordant
$[J_\GG]/[\car/\GG_m]$ par $\tau_G$, on obtient un schéma en
groupes lisses $[J]/[\car/\GG_m]$ qui vient avec un homomorphisme
$\chi^*[J]\rta [I]$ entre schémas en groupes sur
$[\frakg/G\times\GG_m]$. \findem

\begin{proposition}
Le morphisme $[\mathbbm g^{\rm reg}/\GG] \rta \car_\GG$ est une
gerbe liée par $J$. De plus, cette gerbe est neutralisée par la
section de Kostant et on a un isomorphisme $\Out(\GG)$-équivariant
$$[\mathbbm g^{\rm reg}/\GG] \isom [\car_\GG/J].$$
\end{proposition}

\dem Il revient au même de démontrer l'isomorphisme
$$[\mathbbm g^{\rm reg}/\GG] \isom [\car_\GG/J].$$
Pour cela, on utilise la section de Kostant $\kappa:\car_\GG\rta
\mathbbm g^{\rm reg}$ et le morphisme $\GG\times \car_\GG \rta
\mathbbm g^{\rm reg}$ donné par $(g,a)\mapsto \ad(g) \kappa(a)$.
Ce morphisme est lisse et surjectif. Puisque $J=\kappa^* I$ où $I$
est le schéma en groupes des centralisateurs, on a un isomorphisme
$$(\GG\times \car_\GG)/J \isom \mathbbm g^{\rm reg}$$
d'où on déduit l'isomorphisme $[\car_\GG/J]\isom [\mathbbm g^{\rm
reg}/\GG]$ en divisant par l'action de $\GG$. Puisque la section
de Kostant est $\Out(\GG)$-équivariant, la neutralisation de la
$J$-gerbe $[\mathbbm g^{\rm reg}/\GG]$ est
$\Out(\GG)$-équivariante.  \hfill $\square$

\bigskip Notons que la neutralisation définie par la section de
Kostant est $\Out(G)$-équivariante mais n'est pas
$\Aut(\GG)$-équivariante pour l'action de $\Aut(\GG)$ sur $\frakg$
déduite de l'action tautologique de $\Aut(\GG)$ sur $\GG$. La
structure de $[\mathbbm g^{\rm reg}/\GG]$ comme $J$-gerbe est
quant à elle $\Aut(\GG)$-équivariante. On en déduit le corollaire
suivant après la torsion par le $\Aut(\GG)$-torseur $\tau_G$.

\begin{proposition}
Au-dessus de $X$, $[\frakg^{\rm reg}/G]$ est une $J$-gerbe
au-dessus de $\car$. Cette gerbe est neutre si $G$ est
quasi-déployée c'est-à-dire si $\tau_G=\tau^{\rm ext}_G$.
\end{proposition}

Le centralisateur régulier $J_\GG$ au-dessus de $\car_\GG$ peut
être décrit d'une autre façon d'après un résultat de Donagi et
Gaitsgory que nous allons passer en revue (\cf \cite{DG}).
Consi\-dérons le revêtement fini galoisien $\mathbbm
t\rta\car_\GG$ de groupe de Galois $W$, ramifié le long du
diviseur
$$D=\bigcup_{\alpha\in\Phi}D_\alpha$$
où $D_\alpha$ est l'hyperplan de $D$ est défini par la dérivée
$d\alpha:\mathbbm t\rta\GG_a$, et est aussi le lieu des points
fixes de l'involution $\mathbf s_\alpha\in W$ associée à la racine
$\alpha$. Considérons le faisceau en groupes $\calJ^+$ sur $\car$
qui associe à tout $\car_\GG$-schéma $S$, le groupe des morphismes
$$t:\widetilde S=S\times_{\car_\GG} \mathbbm t \rta \mathbbm T$$
qui sont $W$-équivariants. Soit $\widetilde s \in \widetilde S$ un
point fixe de $\mathbf s_\alpha$, alors $t(\widetilde s)$ est un
point fixe par $\mathbf s_\alpha$ dans $\mathbbm T$. On en déduit
que $\alpha(t(\widetilde s))=\pm 1$.

\begin{definition} Soit $\calJ$ le sous-foncteur ouvert de
$\calJ^+$ qui associe à tout $\car_\GG$-schéma $S$ le groupe des
morphismes $t:\widetilde S\rta \mathbbm T$ qui sont
$W$-équivariants et tels que pour tout $\widetilde s\in \widetilde
S$ fixe par l'involution ${\mathbf s}_\alpha$, on a
$\alpha(t(\widetilde s))\not=-1$.
\end{definition}

La condition subtile $\alpha(t(x))\not=-1$ pour les $x$ tels que
$s_\alpha(x)=x$, fait la distinction entre $\SL_2$ et ${\rm
PGL}_2$. Dans ces deux cas, le tore maximal $\mathbbm T$ est
isomorphe à $\GG_m$, le groupe de Weyl $W=\frakS_2$ agit sur
$\mathbbm T$ par $t\mapsto t^{-1}$ et le sous-groupe des points
fixes de $\mathbbm T$ sous $W$ est $\{\pm 1\}$. Dans le cas
$\SL_2$, la racine positive consiste en le caractère
$\alpha(t)=t^2$ ; dans le cas ${\rm PGL}_2$, c'est le caractère
$\alpha(t)=t$. Ainsi, la condition $\alpha(t(x))\not=-1$ est vide
dans le cas $\SL_2$, mais elle est non vide dans le cas $\PGL_2$.

L'énoncé suivant est du à Donagi et Gaitsgory (\cf \cite{DG}
théorème 11.6).

\begin{proposition} On a un isomorphisme $J_\GG\isom \calJ$ entre
faisceaux en groupes au-dessus de $\car_\GG$.
\end{proposition}

Le théorème 11.6. de loc. cit s'énonce dans un contexte légèrement
diffé\-rent qu'ici. En particulier, leur centralisateur régulier
est un schéma en groupes lisse sur une base différente de
$\car_\GG$. Toutefois, on peut montrer facilement que les deux
énoncés sont équi\-va\-lents en utilisant la résolution simultanée
de Grothendieck qui relie les deux bases.

Après torsion par le $\Out(\GG)$-torseur $\tau^\Out_G$ au-dessus
de $X$, on obtient l'énoncé analogue pour le schéma en groupes $J$
au-dessus de $\car$. On a en effet
$\car=\car_\GG\times^{\Out(\GG)}\tau^{\Out}_G$ et
$J=J_\GG\times^{\Out(\GG)}\tau^\Out_G$.

Rappelons qu'on a un revêtement fini étale $X_\Theta\rta X$,
pointé par un point géométrique $x_\Theta$, de groupe de Galois
$\Theta$ où $\Theta$ est l'image de l'homomorphisme
$\rho_G:\pi_1(X,x)\rta \Out(\GG)$ associé au $\Out(\GG)$-torseur
$\tau^{\Out}_G$. On a un morphisme fini ramifié
$$X_\Theta\times \mathbbm t\rta \car$$
dont la restriction à la partie semi-simple régulière
$X_\Theta\times \mathbbm t^{\rm ssr}\rta \car^{\rm ssr}$ est étale
galoisien de groupe de Galois $W'=W\rtimes\Theta$.

La proposition suivante est une conséquence immédiate de 3.7.

\begin{proposition} Comme faisceau sur la topologie étale de $\car$,
$J$ associe à tout $\car$-schéma $S$ le groupe des flèches
$W'$-équivariantes
$$t:\widetilde S= (X_\Theta\times \mathbbm t)\times_\car S \rta
\mathbbm T
$$
telles que pour toute racine $\alpha$ et pour tout $\widetilde
s\in\widetilde S$ fixe par l'involution ${\mathbf s}_\alpha$, on a
$\alpha(t(\widetilde s))\not=-1$.
\end{proposition}

\section{Action de $P_a$ sur $\calM_a$}

Dans cette section, on construit l'action canonique du champ de
Picard relatif $P$ au-dessus de $\bbA$ sur la fibration de Hitchin
$f:\calM\rta \bbA$. On démontre ensuite que l'action de $P_a$ sur
$\calM_a=f^{-1}(a)$ est simplement transitive pour des
caractéristiques $a$ très régulières. Plus généralement, pour des
caractéristiques génériquement semi-simples régulières, on
démontre une formule de produit pour le quotient $[\calM_a/P_a]$.

Soit $a:S\rta \bbA$ un $k$-schéma au-dessus de $\bbA$. La donnée
de $a$ est équivalente à la donnée d'une flèche $h_a:X\times S\rta
[\car/\GG_m]$ qui relève $h_{D}:X\rta \bfB\GG_m$. En prenant
l'image réciproque de $[J]$, on a un schéma en groupes lisse $J
_a:=h_a^* [J]$ au-dessus de $X\times S$. Notons $P_a$ {\em la
catégorie de Picard des $J_a$-torseurs sur $X\times S$}. La règle
$a\mapsto P_a$ définit un champ de Picard au-dessus de l'espace
affine de Hitchin $\bbA$ qu'on note $P$.

On va maintenant construire une action de $P$ sur le champ $\calM$
relativement au-dessus de $\bbA$. Soit $S$ un $\FF_q$-schéma et
$a\in\bbA(S)$. La fibre $f^{-1}(a)$ est la catégorie des paires
$(E,\varphi)\in \calM(S)$ ayant la caractéristique $a$. Il s'agit
de faire agir la catégorie de Picard $P_a$ sur cette catégorie
fibre $\calM_a$.

Pour tout objet $(E,\varphi)$ de $f^{-1}(a)$, on a le schéma en
groupes $I_{E,\varphi}$ qui représente le faisceau des
automorphismes de $(E,\varphi)$ au-dessus de $X\times S$. On a un
homomorphisme de schémas en groupes
$$J_a\rta I_{(E,\varphi)}$$
qui se déduit de l'homomorphisme canonique $\chi^*J\rta I$ sur
$\frakg$. Il s'ensuit qu'on peut tordre la paire $(E,\varphi)$ à
l'aide de tout torseur sous $J_a$ sans changer la caractéristique
de $(E,\varphi)$. Ceci définit {\em l'action de $P_a$ sur
$\calM_a$}.

\medskip
D'après le théorème de restriction de Chevalley-Kostant, rappelé
en §2, on sait que l'anneau des fonctions algébriques sur l'espace
affine $\car_\GG$ est canoniquement isomorphe à l'anneau
$k[\mathbbm t]^W$ des fonctions sur $\mathbbm t$ qui sont
$W$-invariantes. On a donc un morphisme $\pi:\mathbbm t\rta
\car_\GG$ qui est fini, génériquement étale galoisien de groupe de
Galois $W$. Soit $\FrakB_\GG$ le lieu de branchement de $\pi$.
C'est le diviseur de $\car_\GG$ défini par l'annulation de la
fonction discriminante $\prod_{\alpha\in\Phi}d\alpha$. Ici, les
fonctions $d\alpha:\mathbbm t\rta\GG_a$ désignent la dérivation de
la racine $\alpha:T\rta\GG_m$ ; leur produit étant visiblement
invariant par $W$. Après torsion par le $\Out(\GG)$-torseur
$\tau^\Out_G$, on obtient au-dessus de $X$, un morphisme fini
ramifié $\mathbbm t\times^{\Out(\GG)}\tau^\Out_G \rta
\car_\GG\times^{\Out(\GG)}\tau^\Out_G$ que nous allons noter
$\frakt\rta\car$. On notera $\FrakB\subset\car$ le lieu de
branchement de $\frakt\rta\car$ qui se déduit de $\FrakB_\GG$ par
torsion par $\tau^\Out_G$.

Une caractéristique $a\in\bbA(\overline k)$ est une section
$$h_a:\overline X\rta \car\times^{\GG_m}L_D.
$$
Elle est dite {\em très régulière} si $h_a(X)$ coupe
transversalement le lieu de branchement $\FrakB\times^{\GG_m}L_D$
du morphisme fini qui se déduit de $\pi$ en tordant par le
$\GG_m$-torseur $L_D$

\begin{definition}
La caractéristique $a$ est dite {\bf très régulière} si l'image de
la section $s_a(X)$ coupe transversalement le diviseur
$\FrakB\times^{\GG_m} L_D$, c'est-à-dire qu'elle coupe ce diviseur
dans sa partie lisse avec la multiplicité un en chaque point de
l'intersection.
\end{definition}

Lorsque le diviseur $D$ est très ample, les caractéristiques très
régulières forment un ouvert dense de $\bbA$ d'après le théorème
de Bertini (\cf \cite{Fa}).

\begin{proposition}
Soit $a\in\bbA(\bar k)$ une caractéristique très régulière. Soit
$(E,\varphi)$ une paire de Hitchin de caractéristique $a$. Alors
l'application induite
$$h_{E,\varphi}:X\rta [\frakg/G\times \GG_m]$$
se factorise par l'ouvert $[\frakg^{\rm reg}/G\times \GG_m]$.
\end{proposition}

\dem Le problème étant local, on peut remplacer $X$ par
$\Spec(k[[t]])$, $E$ par un $G$-torseur trivial et $L_D$ par un
$\GG_m$-torseur trivial. La section $\varphi$ est alors une flèche
$$\varphi:\Spec(\bar k[[t]]) \rta \Spec(\calO_{\frakg,x})$$
où $\calO_{\frakg,x}$ est le localisé de $\frakg$ en un point
$x\in\frakg$. Notons $\varphi^\# :  \calO_{\frakg,x} \rta k[[t]]$
l'homomorphisme local associé. Notons ${\frak m}_{x}$ l'idéal
maximal de $\calO_{\frakg,x}$ ; on a $(\varphi^\#)^{-1}(tk[[t]])
={\frak m}_x$ en vertu de la localité de $\varphi^\#$.

On doit démontrer que $x\in\frakg^{\rm reg}$.
Si ce n'est pas le cas, on a
$${\rm discr}=\prod_{\alpha\in\Phi} d\alpha \in {\frak m}_x^2$$
puisque la fonction discriminante s'annule à l'ordre au moins $2$
à un point $x\in\frakg$ non régulier. Il s'ensuit que l'image de
${\rm discr}$ dans $k[[t]]$ est de valuation au moins $2$. Ceci
contredit l'hypothèse que la caractéristique $a$ est très
régulière. \hfill $\square$

\bigskip
L'énoncé suivant est une variante d'un théorème de Faltings (\cf
\cite{Fa} III.2).

\begin{proposition}
Si $a\in \bbA(\bar k)$ est une caractéristique très régulière,
l'action du champ de Picard $P_a$ sur la fibre de Hitchin
$\calM_a$ est simplement transitive. Autrement dit, $\calM_a$ est
une gerbe liée par le champ de Picard $P_a$.

Sous l'hypothèse supplémentaire de 2.5, cette gerbe est neutre.
\end{proposition}

\dem Soit $a\in\bbA(\bar k)$ une caractéristique très régulière.
La donnée de $a$ est équivalente à la donnée d'une flèche
$h_a:\ovl X\rta [\car/\GG_m]$ au-dessus de $h_D:\ovl X\rta
\bfB\GG_m$. La donnée d'un $S$-point $(E,\varphi)$ dans la fibre
$\calM_a$, est équivalente à la donnée d'une flèche
$h_{E,\varphi}$ qui, dans le diagramme
$$\xymatrix{
{X\times S} \ar[r]^{h_{E,\varphi}\hbox{\ \ \ }} \ar[dr]_{h_a} &
{[\frakg/ (G\times\GG_m)]} \ar[d]^{[\chi]} \\
& [\car/\GG_m]}
$$
relève le morphisme $h_a$ constant sur le facteur $S$.

L'hypothèse $a$ très régulière force la flèche $h_{E,\varphi}$ à
se factoriser par l'ouvert $[\frakg^{\rm reg}/ G\times \GG_m]$.
Or, $[\frakg^{\rm reg}/( G\times \GG_m)]$ est une $J$-gerbe
au-dessus de $[\car/\GG_m]$ de sorte que la fibre $\calM_a$ est un
torseur sous la catégorie de Picard $P_a$ des $J_a$-torseurs.
\findem

\bigskip Lorsque la caractéristique $a$ n'est plus très régulière, l'action
de $P_a$ sur $\calM_a$ n'est plus simplement transitive. On peut
encore dire quelque chose sur le quotient $[\calM_a/P_a]$ si la
caractéristique de $a$ est {\em génériquement semi-simple
régulière}. Dans le cas où $G=\GL_n$ et $p>n$, cette hypothèse
signifie que la courbe spectrale est réduite.

\begin{definition}
Une caractéristique $a\in\bbA(\bar k)$ est dite génériquement
semi-simple régulière, si l'image du morphisme $h_a:X\rta
\car\times^{\GG_m}L_D$ qui lui est associé, n'est pas contenue
dans $\FrakB\times^{\GG_m}L_D$ où $\FrakB\subset \car$ est le lieu
de branchement du morphisme $\pi:\frakt \rta\car$.
\end{definition}

On note $\bbA^{\rmgssr}$ l'ouvert de $\bbA$ des $a$ génériquement
semi-simple régulières. Soit $a\in\bbA^{\rmgssr}(\bar k)$. Soit
$U_a$ l'image réciproque de $\ovl X$ du complémentaire de
$[\FrakB/\nu_\chi(\GG_m)]$. Le complémentaire $\ovl X-U_a$ est
alors un ensemble fini de points. Soit $v\in \ovl X-U_a$. Notons
$\ovl X_v$ le complété de $\ovl X$ en $v$, le disque formel autour
de $\ovl X$, et $\ovl X_v^\bullet$ le disque formel pointé.

Considérons la catégorie $\calM_{a,v}$ des paires
$(E_v,\varphi_v)$ où $E_v$ est un $G$-torseur sur $\ovl X_v$ et
$\varphi_v$ est une section de $\varphi_v\in \rmH^0(\ovl
X_v,\ad(E_v)\otimes \calO_X(D))$ ayant pour caractéristique $a$.
Autrement dit $\calM_{a,v}$ est la catégorie des flèches
$$h_{\varphi_v}:\ovl X_v\rta [\frakg/G\times \GG_m]$$
qui relève la flèche $h_{a_v}:\ovl X_v\times [\car/\GG_m]$
déduites de $a:\ovl X\rta [\car/\GG_m]$. Notons $J_{a_v}$ la
restriction de $J_a$ à $\ovl X_v$; on a $J_{a_v}=h_{a_v}^* [J]$.
Soit $P_{a,v}$ la catégorie de Picard des $J_{a,v}$-torseurs sur
$\ovl X_v$. On a une action de la catégorie de Picard $P_{a,v}$
sur la catégorie $\calM_{a,v}$ définie de la même façon que
l'action de $P_a$ sur $\calM_a$. On peut former des $2$-catégories
quotients $[\calM_{a,v}/P_{a,v}]$ analogues locaux du $2$-quotient
$[\calM_a(\ovl k)/P_a(\ovl k)]$, voir le complément qui suit pour
un petit résumé sur les $2$-catégories quotients.

\begin{lemme} Soit $a\in\bbA^\heartsuit(\ovl k)$ une caractéristique
génériquement semi-simple régulière. Les $2$-catégories quotients
$[\calM_a(\ovl k)/P_a(\ovl k)]$ et $[\calM_{a,v}/P_{a,v}]$ sont
alors équivalentes à des $1$-catégories.
\end{lemme}

\dem Il suffit de vérifier que le critère de Dat pour qu'un
$2$-quotient soit équivalent à une $1$-catégorie, \cf lemme 4.7
qui suit. Dans le cas local, il s'agit de démontrer que pout tout
objet $(E_v,\varphi_v)$ de $\calM_{a,v}$, l'homomorphisme
$${\rm Aut}_{P_{a,v}}(1_{P_{a,v}})\rta {\rm Aut}_{\calM_{a,v}}
(E_v,\varphi_v)
$$
du groupe des automorphismes de l'objet neutre de $P_{a,v}$ dans
celui de $(E_v,\varphi_v)$ qui se déduit de l'action de $P_{a,v}$
sur $\calM_{a,v}$ est un homomorphisme injectif.

Pour démontrer cette injectivité, notons que la restriction du
disque formel $\ovl X_v$ au disque formel pointé $\ovl
X^\bullet_v$ induit des homomorphismes injectifs lesquels sont
représentés par les flèches verticales dans le diagramme suivant
$$
\xymatrix{
  {\rm Aut}_{P_{a,v}}(1_{P_{a,v}}) \ar[d]_{} \ar[r]^{} &
  {\rm Aut}_{\calM_{a,v}}(E_v,\varphi_v)
   \ar[d]^{} \\
  {\rm Aut}(1_{P_{a,v}}|_{\ovl
X^\bullet_v})  \ar[r]^{} & {\rm Aut}((E_v,\varphi_v)|_{\ovl
X^\bullet_v} )  }
$$
De plus, l'hypothèse que $a$ est génériquement semi-simple
régulière implique que la flèche horizontale en bas est un
isomorphisme. Il s'ensuit que la flèche horizontale en haut est
également injective. La même démonstration vaut pour le quotient
$[\calM_a(\ovl k)/P_a(\ovl k)]$. \findem

\bigskip Supposons que la catégorie $\calM_a(\ovl k)$ est non vide.
C'est le cas notamment si les hypothèses de 2.5 soient
satisfaites. Soit $(E^*,\varphi^*)$ un point dans $\calM_a(\ovl
k)$ par exemple le point de Kostant construit dans 2.5.  Soit
$\calM^\bullet_{a,v}$ l'ensemble des
$(E_v,\varphi_v,\iota^\bullet_v)$ où $(E_v,\varphi_v)$ est un
objet de $\calM_v(a)$ et où $\iota_v^\bullet$ est un isomorphisme
au-dessus du disque pointé $\ovl X^\bullet_v$ entre
$(E_v,\varphi_v)$ et la paire de base $(E^*,\varphi^*)$.
Considérons aussi le groupe $P^\bullet_{a,v}$ des
$J_{a_v}$-torseurs sur $\ovl X_v$ munis d'une trivialisation sur
$\ovl X^\bullet_v$. On a encore une action de $P^\bullet_{a,v}$
sur $\calM^\bullet_{a,v}$.

\begin{theoreme} Supposons que la catégorie $\calM_a(\ovl k)$ est
non vide. Alors, on a des équivalences de catégories
$$[\calM_a(\ovl k)/P_a(\ovl k)]= \prod_{v\in\ovl X-\ovl U_a}
[\calM_{a,v}/P_{a,v}]=\prod_{v\in\ovl X-\ovl U_a}
[\calM^\bullet_{a,v}/P^\bullet_{a,v}].
$$
\end{theoreme}

\dem On a des morphismes naturels qui forment un triangle :
$$\xymatrix{
\prod_v[\calM^\bullet_{a,v}/P^\bullet_{a,v}]  \ar[rr]^{\alpha}
\ar[dr]_{\beta}
                &  &  [\calM_a/P_a]   \ar[dl]^{\gamma}    \\
                &
\prod_v[\calM_{a,v}/P_{a,v}]                 }
$$
où :
\begin{itemize}
\item la flèche $\alpha$ est une flèche de recollement formel à la
Beauville-Laszlo (\cf \cite{BL}).  \`A l'aide des isomorphismes
sur les disques pointés $\ovl X_v^\bullet$, elle recolle les
$G$-torseurs $E_v$ sur les disque $\ovl X_v$ à la restriction de
$E^*$ à l'ouvert complémentaire de la réunion des $v$, en un
$G$-torseur $E$ sur $\ovl X$. Les sections $\varphi_v$ et
$\varphi^*$ se correspondant via l'isomorphisme sur le disque
pointé, se recollent donc en une section $\varphi$ de
$\ad(E)\otimes\calO_X(D)$ au-dessus de $X$. Ce foncteur de
recollement est compatible à l'action de $P^{\bullet}_{a_v}$ et de
$P_a$ et définit donc une flèche entre les quotients. \item la
flèche $\beta$ est l'oubli des trivialisations génériques. \item
la flèche $\gamma$ est la restriction de $\ovl X$ à $\ovl X_v$.
\end{itemize}
Il est clair que $\beta=\gamma\circ\alpha$.

Il suffit de démontrer que $\alpha$ et $\beta$ sont des
équivalences de catégorie. Montrons d'abord que le foncteur
$$
\alpha:\prod_{v\in \ovl X-\ovl U_a}
[\calM^\bullet_{a,v}/P^\bullet_{a,v}] \rta [\calM_a(\ovl
k)/P_a(\ovl k)]
$$
est pleinement fidèle. Donnons-nous $m=(m_v)$ et $n=(n_v)$ avec
$m_v,n_v$ dans $\calM^\bullet_{a,v}$ pour toutes places $v\in \ovl
X-\ovl U_a$ et un isomorphisme de $\alpha(m)$ dans $\alpha(n)$
dans $[\calM_a/P_a]$. Cet isomorphisme consiste en un objet $p$ de
$P_a$ qui envoie $\alpha(m)$ sur $\alpha(n)$ vus comme objets de
$\calM_a$. Or les restrictions de $\alpha(m)$ et $\alpha(n)$ à
l'ouvert $\ovl U_a$ sont isomorphes à la restriction de
$(E^*,\varphi^*)$ de sorte que le $J_a$-torseur $p$ est muni d'une
trivialisation sur l'ouvert $\ovl X-\ovl U_a$. Il définit donc un
objet de $\prod_v P_{a,v}$. Ceci prouve que le foncteur $\alpha$
est pleinement fidèle.

Montrons maintenant que $\alpha$ est une foncteur essentiellement
surjectif. Soit $m$ un objet de $\calM_a$. Sur l'ouvert $\ovl
U_a$, il existe un $J_a$-torseur $p'$ tel que $p'm|_{\ovl X}=
(E^*,\varphi^*)|_{\ovl U_a}$. On peut prolonger le $J_a$-torseur
$p'$ sur $\ovl X^0$ en un $J_a$-torseur $p$ sur $\ovl X$ parce que
sur les disques pointés, n'importe quel $J_a$-torseur est trivial.
Les deux points $pm$ et $(E^*,\varphi^*)$ sont munis d'un
isomorphisme au-dessus de l'ouvert $\ovl U_a$. On en déduit les
points $m_v$ de $\calM_{a,v}^\bullet$ tels qu'en recollant les
$m_v$ avec le point base $(E^*,\varphi^*)$, on obtient $pm$.

On démontre de façon similaire que le foncteur $\beta$ est
pleinement fidèle et essentiellement surjectif. \findem

\bigskip
\noindent{\bf Complément : $2$-catégorie quotient}

\bigskip
\begin{small}

Rappelons la construction de la $1$-catégorie quotient d'un
ensemble $X$ par l'action d'un groupe $G$. La catégorie quotient
est la catégorie qui résout un problème universel
\begin{itemize}
\item $Q$ est une catégorie muni d'un foncteur $\pi:X\rta Q$, $X$
étant considéré comme une catégorie sans autres flèches que les
identités.

\item $\iota$ est un isomorphisme de foncteur
$\pi\circ \pr_X \isom \pi\circ \act$ où $\pr_X,{\rm act}:G\times
X\rightrightarrows X$ sont les application respectivement de
projection sur $X$ et de l'action de $G$ sur $X$.

\item $\iota$ vérifie une égalité de $1$-cocycle sur $G\times G\times
X$.
\end{itemize}
La catégorie quotient $Q$ est construit comme la catégorie ayant
l'ensemble des objets l'ensemble ${\rm ob}(Q)=X$, et pour tous
$x_1,x_2\in {\rm ob}(Q)$,
$$\Hom_Q(x_1,x_2)=\{q_g\mid g\in G\mbox{ tel que } gx_1=x_2\}.$$
La structure de composition dans la catégorie $Q$ vient de la
multiplication dans le groupe $G$. Pour tout $g\in G$, $x\in X$,
on a donc une flèche canonique $q_g\in \Hom_Q(x,gx)$ qui définit
l'isomorphisme entre les deux foncteurs $\pi\circ \pr_X \isom
\pi\circ \act$.

Soient maintenant $X$ une catégorie en groupoïdes et $G$ une
catégorie de Picard agissant sur $X$. Le quotient $Q$ de $X$ par
$G$ est une $2$-catégorie qui résout un problème universel :
\begin {itemize}

\item $Q$ est une $2$-catégorie munie d'un $2$-foncteur $\pi:X\rta
Q$, $X$ étant considéré comme une $2$-catégories n'ayant que des
$2$-flèches identités.

\item un objet donné dans la $1$-catégorie $\Hom(\pi\circ
\pr_X,\pi\circ \act)$ où $\pr_X,{\rm act}:G\times
X\rightrightarrows X$ sont respectivement les foncteurs de
projection sur $X$ et de l'action de $G$ sur $X$.

\item une $2$-flèche donnée au niveau de
$G^2\times X$.

\item une égalité de $2$-cocycle sur la $2$-flèche
au niveau de $G^3 \times X$.
\end{itemize}
Cette $2$-catégorie quotient $Q$ peut être construite comme suit.
\begin{itemize}
\item L'ensemble des objets de $Q$ est l'ensembles des objets de $X$.
\item Soit $x_1,x_2\in {\rm ob}(X)$, les objets de la $1$-catégorie
$\Hom_Q(x_1,x_2)$ sont les paires $q_{g,\alpha}=(g,\alpha)$ où $g$
est un objet de $G$ et $\alpha\in \Hom_X(gx_1,x_2)$.
\item Une $2$-flèche $q_{g,\alpha}\Rightarrow q_{g',\alpha'}$ est
un élément $\beta\in \Hom_G(g,g')$ tel que le triangle formé de
${\rm act}(\beta,1_{x_1}):gx_1 \rta g'x_1$, $\alpha: gx_1 \rta
x_2$ et $\alpha':g'x_1 \rta x_2$, commute.
\end{itemize}

\medskip
J.-F. Dat m'a expliqué le critère suivant pour qu'un $2$-quotient
soit équivalent à une $1$-catégorie.

\begin{lemme}
Pour que le $2$-quotient d'une catégorie en groupoïdes $X$ par
l'action d'un catégorie de Picard $Q$ soit équivalente à une
$1$-catégorie, il faut et il suffit que pour tout objet $x$ de
$X$, l'homomorphisme
$${\rm Aut}_Q(1_Q)\rta {\rm Aut}_X(x)$$
déduit de l'action de $Q$ sur $X$, soit un homomorphisme injectif.
\end{lemme}

\dem Pour que le $2$-quotient $[X/Q]$ soit équivalente à une
$1$-catégorie, il faut et il suffit que pour tous objets $x_1,x_2$
de $[X/Q]$, la $1$-catégorie des flèches de $x_1$ dans $x_2$ soit
une catégorie discrète. Par construction du $2$-quotient, ceci
revient à l'injectivité mentionnée dans l'énoncé. \findem

\end{small}

\section{Lissité}

Dans la section 1, on a un homomorphisme
$\rho_G:\pi_1(X,x)\rta\Out(\GG)$ dont l'image est $\Theta$. Soit
$\Theta_\geom$ le sous-groupe de $\Theta$ l'image du groupe
fondamental géométrique de $X$
$$\rho_G^\geom:\pi_1(\ovl X,x)\rta\Out(\GG)$$
où $\ovl X=X\otimes_k\ovl k$. Puisque $\pi_1(\ovl X,x)$ est un
sous-groupe distingué de $\pi_1(X,x)$, $\Theta_\geom$ est un
sous-groupe distingué de $\Theta$. Notons $\ovl X_\Theta^\geom$ le
revêtement fini étale galoisien de $\ovl X$ de groupe de Galois
$\Theta_\geom$ associé. Il est pointé par $x_\Theta^\geom$
au-dessus de $x$.

Nous supposons dans cette section que le groupe des caractères
$\bbX$ de $\mathbbm T$ n'a pas de vecteurs non triviaux invariants
sous $W'_\geom=W\rtimes\Theta_\geom$. C'est en particulier le cas
si $\GG$ est un groupe semi-simple. Nous allons dans ce cas
dé\-mon\-trer que le produit fibré
$\calM\times_\bbA\bbA^{\rmgssr}$ est lisse et que le morphisme
$P\times_\bbA\bbA^{\rmgssr}\rta \bbA^{\rmgssr}$ est lisse.

Soit $a\in\bbA(\ovl k)$ une caractéristique génériquement
semi-simple régulière correspondant à une section $h_a:\ovl X\rta
\car\times^{\GG_m} L_D$. Soit $\wt X_a\rta \ovl X$ le revêtement
fini plat, génériquement étale galoisien de groupe de Galois
$W'_\geom$ obtenu par l'image réciproque du morphisme $h_a:\ovl
X\rta[\car\times^{\GG_m} L_D]\otimes_k{\ovl k}$ du revêtement
$[\ovl X_\Theta^\geom\times\mathbbm t\times^{\GG_m} L_D]\rta
[\car\times^{\GG_m} L_D]\otimes_k{\ovl k}$. Notons $J_a=h_a^*[J]$
et $\Lie(J_a)$ son schéma en algèbres de Lie.

\begin{lemme}
Pour toutes caractéristiques $a\in\bbA^{\rmgssr}(\ovl k)$, on a
l'annulation $$\rmH^0(\ovl X,\Lie(J_a))=0.$$
\end{lemme}

\dem D'après 3.3, pour tout $\ovl X$-schéma $S$, les sections du
faisceau $\Lie(J_a)$ au-dessus de $S$ sont les sections de $\wt
S=S\times_{\ovl X}\wt X_a$ à valeurs dans $\mathbbm t$ qui sont
$W'_\geom$-équivariants. Ainsi $\rmH^0(\ovl X,\Lie(J_a))$ est le
groupe des sections $\wt X_a\rta \mathbbm t$ qui sont
$W'_\geom$-équivariantes. Puisque $\wt X_a$ est une courbe propre
géométriquement connexe réduite d'après l'hypothèse génériquement
semi-simple régulière $\rmH^0(\wt X,\mathbbm t)=\mathbbm t$. Or,
le sous-espace des sections $W'_\geom$-équivariantes $\mathbbm
t^{W'_\geom}=0$ s'annule de sorte que $\rmH^0(\ovl
X,\Lie(J_a))=0$. \findem

\begin{proposition}
La restriction $P\times_\bbA \bbA^{\rmgssr}$ de $P$ à l'ouvert des
caractéristiques génériquement semi-simples régulières, est lisse
au-dessus de $\bbA^{\rmgssr}$.
\end{proposition}

\dem Pour tous $a\in\bbA(\ovl k)$ l'algèbre de Lie de $P_a$ est
l'espace vectoriel $\rmH^1(\ovl X,\Lie(J_a))$. Au-dessus de
l'ouvert $\bbA^{\rmgssr}$, la dimension de cet espace vectoriel
reste constant parce que $\rmH^0(\ovl X,\Lie(J_a))$ s'annule
d'après le lemme précédent. \findem

\bigskip
D'après Biswas et Ramanan (\cf \cite{B-R}), l'espace tangent de
$\calM$ en un point $(E,\varphi)\in\calM(\ovl k)$ se calcule comme
suit. Formons le complexe à deux crans de faisceaux localement
libres sur $\ovl X$
$$
\ad(E,\varphi)=[0\rta\ad(E)\ \hfld{\ad(\varphi)}{}\ \ad(E)\otimes
\calO_X(D) \rta 0]
$$
où la flèche $\ad(\varphi)$ est donnée par $x\mapsto [x,\varphi]$
pour toutes sections locales $x$ de $\ad(E)$. L'espace tangent de
$\calM$ en le point $(E,\varphi)$ est alors le premier groupe
d'hypercohomologie $\rmH^1(\ovl X,\ad(E,\varphi))$.

\begin{proposition}
Supposons que $f(E,\varphi)=a\in\bbA^{\rmgssr}(\bar k)$
c'est-à-dire que le polynôme caractéristique de $\varphi$ soit
génériquement semi-simple régulière . Alors, $\rmH^2(\ovl
X,\ad(E,\varphi))=0$ dès que $\deg(D)>2g-2$.
\end{proposition}

\dem L'existence de la forme de Killing implique que le fibré
vectoriel $\ad(E)$ est auto-dual. Le complexe parfait $[\ad(E)\
\hfld{\ad(\varphi)}{}\ \ad(E)\otimes \calO_X(D)]$ est donc dual au
complexe $[\ad(E)\otimes \calO_X(-D)\otimes \Omega_X
\hfld{\ad(\varphi)}{}\ad(E)\otimes \Omega_X]$. Par dualité de
Serre, il suffit de démontrer que
$$\rmH^0(\ovl X,[\ad(E)\otimes \calO_X(-D)\otimes \Omega_X
\hfld{\ad(\varphi)}{}\ad(E)\otimes \Omega_X])=0
$$
dès que $\deg(D)>2g-2$.

Notons $K={\rm ker}[\ad(E)\rta \ad(E)\otimes\calO_X(D)]$ le noyau
de le flèche $\ad(\varphi):x\mapsto [x,\phi]$ ; $K$ est un
$\calO_{\ovl X}$-Module localement libre puisque qu'il est une
sous-$\calO_X$-Module de $\ad(E)$. De même, le quotient
$M=\ad(E)/K$, étant isomorphe à une sous-$\calO_{\ovl X}$-Module
de $\ad(E)\otimes \calO_{\ovl X}(D)$, est aussi localement libre.

Le commutant $K$ de $\varphi$ est un faisceau en sous-algèbres de
Lie de $\ad(E)$. Sur l'ouvert dense $U_a$ de $\ovl X$, $K$ est un
faisceau en sous-algèbres de Cartan. Puisque $K$ est localement
libre, le crochet est partout nul sur $K$. L'action de $K$ sur
$\ad(E)$ par dérivation induit donc une action de $K$ sur
$\ad(E)/K=M$ c'est-à-dire on a un homomorphisme de faisceaux en
algèbres de Lie
$$K\rta {\rm End}(M).$$
Cet homomorphisme est injectif sur un ouvert dense, donc partout
injectif puisque $K$ est localement libre. Notons $A$ le
sous-faisceau de ${\rm End}(M)$ des sections commutantes à l'image
de $K$ qui contient donc $K$ puisque le crochet est nul sur $K$.
La $\calO_X$-Algèbre $A$ étant un sous-$\calO_X$-Module de ${\rm
End}(M)$, est nécessairement localement libre comme
$\calO_X$-Module. Puisqu'elle est génériquement commutatif, elle
est partout commutative. On a donc un homomorphisme injectif de
$\calO_{\ovl X}$-Modules localement libres $K\rta A$ où $A$ a de
plus une structure de $\calO_{\ovl X}$-algèbres. On en déduit une
flèche injective
$$\rmH^0(\ovl X,K\otimes \calO_X(-D)\otimes\Omega_{\ovl X})
\rta \rmH^0(\ovl X,A\otimes \calO_X(-D)\otimes\Omega_{\ovl X}))$$
de sorte qu'il suffira de démontrer l'annulation du dernier
$\rmH^0$.

Considérons le revêtement $Y=\Spec_{\calO_{\ovl X}}(A)$ fini et
plat de $\ovl X$. Il suffit de démontrer que pour tout faisceau
inversible $L$ de degré négatif sur $X$, $\rmH^0(Y,L)=0$. Ceci est
évident car la restriction de $L$ à toute composante connexe du
normalisé de $Y$ est encore de degré strictement négatif. \findem

\section{Description de $\pi_0(P/\bbA^\heartsuit)$}

Dans cette section, on va décrire le faisceau
$\pi_0(P/\bbA^\heartsuit)$ dont la fibre au-dessus d'un point
géométrique $a\in\bbA^\heartsuit(\bar k)$ est $\pi_0(P_a)$. La
construction du fai\-sceau $\pi_0(P/\bbA^\heartsuit)$ est basée
sur l'énoncé suivant de Grothendieck (\cf [EGA IV.3] 15.6.4)

\begin{proposition}
Soient $Y$ un schéma noethérien, $f:X\rta Y$ un morphisme de type
fini, plat, à fibres géométriquement réduites. Soit $g:Y\rta X$
une section de $f$. Pour tous points $y\in Y$, notons $X_y^0$ la
composante connexe de $X_y$ contenant $g(y)$. Alors, la réunion
$X^0$ des $X_y^0$ est un ouvert de Zariski de $X$.
\end{proposition}

Voici une conséquence de cet énoncé.

\begin{proposition}
Soient $Y$ un schéma noethérien, $f:X\rta Y$ un morphisme de type
fini, plat, à fibres géométriquement réduites. Pour tout ouvert
étale $Y'\rta Y$, considérons la relation d'équivalence sur
$X(Y')={\rm Mor}_Y(Y',X)$ définie comme suit : deux sections
$$g_1,g_2:Y'\rta X'=X\times_Y Y'
$$
sont équivalentes si pour tous $y'\in Y'$, $g_1(y')$ et $g_2(y')$
appartiennent à la même composante connexe de la fibre $X_{y'}$.
Soit $\pi_0(X/Y)$ le faisceau associé au préfaisceau des sections
de $f:X\rta Y$ modulo cette relation d'équivalence. Ce faisceau
est un faisceau constructible. De plus, pour tout point $y\in Y$,
la fibre en $y$ de $\pi_0(X/Y)$ est $\pi_0(X_y)$. Il en est de
même pour les points géométriques.

\end{proposition}

\dem Pour toute section $g:Y'\rta X$ de $f:X\rta Y$ au-dessus d'un
ouvert étale $Y'$ de $Y$, on a un ouvert de Zariski $U_g$ de
$X'=X\times_Y Y'$ dont la trace sur la fibre $X_{y'}$ au-dessus de
chaque point $y'\in Y'$ est la composante connexe de $X_{y'}$
contenant le point $g(y')$. De plus, on peut trouver une telle
section $g$ passant par n'importe quel point $x$ dans le lieu de
lissité de $f$ quitte à rétrécir l'ouvert $Y'$. Par la
quasi-compacité du lieu lisse de $f$, il existe un nombre fini
d'ouverts étales $\{Y'_i\}_{i=1}^n$ et des sections $g_i:Y'_i\rta
X$ telles que la réunion disjointe des ouverts $U_{g_i}$ s'envoie
surjectivement sur le lieu lisse de $f$. On en déduit la
constructibilité de $\pi_0(X/Y)$. Cette constructibilité peut
aussi se déduire de la constructiblité de $\rmR^{2d} f_! \QQ_\ell$
où $d$ est la dimension relative de $f$.

Soit $y$ un point de $Y$. Il faut montrer que $\pi_0(X_y)$ est la
limite inductive des ensembles des classes d'équivalence de
$X(Y')$ quand $Y'$ parcourt l'ensemble des voisinages étales
pointés par le point $y$.

Pour toute section $g:Y'\rta X$ de $f$  au-dessus d'un voisinage
étale $Y'$ pointé par $y$, le point $g(y)$ appartient à une unique
composante connexe de la fibre $X_y$. Si $g_1, g_2:Y'\rta X$ sont
équivalentes, les points $g_1(y)$ et $g_2(y)$ appartiennent à la
même composante connexe de $X_y$ de sorte qu'on a bien défini une
application de l'ensemble des sections de $f$ au-dessus du
voisinage $Y'$ modulo l'équivalence dans l'ensemble $\pi_0(X_y)$.
Ces applications sont compatibles à la restriction d'un voisinage
$Y'$ pointé par $y$ à un voisinage $Y''$ plus petit si bien
qu'elles s'organisent en une application de la limite inductive
dans $\pi_0(X_y)$.

Cette application est surjective. En effet, soit
$\alpha\in\pi_0(X_y)$. Puisque la fibre $X_y$ est réduite, il
existe un point lisse $x\in X_y$ dans la composante connexe
$\alpha$. Puisque c'est un point lisse, on peut l'étendre à un
voisinage étale $Y'$ pointé par $y$. L'image de l'application de
$X(Y')$, modulo la relation d'équivalence, dans $\pi_0(X_y)$
contient donc $\alpha$.

Soient maintenant $g_1,g_2$ deux sections de $f$ au-dessus d'un
voisinage étale $Y'$ pointé $y$ telles que $g_1(y)$ et $g_2(y)$
appartiennent à la même composante connexe de $X_y$. Soit $X'_1$,
respectivement $X'_2$, l'ouvert de Zariski de $X'=X\times_Y Y'$,
dont la trace sur toutes les fibres $X_{y'}$ de $f$ sont la
composante connexe de $g_1(y')$, respectivement $g'_2(y')$. Les
deux ouverts $X'_1$ et $X'_2$ ont donc la même trace sur $X_y$ si
bien que leur intersection est un ouvert non-vide de $X$. Puisque
$f$ est une application ouverte, l'image de $X'_1\cap X'_2$ est un
ouvert de $Y$ contenant $y$ qu'on va noter $Y''$. Pour tout $y'\in
Y''$, les traces de $X'_1$ et $X'_2$ sur la fibre $X_{y'}$ en sont
deux composantes connexes ayant une intersection non vide. On a
donc $X'_1\cap X_{y'}=X'_2\cap X_{y'}$ pour tout $y'\in Y''$. On
en déduit que $X'_1\cap f^{-1}(Y'')=X'_2\cap f^{-1}(Y'')$
c'est-à-dire les restrictions de $g_1$ et $g_2$ à $Y''$ sont
équivalentes. L'application de la limite inductive dans
$\pi_0(X_y)$ est injective et donc bijective.

La même démonstration vaut si on remplace un point $y\in Y$ par un
point géométrique $\ovl y$ et les voisinages pointés par $y$ par
les voisinages pointés par $\ovl y$. La fibre géométrique
$\pi_0(X/Y)$ en $\bar y$ est alors le groupe des composantes
connexes $\pi_0(X_{\ovl y})$ de la fibre géométrique $X_{\ovl
y}$.\findem

\bigskip
D'après la proposition 5.2, $P$ est un champ de Picard relatif
lisse au-dessus de $\bbA^\heartsuit$. Il résulte alors de la
proposition précédente qu'il existe un faisceau
$\pi_0(P/\bbA^{\rmgssr})$ pour la topologie étale de
$\bbA^\rmgssr$ tel que pour tout point (géométrique) $a$ de
$\bbA^{\rmgssr}$, la fibre de $\pi_0(P/\bbA^{\rmgssr})$ en $a$ est
le groupe des composantes connexes $\pi_0(P_a)$ de la fibre de $P$
en $a$.

On commence par décrire les fibres du faisceau
$\pi_0(P/\bbA^{\rmgssr})$. Soit $a$ un point géométrique de
$\bbA^{\rmgssr}$. Soit $J_a=h_a^*[J]$ le schéma en groupes lisse
image réciproque par $h_a:\ovl X\rta[\car/\nu_\chi(\GG_m)]$ du
centralisateur régulier $[J]$. Notons $U_a$ l'image inverse du
lieu des caracté\-ristiques semi-simples régu\-lières ; cet ouvert
est non vide d'après l'hypothèse génériquement semi-simple
régulière. La restriction de $J_a$ à $U_a$ est un schéma en tores.
Le schéma en groupes $J_a$ étant lisse, il existe un sous-schéma
ouvert $J_a^0$ qui en tous points géométriques $x$ de $\ovl X$, la
fibre de $(J_a^0)_x$ est la composante neutre de $(J_a)_x$.
L'immersion ouverte $J_a^0 \rta J_a$ est un isomorphisme au-dessus
de l'ouvert $\ovl U_a$ puisque les tores sont connexes. Le
quotient $J_a/J_a^0$ est donc un faisceau gratte-ciel de fibres
finies, concentré dans le fermé de dimension zéro $\ovl X-U_a$.

Soit $P'_a$ le champ de Picard des $J_a^0$-torseurs sur $\ovl X$.
Quand $a$ varie, les $P'_a$ s'organisent aussi  en un champ de
Picard relatif $P'$ au-dessus de $\bbA^\heartsuit$. En effet,
au-dessus de $\bbA^\heartsuit\times X$, on a le schéma en groupes
affines commutatifs lisse $J^0$ qui est la composante neutre de
$J$. La construction de $P'$ est alors la même que celle de $P$ en
remplaçant $J$ par $J^0$. Par ailleurs, $P'$ et $P$ ont le même
espace tangent à cause de la suite exacte
$$\rmH^0(\ovl X,J_a/J_a^0)\rta \rmH^1(\ovl X,J_a^0) \rta \rmH^1(\ovl
X,J_a)\rta \rmH^1(\ovl X,J_a/J_a^0)
$$
où $\rmH^0(\ovl X,J_a/J_a^0)$ est un groupe fini et où
$\rmH^1(\ovl X,J_a/J_a^0)=0$. On en déduit que $P'$ est également
lisse au-dessus de $\bbA^\heartsuit$ de sorte qu'on dispose d'un
faisceau $\pi_0(P'/\bbA^\heartsuit)$ dont les fibres sont les
groupes $\pi_0(P'_a)$.

\begin{proposition}
On a la suite exacte
$$\rmH^0(\ovl X,J_a/J_a^0)\rta \pi_0(P'_a)\rta \pi_0(P_a)\rta 0.$$
\end{proposition}

\dem Dans la suite exacte
$$
\rmH^0(\ovl X,J_a/J_a^0)\rta \rmH^1(\ovl X,J_a^0) \rta \rmH^1(\ovl
X,J_a)\rta \rmH^1(\ovl X,J_a/J_a^0)=0
$$les groupes $\rmH^1(\ovl
X,J_a^0)$ et $\rmH^1(\ovl X,J_a)$ sont les groupes des classes
d'isomorphisme des catégories de Picard $P'_a(\ovl k)$ et
$P_a(\ovl k)$. La surjectivité essentielle du foncteur $P'_a(\ovl
k)\rta P_a(\ovl k)$ implique la surjectivité essentielle de sa
restriction aux composantes neutres ${P'}_a^0(\ovl k)\rta
P_a^0(\ovl k)$. On en déduit la suite exacte
$$\rmH^0(\ovl X,J_a/J_a^0)\rta \pi_0(P'_a)\rta \pi_0(P_a)\rta 0$$
qu'on voulait. \findem

\bigskip Cette suite exacte permet en principe de
ramener le calcul de $\pi_0(P_a)$ à celui de $\pi_0(P'_a)$ du
moment qu'on sait décrire la flèche $\rmH^0(\ovl X,J_a/J_a^0)\rta
\pi_0(P'_a)$. Il est possible de la calculer dans le cas des
groupes classiques, et nous allons présenter en détail ce calcul
dans le cas $\SL(2)$ dans §11. Malheureusement, nous ne savons pas
formuler une description générale. Cette difficulté va
réapparaître dans la démonstration du théorème principal 10.4 où
heureusement, il est possible de la contourner à l'aide d'un lemme
de Kottwitz. Il serait plus satisfaisant d'avoir une description
directe de cette flèche.

En revanche, il est possible de décrire $\pi_0(P'_a)$ de façon
générale à l'aide d'un autre lemme de Kottwitz. Le schéma en
groupes $J^0_a$ est un schéma en groupes lisse à fibres connexes
sur $\ovl X$ et sa restriction $J_a^0|_{U_a}=J_a|_{U_a}$ à
l'ouvert $U_a$ est un tore. Ce tore admet un modèle de Néron
connexe canonique ${\rm Ner}^0(J_a|_{U_a})$ qui est une schéma en
groupes lisse à fibre connexe sur $\ovl X$, prolongeant le schéma
en tores $J_a|_{U_a}$, et universel pour cette propriété. En
particulier, on a un homomorphisme canonique
$$J^0_a\rta {\rm Ner}^0(J_a|_{U_a}).$$

\begin{proposition}
L'homomorphisme du champ de Picard des $J^0_a$-torseurs dans celui
des ${\rm Ner}^0(J_a|_{U_a})$-torseurs induit un isomorphisme sur
leurs groupes des composantes connexes.
\end{proposition}

\dem L'homomorphisme $J^0_a\rta {\rm Ner}^0(J_a|_{U_a})$ est
injectif comme un homomorphisme de faisceaux. Considérons la suite
exacte
$$0\rta J^0_a\rta {\rm Ner}^0(J_a|_{U_a}) \rta K\rta 0$$
où $K$ est un faisceau supporté par $\ovl X-U_a$. On a
$\rmH^1(\ovl X,K)=0$ et $\rmH^0(\ovl X,K)$ est un groupe
algébrique affine connexe. La proposition se déduit donc de la
suite exacte longue
$$\rmH^0(\ovl X,K)\rta \rmH^1(\ovl X,J^0_a)\rta
\rmH^1(\ovl X,{\rm Ner}^0(J_a|_{U_a})) \rta \rmH^1(\ovl X,K)=0$$
qui se déduit de la suite exacte courte précédente de la même
manière que 6.2. \findem

\bigskip
On se ramène au calcul du groupe des composantes connexes du champ
de Picard des torseurs sous le modèle de Néron connexe ${\rm
Ner}^0(J_a|_{U_a})$. Pour cela, nous utilisons un lemme général du
à Kottwitz (\cf \cite{K-Iso1}).

Fixons un point géométrique $u$ de l'ouvert $U_a$ et notons
$\Gamma=\pi_1(U_a,u)$ le groupe fondamental de $U_a$ pointé en
$u$. La catégorie des tores sur $U_a$ est alors équivalente à la
catégorie des $\ZZ$-modules libres munis d'une action finie de
$\Gamma$. On associe à un tore $\calT$ sur $U_a$ la fibre en $u$
du faisceau $\bbX_*(\calT)$ des cocaractères de $\calT$, munie de
l'action naturelle de $\pi_1(U_a,u)$. Un tore $\calT$ est dit
induit s'il existe une base du $\ZZ$-module $\bbX_*(\calT)_{u}$
telle que l'action de $\Gamma$ sur $\bbX_*(\calT)_{u}$ se déduit
d'une action transitive de $\Gamma$ sur cette base. Le lemme 2.2
de \cite{K-Iso1} peut s'énoncer comme suit.

\begin{lemme} Soit $A$ un foncteur de la catégorie des
tores sur $U_a$ dans la catégorie des groupes abéliens vérifiant
les propriétés suivante :
\begin{enumerate}

\item le foncteur $A$ est exact à droite c'est-à-dire il
transforme une suite exacte de tores $1\rta \calT_1\rta
\calT_2\rta \calT_3\rta 0$ en une suite exacte à droites de
groupes abéliens
$$A(\calT_1)\rta A(\calT_2) \rta A(\calT_3)\rta 0,$$
\item si $\calT$ est un tore induit, il existe un isomorphisme
canonique $A(\calT)\isom \ZZ$.
\end{enumerate}
Soit $u$ un point géométrique de $U_a$. Il existe un isomorphisme
canonique du foncteur $\calT\mapsto A(\calT)$ dans le foncteur
$\calT\mapsto [\bbX_*(\calT)_u]_\Gamma$ qui associe à un tore
$\calT$ le groupe abélien des coinvariants de
$\Gamma=\pi_1(U_a,u)$ dans le groupe des cocaractères
$\bbX_*(\calT)_u$.
\end{lemme}

On va appliquer ce lemme de Kottwitz au foncteur $A$ construit
comme suit. Un tore $\calT$ sur $U_a$ se prolonge canoniquement en
un $\ovl X$-schéma en groupes lisses à fibres connexe ${\rm
Ner}^0(\calA)$, le modèle de Néron connexe de $\calT$. On pose
$$A(\calT)=\pi_0({\rm Tors}({\rm Ner}^0(\calT)))$$
le groupe des composantes connexes du champ de Picard des torseurs
sous ${\rm Ner}^0(\calT)$ sur $\ovl X$.

\begin{lemme}
Le foncteur $\calT\mapsto \pi_0({\rm Tors}({\rm Ner}^0(\calT)))$
est un foncteur exact à droite. De plus, pour les tores induits
$\calT$, on a un isomorphisme canonique $\pi_0({\rm Tors}({\rm
Ner}^0(\calT)))\isom\ZZ$.
\end{lemme}

\dem Soit $1\rta \calT_1\rta \calT_2\rta \calT_3\rta 1$ une suite
exacte de tores sur $U_a$. On en déduit une suite exacte entre les
modèles de Néron localement de type fini
$$
1\rta {\rm Ner}(\calT_1)\rta {\rm Ner}(\calT_2)\rta {\rm
Ner}(\calT_3)\rta 1.
$$
Soit $x\in \ovl X-U_a$ un point dans le complémentaire de
$U_a$, soit $I_x$ le sous-groupe d'inertie en $x$. D'après
\cite{Ra}, le groupe de composantes connexes de la fibre en $x$ du
modèle de Néron est
$$\pi_0({\rm Ner}(\calT_\alpha))=(\bbX_*(\calT_\alpha))_{I_x}.$$
pour tous $\alpha\in\{1,2,3\}$. On en déduit une suite exacte de
faisceaux
$$
0\rta {\rm Ner}^0(\calT_\alpha) \rta {\rm Ner}(\calT_\alpha)\rta
\bigoplus_{x\in \ovl X-U_\alpha} (\bbX_*(\calT_\alpha))_{I_x} \rta
0
$$
où $(\bbX_*(\calT_\alpha))_{I_x}$ est vu comme un gratte-ciel
planté en $x$. On a maintenant une suite exacte à droite
$$
\bbX_*(\calT_1)_{I_x}\rta \bbX_*(\calT_2)_{I_x}\rta
\bbX_*(\calT_3)_{I_x} \rta 0
$$
puisque le foncteur des coinvariants est exacte à droite.
En utilisant le diagramme du serpent, on voit que parmi les
homomorphismes entre leurs modèles de Néron connexes
$$
{\rm Ner}^0(\calT_1)\hfld{\phi_1}{} {\rm Ner}^0(\calT_2)
\hfld{\phi_2}{} {\rm Ner}^0(\calT_3)
$$
$\phi_1$ est injectif, $\phi_2$ est surjectif mais ${\rm
im}(\phi_1)\not= {\rm ker}(\phi_2)$. Plus précisément, on a
$$
{\rm ker}(\phi_2)/{\rm im}(\phi_1)=\bigoplus_{x\in\ovl X-U_\alpha}
{\rm ker}(\bbX_*(\calT_1)_{I_x}\rta \bbX_*(\calT_2)_{I_x}).
$$
Posons $\calT'_1={\rm ker}(\phi_2)$. On a alors deux suites
exactes :
$$
0\rta \calT'_1 \rta {\rm Ner}^0(\calT_2) \hfld{\phi_2}{} {\rm
Ner}^0(\calT_3)\rta 0
$$
et
$$
0\rta {\rm Ner}^0(\calT_1) \rta \calT'_1 \rta \bigoplus_{x\in\ovl
X-U_\alpha} {\rm ker} (\bbX_*(\calT_1)_{I_x}\rta
\bbX_*(\calT_2)_{I_x}) \rta 0.
$$
Dans la première suite exacte, $\calT'_1$ est un schéma en groupes
sur $\ovl X$ à fibre générique géométriquement connexe de sorte
que $\rmH^2(\ovl X,\calT'_1)=0$ d'après un théorème de Tsen
\cite{Se}. On en déduit la suite exacte
$$\rmH^1(\ovl X,\calT_1)\rta \rmH^1(\ovl X,{\rm Ner}^0(\calT_2))\rta
\rmH^1(\ovl X,{\rm Ner}^0(\calT_3))\rta 0.
$$
Dans la seconde suite exacte, $\bigoplus_{x\in\ovl X-U_\alpha}
{\rm ker}(\bbX_*(\calT_1)_{I_x}\rta \bbX_*(\calT_2)_{I_x})$ est
supporté par un schéma de dimension nulle de sorte que son
$\rmH^1$ est nul. On en déduit la surjectivité de la flèche
$$\rmH^1(\ovl X,{\rm Ner}^0(\calT_1))\rta\rmH^1(\ovl
X,\calT'_1).
$$
On en déduit une suite exacte à droite
$$
\rmH^1(\ovl X,{\rm Ner}^0(\calT_1))\rta\rmH^1(\ovl X, {\rm
Ner}^0(\calT_2))\rta \rmH^1(\ovl X,{\rm Ner}^0(\calT_3))\rta 0
$$
qui induit une suite exacte à droite des $\pi_0$
$$
\pi_0({\rm Tors}({\rm Ner}^0(\calT_1)))\rta \pi_0({\rm Tors} ({\rm
Ner}^0(\calT_2)))\rta \pi_0({\rm Tors}({\rm Ner}^0(\calT_3)))\rta
0.
$$
Le foncteur $\calT\mapsto \pi_0({\rm Tors}({\rm Ner}^0(\calT)))$
est donc un foncteur exact à droite.

Considérons maintenant un tore induit $\calT$. Il existe donc une
base de $\bbX_*(T)_{\ovl u}$ munie d'une action transitive de
$\pi_1(U_a,\ovl u)$. Cet ensemble fini d'une action de
$\pi_1(U_a,\ovl u)$, définit un revêtement fini étale connexe
$\pi_U:\widetilde U\rta U_a$. On vérifie alors que
$\calT=(\pi_U)_*\GG_m$. Par normalisation on obtient un revêtement
fini ramifié $\widetilde X \rta\ovl X$. On vérifie que ${\rm
Ner}^0(\calT)=\pi_*\GG_m$ de sorte que
$${\rm Tors}({\rm Ner}^0(\pi_*\GG_m))={\rm Pic}(\widetilde X)$$
La courbe $\widetilde X$ étant irréductible, on a $\pi_0({\rm
Pic}(\widetilde X))=\ZZ$. \findem

\begin{corollaire} Soit $a\in \bbA^\heartsuit(\ovl k)$ et $u\in
U_a(\ovl k)$. On a un isomorphisme canonique $\pi_0(P'_a)\simeq
[\bbX_*(J_a)_u]_\Gamma$ où $\bbX_*(J_a)_u$ est la fibre en $u$ du
faisceau des cocaractères du tore $J_a|_{U_a}$ et où
$\Gamma=\pi_1(U_a,u)$.
\end{corollaire}

Le corollaire précédent décrit les fibres géométriques
$\pi_0(P'_a)$ du faisceau $\pi_0(P'/\bbA^\heartsuit)$ à l'aide de
la monodromie du tore $J_a|_{U_a}$. Nous allons maintenant
regarder comment cette monodromie varie en fonction de $a$.

Notons $U$ l'image réciproque de $\car^{\rm ssr}\times^{\GG_m}
L_D$ par la flèche tautologique
$$X\times \bbA^\heartsuit\rta \car\times^{\GG_m}L_D.$$
Par définition de $\bbA^\heartsuit$, les fibres du morphisme
$U\rta \bbA^\heartsuit$ sont non vides et connexes.

Soit $X_\Theta$ le revêtement fini étale galoisien de $X$ de
groupe de Galois $\Theta$ associé à $\rho_G:\pi_1(X,x)\rta\Theta$
comme dans la section 1. On a alors un morphisme fini
$$(X_\Theta\times{\mathbbm t})\times^{\GG_m}L_D \rta
\car\times^{\GG_m}L_D
$$
qui au-dessus de l'ouvert $\car^{\rm ssr}\times^{\GG_m}L_D$ est
étale galoisien de groupe de Galois $W'=W\rtimes \Theta$. Par
image réciproque, on obtient ainsi un $W'$-torseur
$$\widetilde U\rta U$$
au-dessus de la $\bbA^\heartsuit$-courbe $U$. Le morphisme
$\widetilde U\rta \bbA^\heartsuit$ est un morphisme lisse de sorte
qu'on dispose d'un faisceau d'ensembles $\pi_0(\widetilde
U/\bbA^\heartsuit)$ sur la topologie étale de $\bbA^\heartsuit$
dont la fibre au-dessus de $a\in\bbA^\heartsuit(\ovl k)$ est
l'ensemble des composantes connexes de $\widetilde U_a$ d'après la
proposition 6.2. Comme $\widetilde U$ est un $W'$-torseur
au-dessus de $U$ et $U$ ont des fibres non vides et connexes, $W'$
agit sur le faisceau $\pi_0(\widetilde U/\bbA^\heartsuit)$ et
cette action est transitive fibres à fibres.

L'énoncé suivant m'a été suggéré par Drinfeld.

\begin{proposition} Il existe un isomorphisme canonique de
faisceaux
$$\ZZ[\pi_0(\widetilde U/\bbA^\heartsuit)]\otimes_{W'}
{\mathbbm X}^\sv \rta \pi_0(P'/\bbA^\heartsuit).$$
\end{proposition}

\dem Commençons par construire la flèche canonique. Localement
pour la topologie étale, une section de $\pi_0(\widetilde
U/\bbA^\heartsuit)$ est représentée par une section de $\widetilde
U$. Il suffit donc de construire à partir d'une section
$$\widetilde u_S:S\rta \widetilde U_{a_S}=\widetilde U\times_{\bbA^\heartsuit,a_S} S$$
au-dessus d'un ouvert étale $a_S:S\rta \bbA^\heartsuit$, un
homomorphisme de faisceaux dépendant de $\widetilde u_S$
$$\mathbbm X^\sv_S \rta \pi_0(P'/S)$$
puis de vérifier que cet homomorphisme est invariant pour la
relation d'équi\-valence 6.2 et qu'il vérifie une certaine
propriété $W'$-équivariante qui sera précisée.

D'après 3.8, l'image réciproque de $J^0$ par le morphisme
$\widetilde U_{a_S}\rta [\car/\GG_m]$ est canoniquement isomorphe
au tore constant $\widetilde U_{a_S}\times {\mathbbm T}$ de sorte
qu'on a un isomorphisme entre leurs faisceaux des cocaractères
$$\Hom_{\widetilde U_{a_S}}(\GG_m,J^0)\isom \bbX^\sv_{\widetilde U_{a_S}}.$$
Avec la donnée de la section $\widetilde u_S$, on a un
isomorphisme canonique
$${\widetilde u}_S^*\Hom_{U_{a_S}}(\GG_m,J^0)\isom \bbX^\sv_S.$$
Soit $J^0_{a_S}$ l'image réciproque de $J^0$ par le morphisme
$X\times S\rta [\car/\GG_m]$. Soit $u_S:S\rta U_{a_S}$ l'image de
$\widetilde u_S$. On a alors une trivialisation de
$\Hom_{U_a}(\GG_m,J^0_a)$ sur le complété formel $U_{a_S}$ le long
de la section $u_S$.

Un cocaractère $\mu\in\bbX^\sv$ induit alors sur le complété
formel $X_{S,u_S}$ de $X\times S$ le long de la section $u_S$, un
homomorphisme de $\GG_m$ dans la restriction de $J^0_{a_S}$. Le
fibré inversible $\calO_{X\times S}(u_S)$ fournit un
$\GG_m$-torseur sur $X\times S$ trivialisé sur le complémentaire
de la section $u_S$, de sorte qu'en poussant par le caractère
$\mu$, on obtient un $J^0_{a_S}$-torseur sur $X_{S,u_S}$ muni
d'une trivialisation sur le complémentaire de la section $u_S$.
D'après le théorème de recollement formel de Beauville-Laszlo (\cf
\cite{BL}), on peut recoller ce $J^0_{a_S}$-torseur sur
$X_{S,u_S}$ avec le $J^0_{a_S}$-torseur neutre sur le
complémentaire de $u_S$ pour produire un $J^0_{a_S}$-torseur sur
$X\times S$. On a donc produit un $S$-point de $P'_S$ à partir de
$\mu\in\mathbbm X^\sv$.

Cette construction fournit un homomorphisme
$$\lambda(\widetilde u_S):\bbX^\sv_S\rta P'_S$$
dépendant de la section $\widetilde u_S$. En composant avec
l'homomorphisme évident $P'_S \rta \pi_0(P'/S)$, on a un
homomorphisme de faisceaux
$$\pi_0(\lambda(\widetilde u_S)):\bbX^\sv_S \rta \pi_0(P'/S).$$

Soit $\widetilde u'_S:S\rta \widetilde U_{a_S}$ une section
équivalente à $\widetilde u_S$ au sens de 6.2. Les deux
homomorphismes
$$\pi_0(\lambda(\widetilde u_S)), \pi_0(\lambda(\widetilde u'_S)):
\bbX^\sv_S \rta \pi_0(P'/S)$$ sont alors égaux. Ceci se démontre
par un argument d'homotopie comme suit. D'après 6.1, il existe
l'ouvert de Zariski $\widetilde U_S(\widetilde u_S)$ de
$\widetilde U_S$ dont la trace sur chaque fibre $\widetilde U_s$
de $\widetilde U_S$ au-dessus d'un point $s\in S$, est la
composante neutre de $U_s$ contenant $\widetilde u_S(s)$. Du fait
que $\widetilde u_S$ et $\widetilde u'_S$ sont équivalentes,
$\widetilde u'_S$ se factorise aussi à travers l'ouvert
$\widetilde U_S(\widetilde u_S)$. Au-dessus de $\widetilde
U_S(\widetilde u_S)$, on a une flèche canonique
$$\bbX^\sv_{\widetilde U_S(\widetilde u_S)} \rta \pi_0(P'/
\widetilde U_S(\widetilde u_S))
$$
à cause de la projection vers $\widetilde U_S$. Les flèches
$\pi_0(\lambda(\widetilde u_S)), \pi_0(\lambda(\widetilde u'_S))$
s'obtiennent en restreignant cette dernière aux sections
$\widetilde u_S$ et $\widetilde u'_S$. Puisque le morphisme
$\pr_S: \widetilde U_S(\widetilde u_S)\rta S$ est un morphisme
lisse à fibres géométriquement connexes, $\pr_S^*$ induit un
isomorphisme
$$\Hom\left(\bbX^\sv_S,\pi_0(P'/S)\right)\rta
\Hom\bigl(\bbX^\sv_{\widetilde U_S(\widetilde u_S)},
\pi_0(P'/\widetilde U_S(\widetilde u_S))\bigr)
$$
d'après le lemme 4.2.5
de \cite{BBD}. L'égalité de $\pi_0(\lambda(\widetilde u_S)),
\pi_0(\lambda(\widetilde u'_S))$ s'en déduit.

Par ailleurs, 3.8 implique que pour tout $w'\in W'$, pour tout
$\mu\in \mathbbm X^\sv$, on a
$$\lambda(\widetilde u_S)(\mu)=\lambda(w'\widetilde u_S)(w'\mu).$$

La conjonction de ce qui précède implique qu'on a un homomorphisme
$$\ZZ[\pi_0(\widetilde U/\bbA^\heartsuit)]\otimes_{W'} \mathbbm X^\sv
\rta \pi_0(P'/\bbA^\heartsuit).
$$
Pour vérifier qu'il est un isomorphisme, il suffit de le faire
fibres à fibres et c'est ce que dit le corollaire 6.7. \findem

\section{Stratification}

Soit $a\in\bbA^{\rmgssr}(\bar k)$. La fibre $\pi_0(\widetilde
U/\bbA^\heartsuit)_a$ est un espace homogène sous le groupe $W'$.
Les sous-groupes d'isotropie de $W'$ agissant sur cette fibre
forment une classe de $W'$-conjugaison de sous-groupes de $W'$
qu'on va noter $[\Sigma'_a]$.

Le faisceau $\pi_0(\widetilde U/\bbA^\heartsuit)$ étant
constructible, il existe une stratification de $\bbA^\heartsuit$
telle que sa restriction aux strates est localement constante. Si
les strates sont connexes, ce qu'on peut supposer, les classes de
$W$-conjugaison $[\Sigma'_a]$ sont constantes quand le point
géométrique $a$ varie le long de l'une de ces strates de sorte
qu'à chaque strate est associée une classe de $W'$-conjugaison de
sous-groupes de $W'$. Dans cette section, on va construire
explicitement cette stratification dont les strates sont indexées
par certaines classes de $W'$-conjugaison de sous-groupes de $W'$
vérifiant une contrainte qui sera précisée dans le lemme 7.1.

Soit $x$ un point géométrique de $X$. Dans la section 1, on a un
homomorphisme $\rho_G:\pi_1(X,x)\rta{\rm Out}(\mathbbm G)$ dont
l'image est un groupe fini $\Theta$. On a aussi un revêtement
étale galoisien $X_\Theta$ de $X$ de groupe de Galois $\Theta$
pointé par un point géométrique $x_\Theta$ au-dessus de $x$
associé à $\rho_G$. Nous avons noté $\Theta_\geom$ l'image du
groupe fondamental géométrique $\pi_1(\ovl X,x)$ dans $\Theta$ et
$W'_\geom$ le sous-groupe $W\rtimes\Theta_\geom$ de
$W'=W\rtimes\Theta$.

Soit $a\in\bbA^\heartsuit(\ovl k)$ et supposons que $x\in U_a$ où
$U_a$ est la fibre de $U$ en $a$ avec $U$ et $\widetilde U$ comme
dans 6.8. Choisissons un point $\widetilde x$ de $\widetilde U_a$
au-dessus de $x$. Le groupe $W'$ agit transitivement sur
l'ensemble des composantes connexes de $\widetilde U_a$. Soit
$\Sigma_a(\widetilde x)$ le sous-groupe de $W'$ qui stabilise la
composante connexe de $\widetilde U_a$ contenant le point
$\widetilde x$. La classe de $W'$-conjugaison $[\Sigma_a]$ de
$\Sigma_a(\widetilde x)$ ne dépend pas du choix de $\widetilde x$.
Par construction, on a un morphisme $\widetilde U_a\rta X_\Theta$.

Supposons que $\widetilde x$ est au-dessus du point $x_\Theta$ de
$X_\Theta$. L'homomorphisme $\pi_1(U_a,x) \rta W'$ donné par
$\widetilde x$ relève alors l'homomorphisme
$\rho_G^\geom:\pi_1(\ovl X,x) \rta \Theta$ donné par $\widetilde
x_\Theta$. En particulier, l'image de $\Sigma'_a(\widetilde x)$
dans $\Theta$ est $\Theta_\geom$. Puisque $\Theta_\geom$ est un
sous-groupe distingué de $\Theta$, tous les $W'$-conjugués de
$\Sigma'_a(\widetilde x)$ ont l'image $\Theta_\geom$ dans
$\Theta$. On retient la contrainte suivante sur la classe de
$W'$-conjugaison $[\Sigma'_a]$.

\begin{lemme}
Pour tout $a\in\bbA^{\heartsuit}(\ovl k)$, tous les membres de la
classe de $W'$-conjugaison $[\Sigma'_a]$ ont l'image
$\Theta_\geom$ dans $\Theta$.
\end{lemme}

Pour tout sous-groupe $\Sigma'$ de $W'$ dont l'image dans $\Theta$
est $\Theta_\geom$, considérons la $\calO_X$-Algèbre des fonctions
$\Sigma'$-invariantes sur $X_\Theta\times \mathbbm t$. Son
spectre, désigné par $(X_\Theta\times \mathbbm t)/\Sigma'$, est
muni d'une action de $\GG_m$ induite de l'action de $\GG_m$ sur
$\mathbbm t$. Soit $\mathbbm B_{\Sigma'}$ le schéma qui représente
l'ensemble des sections de
$$h_b:X\rta ((X_\Theta\times \mathbbm t)/\Sigma')\times^{\GG_m}
L_D
$$
où $L_D$ est le $\GG_m$-torseur associé à $\calO_X(D)$. On a un
morphisme
$$\pi_{\Sigma'}:((X_\Theta\times \mathbbm t)/\Sigma')\times^{\GG_m}
L_D\rta ((X_\Theta\times \mathbbm t)/W')\times^{\GG_m} L_D$$ qui
induit un morphisme
$$\pi_{\Sigma'}^X:\mathbbm B_{\Sigma'} \rta \bbA.$$

Considérons l'ouvert ${\mathbbm B}^{\rmgssr}_{\Sigma'}$ l'image
récipro\-que de $\bbA^{\rmgssr}$ dans ${\mathbbm B}_{\Sigma'}$.
Pour tout $b\in {\mathbbm B}_{\Sigma'}^\rmgssr(\ovl k)$, le
revêtement
$$X_\Theta\times \mathbbm t\rta (X_\Theta\times \mathbbm
t)/\Sigma'$$ induit par image réciproque un revêtement fini
$\widetilde X_b\rta \ovl X$ qui est génériquement étale. Nous
allons noter ${\mathbbm B}_{\Sigma'}^{\rm max}$ l'ouvert de
${\mathbbm B}^\rmgssr$ qui consiste en les points $b\in {\mathbbm
B}_{\Sigma'}^\rmgssr(\ovl k)$ tels que le revêtement $\widetilde
X_b\rta \ovl X$ est irréductible.

\begin{proposition}
Le morphisme $\pi^X_{\Sigma'}:{\mathbbm B}^{\rmgssr}_{\Sigma'}\rta
\bbA^{\rmgssr}$ est un morphisme fini et net. En particulier,
${\mathbbm B}^{\rmgssr}_{\Sigma'}$ est représentable.

Soit $\ovl \bbA_{\Sigma'}$ le sous-schéma fermé de $\bbA^\rmgssr$
image du morphisme fini $\pi^X_{\Sigma'}$. Il existe alors un
ouvert $\bbA_{\Sigma'}$ de $\ovl \bbA_{\Sigma'}$ tel que
$(\pi^X_{\Sigma'})^{-1}(\bbA_{\Sigma'})$ est l'ouvert ${\mathbbm
B}_{\Sigma'}^{\rm max}$ de $\mathbbm B^\rmgssr_{\Sigma'}$. De
plus, le morphisme restreint ${\mathbbm B}_{\Sigma'}^{\rm max}\rta
\mathbbm A_{\Sigma'}$ est un morphisme fini étale de degré $|{\rm
Nor}_{W'}(\Sigma')|/|\Sigma'|$.
\end{proposition}

\dem Le lemme suivant va être utilisé à plusieurs reprises dans la
démonstration de la proposition.

\begin{lemme}
Soit $S$ un $k$-schéma normal intègre. Soient $V$ et $\widetilde
V$ des $S$-schémas et $v:\widetilde V\rta V$ un $S$-morphisme
fini. Soit $h:S\rta V$ une section de $V$. Supposons qu'il y a un
ouvert non vide $S'$ de $S$ sur lequel $h'=h|_{S'}$ se relève en
une section $\widetilde h':S'\rta \widetilde V\times_S S'$. Alors
la section $\widetilde h'$ se prolonge de façon unique en une
section $\widetilde h:S\rta \widetilde V$ qui relève $h$.
\end{lemme}

\dem Soit $h(S)$ l'image de la section $h$ et $v^{-1}(h(S))$ son
image réciproque dans $\widetilde V$. Puisque $v$ est un morphisme
fini, $v^{-1}(h(S))$ est un $S$-schéma fini. La section
$\widetilde h'$ de $\widetilde V$ au-dessus de l'ouvert $S'$,
relève $h'=h|_S$ si bien qu'elle est à l'image dans
$v^{-1}(h(S))$. Notons $Z'$ l'image de $\widetilde h'$ dans
$v^{-1}(h(S))$ et $Z$ l'adhérence schématique de $Z'$ dans
$v^{-1}(h(S))$. Puisque $v^{-1}(h(S))$ est un $S$-schéma fini, $Z$
est un $S$-schéma fini intègre bi-rationnel à $S$. Puisque $S$ est
supposé normal et intègre, ceci implique que $Z\rta S$ est un
isomorphisme. En l'inversant, on obtient la section $\widetilde h$
de $\widetilde V$ qui relève la section $h$ de $V$. \findem

\bigskip
Démontrons d'abord que le morphisme $\pi^X_{\Sigma'}:{\mathbbm
B}^{\rmgssr}_{\Sigma'}\rta \bbA^{\rmgssr}$ est {\em quasi-fini}.
Soit $a\in\bbA^{\rmgssr}(\ovl k)$ une caractéristique
génériquement semi-simple régulière et soit $h_a:\ovl X\rta \car
\times^{\GG_m} L_D$ la section associée. Soit $U_a$ l'ouvert non
vide de $\ovl X$ où $a$ est une caractéristique semi-simple
régulière. On doit démontrer qu'il y a au-plus un nombre fini de
relèvements de $h_a$ en une section $h_b$ dans le diagramme
$$ \xymatrix{
& ((X_\Theta\times \mathbbm t)/\Sigma')\times^{\GG_m} L_D \ar[d]\\
 X \ar[ur]^{h_b} \ar[r]_{h_a\ \ \ \ } & \car\times^{\GG_m}L_D }
$$
D'après le lemme précédent, la donnée du relèvement $h_b$ est
équivalente à la donnée de $h_b$ sur l'ouvert $U_a$. La théorie de
Galois usuelle montre que ceci est possible si et seulement si
l'image $\Sigma'_a$ de la représentation de monodromie
$$\rho_a:\pi_1(U_a,u)\rta W'$$
soit conjugué à un sous-groupe de $\Sigma'$. De plus, dans ce cas,
il y a exactement $|{\rm Nor}_{W'}(\Sigma'_a)|/|\Sigma'_a|$
relèvements différents.

Démontrons que le morphisme $\pi^X_{\Sigma'}:{\mathbbm
B}^{\rmgssr}_{\Sigma'}\rta \bbA^{\rmgssr}$ est {\em propre}. Soit
$S$ un trait et soit $a\in\bbA^{\rmgssr}(S)$ un $S$-point de
$\bbA^{\rmgssr}$ donné par un morphisme $h_s:S\times X\rta
\car\times^{\GG_m}L_D$. Supposons qu'au-dessus du point générique
$\eta$ de $S$, $h_a(\eta)$ se relève en $h_b(\eta): \eta\times
X\rta ((X_\Theta\times \mathbbm t)/\Sigma')\times^{\GG_m} L_D$.
Alors $h_b(\eta)$ se prolonge de façon unique en une section
$h_b:S\times X\rta ((X_\Theta\times \mathbbm
t)/\Sigma')\times^{\GG_m} L_D$ relevant $h_a$ d'après le lemme
précédent.

Démontrons que le morphisme $\pi^X_{\Sigma'}:{\mathbbm
B}^{\rmgssr}_{\Sigma'}\rta \bbA^{\rmgssr}$ est {\em net}.
Donnons-nous deux flèches
$$h_{b}, h_{b'}:X\times \Spec(k[\varepsilon]/\varepsilon^2)\rta
M=((X_\Theta\times \mathbbm t)/\Sigma')\times^{\GG_m}L_D$$ qui
induisent la même flèche
$$h_a:X\times \Spec(k[\varepsilon]/\varepsilon^2)
\rta N=\car\times^{\GG_m}L_D.$$ Soit $U_a$ l'ouvert de $X$ image
réciproque de l'ouvert semi-simple régulier de $\car$ et supposons
que $U_a$ soit non-vide. Puisqu'au-dessus de $\car^{\rm ssr}$ le
morphisme $(X_\Theta\times \mathbbm t)/\Sigma' \rta \car$ est
étale, les deux flèches $h_b$ et $h_{b'}$ sont nécessairement
égales au-dessus de $U_a$. Ceci implique que $h_b=h_{b'}$ car la
$\calO_X$-Algèbre $\calO_X[\varepsilon]/(\varepsilon^2)$ est
plate.

Le morphisme $\pi^X_{\Sigma'}:{\mathbbm B}^{\rmgssr}_{\Sigma'}\rta
\bbA^{\rmgssr}$ fini et donc en particulier représentable, a
fortiori ${\mathbbm B}^{\rmgssr}_{\Sigma'}$ est représentable.
L'image de ${\mathbbm B}^{\rmgssr}_{\Sigma'}$ dans
$\bbA^{\rmgssr}$ est un sous-schéma fermé $\ovl{\bbA}_{[\Sigma']}$
de $\bbA^{\rmgssr}$ dont les points géométriques $a$ sont
caractérisés par le fait que la classe de $W'$-conjugaison
$[\Sigma'_a]$ contient un membre qui est un sous-groupe de
$\Sigma'$.

Si $\Sigma_1'$ est un sous-groupe de $\Sigma'$, on a l'inclusion
de sous-schémas fermés $\ovl{\bbA}_{[\Sigma'_1]}\subset
\ovl{\bbA}_{[\Sigma']}$. On a donc un ouvert de $\bbA_{[\Sigma']}$
de $\ovl\bbA_{\Sigma'}$ qui consiste en les points $a\in
\ovl\bbA_{[\Sigma']}(\ovl k)$ tels que $[\Sigma'_a]$ est la classe
de $W'$-conjugaison de $\Sigma'$. Son image réciproque dans
${\mathbbm B}^\rmgssr$ est l'ouvert $\mathbbm B_{\Sigma'}^{\rm
max}$. Au-dessus de $\bbA_{[\Sigma']}$, le morphisme $\mathbbm
B_{\Sigma'}^{\rm max} \rta \bbA_{[\Sigma']}$ est un morphisme
fini, net dont toutes les fibres géométriques ont $|{\rm
Nor}_{W'}(\Sigma')|/|\Sigma'|$ points. Dans un voisinage
hensélien, un morphisme fini et net est la réunion disjointe
d'immersions fermées. Si celui a de plus un nombre constant de
points dans ses fibres géométriques, il doit être réunion
disjointe d'isomorphismes. C'est donc un morphisme fini étale.
\findem

\bigskip
La conjonction de 7.1 et 7.2 implique qu'on a une stratification
de $\bbA^{\rmgssr}$
$${\bbA}^{\rmgssr}=\bigsqcup_{[\Sigma']}\bbA_{[\Sigma']}$$
où $[\Sigma']$ par court l'ensemble des classes de
$W'$-conjugaison de sous-groupes de $W'$ satisfaisant à la
contrainte 7.1.

\begin{corollaire}
Le faisceaux $\pi_0(\widetilde U/\bbA^\heartsuit)$ est localement
constant le long de la strate $\bbA_{[\Sigma']}$ de fibre typique
$W'/\Sigma'$. Le faisceau $\pi_0(P'/\bbA^{\heartsuit})$ est
localement constant le long de la strate $\bbA_{[\Sigma']}$ de
fibre typique $\bbX^\sv_{\Sigma'}$.
\end{corollaire}

\dem La proposition 7.2 montre que $\pi_0(\widetilde
U/\bbA^\heartsuit)$ est localement constant le long des strates
$\bbA_{[\Sigma']}$ de fibre typique $W'/\Sigma'$. La seconde
assertion s'en déduit en vertu de 6.8. \findem

\begin{definition}
Un point $a\in\bbA^{\rmgssr}(\bar k)$ est dit {\em elliptique} si
pour un membre arbitraire $\Sigma'_a$ de la classe de
$W'$-conjugaison $[\Sigma'_a]$, le groupe $\bbX^\vee_{\Sigma'_a}$
des $\Sigma'_a$-coinvariants des cocaractères est un groupe fini.
\end{definition}

Il est clair sur cette définition que l'ensemble des $a$
elliptiques est une réunion des strates $\bbA_{[\Sigma']}$.

\begin{corollaire}
Supposons que le groupe des caractères $\bbX$ de $\mathbbm T$ n'a
pas d'invariants non triviaux sous l'action de $W'$ ce qui est le
cas si $G$ est semi-simple. L'ensemble des caractéristiques
$a\in\bbA^{\rmgssr}(\ovl k)$ elliptiques forment un ouvert non
vide $\bbA^{\rm ell}$ de $\bbA^{\rmgssr}$.

Cet ouvert est caractérisé par la propriété $a\in \bbA^{\rm
ell}(\ovl k)$ si et seulement si $\pi_0(P_a)$ est fini.
\end{corollaire}

\dem D'après la relation d'adhérence décrite ci-dessus, si un
point $a\in\bbA^{\rmgssr}$ est elliptique, toutes ses
généralisations le sont. Il s'ensuit que la réunion des strates
$\bbA_{[\Sigma']}$ pour les classes de conjugaison des
sous-groupes $\Sigma'\subset W'$ tels que le groupe des
coinvariants $\bbX^\sv_{\Sigma'}$ est fini, forme un ouvert de
$\bbA$. Il reste à vérifier que celui-ci est non vide. D'après
\cite{Fa}, il existe des caractéristiques $a$ très régulières,
voir définition 4.1, pour lesquels le groupe de monodromie est
automatiquement tout le groupe $W'_\geom$. Celles-ci sont alors
clairement elliptiques.

Si $a\in\bbA^{\rm ell}(\ovl k)$, le groupe $\pi_0(P'_a)$ est fini.
La suite exacte 6.3 implique alors que $\pi_0(P_a)$ est aussi fini
car il est un quotient de $\pi_0(P'_a)$. Dans cette suite exacte
dont le terme $\rmH^0(\ovl X,J_a/J_a^0)$ est un groupe fini, la
finitude de $\pi_0(P_a)$ implique celle de $\pi_0(P'_a)$. D'après
le corollaire précédent, $\pi_0(P'_a)$ fini si et seulement si
$a\in\bbA^{\rm ell}(\ovl k)$.\findem

\bigskip

Pour $G$ un groupe classique semi-simple, on peut vérifier
qu'au-dessus de l'ouvert elliptique $\bbA^{\rm ell}$, l'espace de
module $\calM^{\rm ell}$ est un champ de Deligne-Mumford et le
morphisme $f^{\rm ell}:\calM^{\rm ell}\rta \bbA^{\rm ell}$ est un
morphisme propre de type fini en suivant les arguments de
\cite{Fa}. Dans le cas unitaire, c'est fait dans \cite{L-N}, nous
en donnerons une démonstration détaillée dans un travail futur
pour les autres groupes classiques.

\section{La $[\kappa]$-décomposition}

L'action de $P$ sur $\calM$ relativement à $\bbA$ induit une
action de $P$ sur la cohomologie de $\calM$ c'est-à-dire pour tout
$\bbA$-schéma $a:S\rta\bbA$, le groupe $P(S)$ agit sur la
restriction à $S$ de $f_*\QQ_\ell$. Considérons cette action
au-dessus de l'ouvert $\bbA^{\rmgssr}$. Sur cet ouvert, $P$ est
lisse de sorte qu'on peut définir $P^0$ sa composante neutre
relative en vertu de 6.1. D'après le lemme d'homotopie (\cf
\cite{L-N} 3.2), si $p$ est une section de $P^0$, $p$ agit
trivialement sur les faisceaux de cohomologie perverse
$^p\rmH^j(f_*\QQ_\ell)$. On en déduit une action du faisceau étale
$\pi_0(P/\bbA^{\rmgssr})$ sur les faisceaux pervers
$^p\rmH^j(f_*\QQ_\ell)$. Sur l'ouvert elliptique $\bbA^{\rm ell}$,
$\pi_0(P/\bbA^{\rm ell})$ est un faisceau en groupes abéliens
finis. Cette action induit localement une décomposition en somme
directe de $^p\rmH^j(f^{\rm ell}_*\QQ_\ell)$.

\begin{proposition} Soient $S$ un schéma de type fini, $A$ un faisceau
en groupes abéliens finis pour la topologie étale de $S$ et $K$ un
faisceau pervers sur $S$. Supposons que $A$ agit sur $K$. En
particulier, le groupe fini $\Gamma(S,A)$ agit sur $K$ qui induit
une décomposition
$$K=\bigoplus_{\kappa\in \Gamma(S,A)\rta\ovl \QQ_\ell^\times}
K_\kappa
$$
où $\kappa$ parcourt l'ensemble des caractères du groupe des
section globales de $S$. Pour tout caractère $\kappa$ de
$\Gamma(S,A)$, $K_\kappa$ est supporté par le fermé $S_\kappa$
constitué des point $\ovl s\in S$ tel que le caractères
$\kappa:\Gamma(S,A)\rta\ovl\QQ_\ell^\times$ se factorise par
l'homomorphisme de restriction aux fibres $\Gamma(S,A)\rta A_{\ovl
s}$.
\end{proposition}

\dem Soit $a\in \Gamma(S,A)$ une section telle que
$\kappa(a)\not=1$. Soit $U(a)$ l'ouvert de $S$ où la section $a$
s'annule. Pour démontrer que $K_\kappa$ est supporté par
$S_\kappa$ qui est le fermé complémentaire de la réunion des
ouverts $U(a)$ avec $\kappa(a)\not=1$, il suffit de démontrer que
$K_\kappa|_{U(a)}=0$ pour des tels $a$. Ceci est clair car $a$,
étant égal à la section neutre sur $U(a)$, agit trivialement sur
$K|_{U(a)}$. Pour tous caractères $\kappa$ de $\Gamma(S,A)$ tels
que $\kappa(a)\not=1$, on a $K_\kappa|_{U(a)}=0$. \findem

\bigskip Soit $a\in\bbA^{\rm ell}(k)$ une caractéristique définie
sur $k$, soit $S_a$ l'hensélisé de $\bbA^{\rm ell}$ en $a$. Soit
$\ovl a$ un point géométrique au-dessus de $a$. Alors, la fibre de
$\pi_0(P_{\ovl a})$ en le point $\ovl a$ est un groupe abélien
fini muni d'une action de $\Gal(\ovl k/k)=\langle \sigma\rangle$.

\begin{corollaire}
Sur $S_a$, on a une décomposition canonique
$$^p\rmH^j(f_{S_a,*}\QQ_\ell)=\bigoplus_{[\kappa]}\,
^p\rmH^j(f_{S_a,*}\QQ_\ell)_{[\kappa]}
$$
où $\kappa$ parcourt l'ensemble des classes de
$\sigma$-conjugaison de caractères de $\pi_0(P_{\ovl a})$ et où
géométriquement $^p\rmH^j(f_{S_a,*}\QQ_\ell)_{[\kappa]}$ est la
somme directe des facteurs isotypiques
$^p\rmH^j(f_{S_a,*}\QQ_\ell)_{\kappa'}$ avec $\kappa'$ dans la
classe de $\sigma$-conjugaison de $\kappa$.

De plus, $^p\rmH^j(f_{S_a,*}\QQ_\ell)_{[\kappa]}$ est supporté par
le sous-schéma fermé de $S_a$ constitué des points $\ovl s\in S$
tels qu'il existe au moins un caractère $\kappa'$ dans la classe
de $\sigma$-conjugaison de $\kappa$ tel que
$\kappa':\pi_0(P_a)^{\langle\sigma\rangle}\rta\ovl\QQ_\ell^\times$
se factorise par l'homomorphisme canonique
$\pi_0(P_a)^{\langle\sigma\rangle}\rta \pi_0(P_{\ovl s})$.
\end{corollaire}

\dem La somme directe $\bigoplus_{\kappa'\in [\kappa]}
^p\rmH^j(f_{S_a,*}\QQ_\ell)_{\kappa'}$ est stable par $\sigma$ et
défini un facteur direct de $^p\rmH^j(f_{S_a,*}\QQ_\ell)$.
L'énoncé sur le support se déduit du lemme précédent. \findem

\bigskip
Comme on verra dans §9, cette décomposition au-dessus d'un
voisinage d'un point rationnel de $\bbA^{\rm ell}$ correspond à la
réécriture habituelle de la somme des intégrales orbitales
globales dans une classe de conjugaison stable donnée, en une
somme de $\kappa$-intégrales orbitales à l'aide d'une
transformation de Fourier sur le groupe fini des obstructions. On
peut aussi considérer la décompo\-sition globale sur l'ouvert
$\bbA^{\rm ell}$ qui correspond alors à la préstabilisation de la
partie elliptique de la formule des traces. Dans cette situation
globale, il nous faudra introduire le formalisme de cofaisceaux,
évident mais inhabituel.

Revenons à la situation générale où on a au-dessus d'une base $S$
un faisceau en groupes abéliens finis $A$ agissant sur un faisceau
pervers $K$. Pour tout ouvert étale $U$ de $S$, considérons le
groupe $A(U)^*$ des caractères d'ordre fini de $A(U)$. Si $U'$ est
un ouvert étale de $U$, la flèche de restriction $A(U)\rta A(U')$
induit une flèche de corestriction $A(U')^*\rta A(U)^*$. On
appelle {\em précofaisceau} sur $S$ un foncteur covariant de la
catégorie des ouverts étales de $S$ dans la catégorie des
ensembles. En particulier, $U\mapsto A(U)^*$ est un précofaisceau.

On appelle {\em cofaisceau} un précofaisceau ${\rm Co}$ qui
vérifie en plus une propriété de recollement : pour tout ouvert
étale $U$ de $S$, pour tout recouvrement étale $U_1$ de $U$, et
pour $U_2=U_1\times_U U_1$, l'ensemble ${\rm Co}(U)$ est
l'ensemble quotient de ${\rm Co}(U_1)$ par la relation
d'équivalence ${\rm Co}(U_2)\rightrightarrows {\rm Co}(U_1)$.

Tout comme pour les faisceaux, on peut associer un cofaisceau
${\rm Co}$ à un précofaisceau ${\rm Preco}$. Pour tout ouvert
étale $U$ de $S$, on prend pour ${\rm Co}(U)$ la limite projective
sur les recouvrements étales $U_1\rta U$ des ensembles des classes
d'équivalence de ${\rm Preco}(U_1)$ par le relation d'équivalence
${\rm Preco}(U_2)\rightrightarrows {\rm Preco}(U_1)$ où
$U_2=U_1\times_U U_1$.

L'énoncé suivant est tautologique. Nous allons néanmoins en
esquisser une démonstration car la notion du cofaisceau est moins
familier que celle du faisceau.

\begin{proposition}
Soit $A$ un faisceau en groupes abéliens finis sur $S$ qui agit
sur un faisceau pervers $K$ sur $S$. Soit $A^{\rm co}$ le
cofaisceau associé au précofaisceau $U\mapsto A(U)^*$ où $A(U)^*$
est l'ensemble des caractères d'ordre fini de $A(U)$. Alors pour
tout ouvert étale $U$ de $S$, on a un décomposition de la
restriction de $K$ à $U$
$$K|_U=\bigoplus_{[\alpha]\in A^{\rm co}(U)}(K|_U)_{[\alpha]}$$
et pour tout $[\alpha]\in A^{\rm co}(U)$ et pour tout ouvert étale
$U'$ de $U$, on a la relation de compatibilité
$$((K|_U)_{[\alpha]})|_{U'}=\bigoplus_{[\alpha']\in A^{\rm co}(U')\atop
[\alpha']\mapsto [\alpha]} (K|_{U'})_{[\alpha']}.$$
\end{proposition}

\dem Il suffit de démontrer que pour tout recouvrement étale $U_1$
de $U$, on a une décomposition canonique
$$K|_U=\bigoplus_{[\alpha]}(K|_U)_{[\alpha]}$$
où $[\alpha]$ parcourt l'ensemble des classes d'équivalence de
$A(U_1)_*$ pour la relation d'équivalence
$A(U_2)^*\rightrightarrows A(U_1)^*$.

En effet, on a une décomposition
$$K|_{U_1}=\bigoplus_{\alpha\in A(U_1)^*}(K|_{U_1})_\alpha$$
à cause de l'action du groupe abélien $A(U_1)$ sur $K|_{U_1}$.
Regroupons les facteurs directs $(K|_{U_1})_\alpha$ selon la
relation d'équivalence $A(U_2)^*\rightrightarrows A(U_1)^*$, la
somme directe
$$\bigoplus_{\alpha\in [\alpha]}(K|_{U_1})_\alpha$$
descend alors à $U$. \findem

\bigskip

Pour le cas particulier du faisceau $\pi_0(P'/\bbA^{\rm ell})$, on
a des renseignements explicites sur les sections globales de
$\pi_0(P'/\bbA^{\rm ell})^{\rm co}$ par la proposition 6.8.

\begin{proposition}
On a une application canonique de l'ensemble des section globales
du cofaisceau $\pi_0(P'/\bbA^{\rm ell})^{\rm co}$, dans l'ensemble
des classes de $W'$-conjugaison $[\kappa]$ des caractères d'ordre
fini $\kappa:\bbX^\sv \rta\QQ_\ell^\times$.
\end{proposition}

\dem Choisissons un recouvrement étale $S_1\rta \bbA^{\rm ell}$
au-dessus du quel $\pi_0(\widetilde U/S_1)$ admet une section
qu'on va noter $\beta$. Grâce à cette section, par la proposition
6.8, on a une flèche surjective
$$\bbX^\sv_{S_1} \rta \pi_0(P'/S_1).$$
On a ainsi une application de $\pi_0(P'/S_1)$ dans l'ensemble des
caractères d'ordre fini $\kappa:\bbX^\sv \rta\QQ_\ell^\times$.

Soient $S_2=S_1\times_{\bbA^\heartsuit} S_1$ et $\pr_1(\beta)$,
$\pr_2(\beta)$ les images réciproques de $\beta$ sur $S_2$.
Puisque $W'$ agit transitivement sur les fibres de
$\pi_0(\widetilde U/\bbA^{\rm ell})$, il existe une section $w'\in
\Gamma(S_2,W')$ tel que $w'\pr_1(\beta)=\pr_2(\beta)$. On en
déduit une application de l'ensemble des section globales du
cofaisceau $\pi_0(P'/\bbA^{\rm ell})^{\rm co}$, dans l'ensemble
des classes de $W'$-conjugaison $[\kappa]$ des caractères d'ordre
fini $\kappa:\bbX^\sv \rta\QQ_\ell^\times$. Cette application ne
dépend ni du choix de $U_1$ ni du choix de la section $\beta$.
\findem

\bigskip
On déduit de 8.3 et 8.4 une décomposition
$$^p\rmH^j(f_*^{\rm ell}\QQ_\ell)=\bigoplus_{[\kappa]}\,^p
\rmH^j(f_*^{\rm ell}\QQ_\ell)_{[\kappa]}
$$
où la somme directe est étendue sur l'ensemble des classes de
$W'$-conjugaison des éléments d'ordre fini de $\hat T$.

\begin{theoreme}
Supposons que la caractéristique de $k$ ne divise pas $|W|$.
Supposons que le groupe des caractères $\bbX$ de $\mathbbm T$ n'a
pas d'invariants non triviaux sous $W'$. C'est notamment le cas si
$G$ est semi-simple.

Sur l'ouvert $\bbA^{\rm ell}$ de $\bbA$, qui est non vide sous les
hypothèses précédentes, on a une décomposition canonique
$$^p\rmH^j(f_*^{\rm ell}\QQ_\ell)=\bigoplus_{[\kappa]}\,^p
\rmH^j(f_*^{\rm ell}\QQ_\ell)_{[\kappa]}
$$
où la somme directe est étendue sur l'ensemble des classes de
$W'$-conjugaison des éléments d'ordre fini de $\hat T$.

Soit $\kappa$ un élément d'ordre fini de $\hat T$. Le facteur
direct $\rmH^j(f_*^{\rm ell}\QQ_\ell)_{[\kappa]}$ est supporté par
le fermé $\ovl {\mathbb S}_{[\kappa]}$ constitué des points $ a\in
\bbA^{\rm ell}(\ovl k)$ tel que pour tout point géométrique
$\widetilde u \in \widetilde U_a$, l'un des caractères $\kappa$ de
la $W'$-classe $[\kappa]$ se factorise à travers le composé de
l'homomorphisme
$$\pi_0(\lambda(\widetilde u)):\bbX^\sv\rta\pi_0(P'_{a})$$
défini dans la démonstration de 6.8 et l'homomorphisme
$\pi_0(P'_{a})\rta \pi_0(P_{a})$.

En particulier,
$$\ovl {\mathbb S}_{[\kappa]}\subset \bigsqcup_{[\Sigma']}\bbA_{[\Sigma']}$$
où $[\Sigma']$ parcourt l'ensemble des classes de $W'$-conjugaison
de sous-groupes de $W'$ contenant au moins un membre $\Sigma'$ qui
stabilise $\kappa$ et où $\bbA_{[\Sigma']}$ est strate associée à
$[\Sigma']$ comme dans la section 7.
\end{theoreme}

\dem La première assertion résulte de 8.3 et 8.4. La seconde
assertion résulte de 8.1. La troisième assertion résulte de la
seconde, de la propositions 6.8 et du corollaire 7.4. \findem

\begin{corollaire} La composante indexée $[\kappa]$, avec $\kappa\in\hat
T$, est non nulle seulement si la $W$-orbite de $\kappa$ est
stable par $\Theta_\geom$.
\end{corollaire}

\dem D'après l'estimation du support, si celle-ci est non vide, il
existe un sous-groupe $\Sigma'$ de $W'=W\rtimes \Theta$ d'image
$\Theta_\geom$ dans $\Theta$ qui fixe $\kappa$. Il revient au même
de dire que la $W$-orbite de $\kappa$ est stable par
$\Theta_\geom$. \findem

\bigskip
L'énoncé de support du théorème 8.6 va être raffiné dans la
section §10 une fois qu'on aura introduit les groupes
endoscopiques. Disons pour l'instant que le théorème 8.6 est déjà
optimal dans le cas où $G$ est un groupe réductif de centre
connexe comme le groupe unitaire où les groupes semi-simple
adjoints. Dans ce cas, $P'_a=P_a$ car $J_a$ a des fibres connexes.

\section{Le lien avec les $\kappa$-intégrales orbitales}

La décomposition du corollaire 8.2 correspond bien à la
décomposition de la somme des intégrales locales dans une classe
de conjugaison stable, en somme des $\kappa$-intégrales orbitales.
Soit $a\in\bbA^{\rm ell}(k)$ une caractéristique elliptique à
valeur dans $k$. D'après §2, le nombre de $k$-points de la fibre
$\calM_a$ est donné par la formule
$$|\calM_a(k)|=\sum_{\alpha\in\ker^1(F,G)} \sum_{\gamma\in\frakg
(\alpha)/{\rm conj.}\atop \chi(\gamma)=a} O_\gamma(1_D).
$$
D'après Langlands et Kottwitz (\cf \cite{RL}, \cite{K-EST}) cette
somme peut être transformée en une somme de $\kappa$-intégrales
orbitales. Nous allons passer en revue cette transformation de
notre point de vue.

Considérons le morphisme quotient
$$\calM_a \rta [\calM_a/P_a]=\prod_v
[\calM^\bullet_{a,v}/P^\bullet_{a,v}]$$ dont le quotient se
décompose en produit des facteurs locaux d'après 4.5. Au-dessus de
chaque $k$-point ${\underline m}=(m_v)$ de $\prod_v
[\calM^\bullet_{a,v}/P^\bullet_{a,v}]$, on a un $P_a$-torseur dont
la classe d'isomorphisme définit un élément ${\rm inv}({\underline
m})\in \rmH^1(k,P_a)$. Cette classe s'annule si et seulement si le
torseur a un $k$-point. Notons que la suite exacte longue associée
à la suite exacte courte
$$1\rta P_a^0\rta P_a \rta \pi_0(P_a)\rta 1$$
induit un isomorphisme
$$0=\rmH^1(k,P_a^0) \rta \rmH^1(k,P_a) \isom \rmH^1(k,\pi_0(P_a))
\rta \rmH^2(k,P_a^0)=0,
$$
en vertu des annulations résultant du théorème de Lang et de la
dimension cohomologique d'un corps fini \cite{Se}. Soit $\ovl a$
un point géométrique au-dessus de $a$. Le groupe $\rmH^1(k,P_a)$
est isomorphe au groupe des $\sigma$-coinvariants de
$\pi_0(P_{\ovl a})$
$$\rmH^1(k,P_a)=\pi_0(P_{\ovl a})_{\langle\sigma\rangle}$$
et donc en particulier, fini.

On peut maintenant compter d'une autre façon le nombre de
$k$-points de $\calM_a$
$$|\calM_a(k)|=|P_a(k)|\sum_{{\underline m}\in \prod_v
[\calM^\bullet_{a,v}/P^\bullet_{a,v}](k)\atop {\rm
inv}({\underline m})=0}{1\over {\rm Aut}({\underline m})(k)}.
$$
En faisant opérer la transformation de Fourier sur le groupe fini
$\pi_0(P_{\ovl a})_{\langle \sigma\rangle}$, on peut réécrire
cette somme à l'aide des $\kappa$-intégrales orbitales
$$|\calM_a(k)|=|P^0_a(k)| \sum_{\kappa:\pi_0(P_a)
 \rta \ovl\QQ_\ell^\times \atop \sigma(\kappa)=\kappa}
\sum_{{\underline m}\in \prod_v
[\calM^\bullet_{a,v}/P^\bullet_{a,v}](k)} \langle\kappa,{\rm
inv}({\underline m}) \rangle {1\over {\rm Aut}({\underline m})(k)}
$$
car $|P_a(k)/P^0_a(k)|=|\pi_0(P_{\ovl a})^{\langle
\sigma\rangle}|=|\pi_0(P_{\ovl a})_{\langle \sigma\rangle}|$. On a
donc
$$|\calM_a(k)|=\sum_{\kappa:\pi_0(P_a)
 \rta \ovl\QQ_\ell^\times \atop \sigma(\kappa)=\kappa}
 O^\kappa_a(1_D)$$
 avec
$$O^\kappa_a(1_D)=|P^0(k)| \sum_{{\underline m}\in \prod_v
[\calM^\bullet_{a,v}/P^\bullet_{a,v}](k)} \langle\kappa,{\rm
inv}({\underline m}) \rangle {1\over {\rm Aut}({\underline
m})(k)}.
$$
Notons que le choix d'un $k$-point de $\calM_a$ est implicite dans
cette définition de $O^\kappa_a(1_D)$ car les
$\calM^\bullet_{a,v}$ dépendent de ce point base. Il est commode
de prendre pour point base le point de Kostant construit dans 2.5
sous l'hypothèse que $G$ est quasi-déployé et qu'on s'est donné
une racine carrée de $L_D$.

Du fait que pour $\underline m=(m_v)$, $\inv(\underline m)$ est
une somme des invariants locaux ${\rm inv}(m_v)$,
$O^\kappa_a(1_D)$ se décompose en produit de $\kappa$-intégrales
orbitales locales
$$O^\kappa_a(1_D)=|P^0(k)| \prod_{v\in |X|}
O^\kappa_{a,v}(1_{D_v}).$$

Cette décomposition du nombre de $k$-points de $\calM_a$ est
essentiellement la même que la $[\kappa]$-décomposition de la
proposition 8.2 en vertu de l'énoncé suivant.

\begin{proposition}
Pour tous $\kappa:\pi_0(P_{\ovl a})\rta\ovl\QQ_\ell^\times$, si
$\sigma(\kappa)\not=\kappa$, on a
$${\rm Tr}(\sigma,^p\rmH^j(f_*\QQ_\ell)_{[\kappa],\ovl a})=0$$
pour tous $j$. Si $\sigma(\kappa)=\kappa$ ou autrement dit si la
$\sigma$-orbite $[\kappa]$ est un singleton, on a
$${\rm Tr}(\sigma,\sum_j (-1)^j[^p\rmH^j(f_*\QQ_\ell)_{[\kappa],\ovl a}])
=O^\kappa_a(1_D).$$
\end{proposition}

\dem Si $\sigma(\kappa)\not=\kappa$, $\sigma_a$ permute
circulairement les facteurs directs
$^p\rmH^j(f_*\QQ_\ell)_{\kappa',\ovl a}$ avec $\kappa'$ dans la
$\sigma$-orbite $[\kappa]$. Un tel opérateur a nécessairement une
trace nulle. Le cas $\sigma(\kappa)=\kappa$ résulte l'appendice
A.3 de \cite{L-N}. \findem

\bigskip
On en déduit le pendant arithmétique du corollaire 8.6.

\begin{corollaire}
Soit $\kappa$ un caractère d'ordre fini de $\mathbbm X^\sv$ dont
la $W$-orbite est stable par $\Theta_\geom$. Soit
$a\in\bbA^\heartsuit(k)$. Pour que
$${\rm Tr}(\sigma,^p\rmH^j(f_*\QQ_\ell)_{[\kappa],\ovl a})\not=0
$$
il est nécessaire que la $W$-orbite de $\kappa$ soit stable par
$\Theta$.
\end{corollaire}

\section{Le lien avec les groupes endoscopiques}

On suppose toujours que le groupe des caractères $\bbX$ de
$\mathbbm T$ n'a pas de $W'_\geom$-invariants non triviaux ce qui
est notamment le cas si $G$ est un groupe semi-simple. Cette
hypothèse implique en particulier l'existence des
$a\in\bbA^\heartsuit(\ovl k)$ elliptiques.

Soit $\hat G$ le groupe dual connexe complexe de $\GG$. Il est
muni d'un épinglage comprenant en particulier d'un tore maximal
$\hat T$ et d'un sous-groupe de Borel $\hat B$ comprenant $\hat
T$. On pose $^L G_\geom=\hat G\rtimes \Theta_\geom$ où
$\Theta_\geom$ est l'image de l'homomorphisme $\rho_G:\pi_1(\ovl
X,x)\rta \Out(\GG)$ associé au torseur $\tau_G^\Out$ et où le
produit semi-direct est construit à partir de l'action de
$\Theta_\geom$ sur $\hat G$ qui fixe l'épinglage.

Soit $\kappa$ un élément d'ordre fini de $\hat T$. Notons $\hat H$
la composante neutre du centralisateur $(\,^L G)_\kappa$ de
$\kappa$ dans $\,^L G$. C'est le sous-groupe réductif connexe de
$\hat G$ engendré par $\hat T$ et les sous-groupes radiciels
$X_{\alpha^\sv}$ de $\hat G$ associées aux racines
$\alpha^\sv:\hat T\rta \GG_m$ telles que $\alpha^\sv(\kappa)=1$.
Le groupe réductif connexe $\hat H$ vient avec un épinglage déduit
de celui de $\hat G$.

Notons $W_{H}$ le groupe de Weyl de $\hat T$ dans $\hat H$ et
$(W'_\geom)_\kappa$ le sous-groupe des éléments de
$W'_\geom=W\rtimes \Theta_\geom$ qui fixent $\kappa$. On a
$W_H=N(\hat T,\hat H)/\hat T$ et $(W'_\geom)_\kappa= N(\hat
T,(\,^L G_\geom)_\kappa)/\hat T$ où $N(\hat T,\hat H)$ et $N(\hat
T,(\,^L G_\geom)_\kappa)$ sont les normalisateurs de $\hat T$ dans
$\hat H$ et dans $(\,^L G)_\kappa$. On en déduit un isomorphisme
$N(\hat T,(\,^L G_\geom)_\kappa)/N(\hat T,\hat H) \rta
(W'_\geom)_\kappa/W_H$ et donc un homomorphisme
$$
(W'_\geom)_\kappa/W_H \rta (\,^L G_\geom)_\kappa/\hat H.$$

\begin{lemme}
\begin{enumerate}
\item
L'homomorphisme $(W'_\geom)_\kappa/W_H \rta (\,^L G)_\kappa/\hat
H$ est un isomorphisme.
\item
Soit ${(W'_\geom)}_\kappa^{\rm out}$ le sous-groupe de
$(W'_\geom)_\kappa$ des éléments qui laissent stable l'ensemble
des racines positives de $\hat H$. Alors
$$(W'_\geom)_\kappa=W_H\rtimes (W'_\geom)_\kappa^{\rm out}.$$
\end{enumerate}
\end{lemme}

\dem
\begin{enumerate}
\item Il revient au même de démontrer que l'homomorphisme
$$N(\hat
T,(\,^L G_\geom)_\kappa)/N(\hat T,\hat H)\rta (\,^L
G_\geom)_\kappa/\hat H
$$
est un isomorphisme. Il est injectif parce que $N(\hat T,\hat H)$
est l'intersection $N(\hat T,(\,^L G_\geom)_\kappa)\cap \hat H$.
Il est surjectif car tout automorphisme de $\hat H$ peut être
modifié par un homomorphisme intérieur pour faire un automorphisme
qui fixe l'épinglage de $\hat H$ et stabilise en particulier $\hat
T$.

\item
Il suffit de démontrer tout élément $w\in (W'_\geom)_\kappa$
s'écrit de façon unique sous la forme $w_\kappa=w_H w_\kappa^{\rm
out}$ avec $w_H\in W_H$ et $w_\kappa^{\rm out}\in
(W'_\geom)_\kappa^{\rm out}$. En effet, $w_\kappa$ agit sur
l'ensemble des racines de $\hat H$ et on peut modifier cette
action par un unique élément $w_H$ pour qu'il stabilise l'ensemble
des racines positives de $\hat H$. \findem
\end{enumerate}

\begin{definition} Une {\bf donnée endoscopique}\footnote{Notre définition
diffère légèrement de la définition de \cite{K-CTT} 7.1.}
géomé\-tri\-que non ramifiée, consiste en un couple
$(\kappa,\rho)$ constitué d'un élément
$$\rho:\pi_1(\ovl X,x)\rta (W'_\geom)_\kappa^{\rm out}$$
tel que le composé
$$\pi_1(\ovl X,x)\rta (W'_\geom)_\kappa^{\rm out}\simeq
(W'_\geom)_\kappa/W_H \rta W'_\geom/W=\Theta_\geom
$$
est l'homomorphisme $\rho_G^\geom$ de la section 5.
\end{definition}

Ceci implique en particulier que le sous-groupe
$(W'_\geom)_\kappa$ des éléments de $W'_\geom$ qui fixent
$\kappa$, se surjecte sur $\Theta_\geom$ autrement dit la
$W$-orbite de $\kappa$ est stable par $\Theta_\geom$. Donc
$\kappa$ vérifie la contrainte 8.6.

Soit $(\kappa,\rho)$ une donnée endoscopique géométrique non
ramifiée comme ci-dessus. Soit $\mathbbm H$ le groupe déployé sur
$k$ dont le dual complexe est $\hat H$. Notons $\rho_H$
l'homomorphisme composé
$$\rho_H:\pi_1(\ovl X,x)\rta (W'_\geom)_\kappa^{\rm out}
\simeq (W'_\geom)_\kappa/W_H\rta \Out(\hat H)=\Out(\mathbbm H).
$$
Notons $\Theta_{H}$ l'image de l'homomorphisme $\rho_H:\pi_1(\ovl
X,x)\rta (W'_\geom)_\kappa^{\rm out}$ et posons $W'_H=W_H\rtimes
\Theta_H$. On a alors les inclusions
$$W'_H=W_H\rtimes
\Theta_H \subset W_H\rtimes (W'_\geom)_\kappa^{\rm
out}=(W'_\geom)_\kappa\subset W'_\geom.$$

L'homomorphisme $\rho_H:\pi_1(\ovl X,x)\rta{\rm Out}(\mathbbm H)$
nous fournit un $\Out(\mathbbm H)$-torseur $\tau^\Out_H$ sur $\ovl
X$ et donc un schéma en groupes quasi-déployé $H$ sur $\ovl X$
obtenu en tordant $\mathbbm H$ muni de l'action de $\Out(\mathbbm
H)$ fixant un épinglage, par le $\Out(\mathbbm H)$-torseur
$\tau^\Out_H$.

Soit $\ovl X_{\Theta_H}$ le revêtement fini étale de $\ovl X$
pointé par un point géométrique $x_{\Theta_H}$ au-dessus du point
$x$, attaché à l'homomorphisme $\rho_H:\pi_1(\ovl X,x)\rta
\Theta_H$. Du fait que le diagramme
$$
\xymatrix{
  \pi_1(\ovl X,x) \ar[dr]_{\,\,\rho_G^\geom} \ar[r]^{\,\,\,\,\,\,\rho_H}
                & \Theta_H \ar[d]^{}  \\
                & \Theta_\geom             }
$$
est commutatif, il existe un $X$-morphisme $\ovl X_{\Theta_H}\rta
\ovl X_\Theta^\geom$ envoyant le point $x_{\Theta_H}$ sur le point
$x_\Theta^\geom$ de $X_\Theta^\geom$.

On a l'espace de module de Hitchin $\calM_H$ associé à $H$ et au
diviseur $D$ ainsi que le morphisme de Hitchin
$$f_H:\calM_H\rta \bbA_H$$
où $\bbA_H$ est l'espace des sections
$$\ovl X\rta ((\ovl X_{\Theta_H}\times \mathbbm t)/W'_H)\times^{\GG_m}L_D.$$
Le morphisme évident
$$(\ovl X_{\Theta_H}\times \mathbbm t)/W'_H
\rta (\ovl X_\Theta^\geom\times\mathbbm t)/W'_\geom
$$
induit un morphisme de l'espace affine de Hitchin pour $H$ vers
l'espace affine de Hitchin pour $G$
$$\bbA_H\rta \ovl \bbA= \bbA\otimes_k \ovl k.$$

Le cas où $\hat G$ a un groupe dérivé simplement connexe est très
agréable. Dans ce cas, le centralisateur $\hat G_\kappa$ est
connexe de sorte que $\hat H=\hat G_\kappa$. Il s'ensuit que
$\Theta_H=\Theta$ et les revêtements étales $X_{\Theta_H}$ et
$X_\Theta$ de $X$ sont les mêmes. Dans ce cas particulier, le
morphisme $\bbA_H\rta\ovl \bbA$ est l'un des morphismes $\mathbbm
B_{\Sigma'} \rta\bbA$ déjà considérés dans la section 7.

\begin{proposition}
Dans l'ouvert $\ovl\bbA^{\rmgssr}$, le morphisme ${\mathbbm
A}_{H}^{\rmgssr}\rta \ovl\bbA^{\rmgssr}$ est un morphisme fini et
net. Soit $\bbA_{H}^{\rm max}$ l'ouvert de $\bbA_{H}$ où la
monodromie de $J_{H,a}$ est maximale en l'occurrence est égale à
$W'_H$. La restriction de ${\mathbbm A}_{H}^{\rmgssr}\rta \ovl
\bbA^{\rmgssr}$ à $\bbA_{H}^{\rm max}$ est finie et étale sur son
image.

Soient $\rho_1$ et $\rho_2$ deux homomorphismes $\pi_1(\ovl X,\ovl
x)\rta (W'_\geom)_\kappa^{\rm out}$ qui ne sont pas conjuguées.
Alors les images de ${\mathbbm A}_{H(\rho_1)}^{\rmgssr}$ et
${\mathbbm A}_{H(\rho_2)}^{\rmgssr}$ dans $\bbA^\heartsuit$ sont
des sous-schémas localement fermés disjoints.
\end{proposition}

\dem L'assertion sur la finitude et la netteté se démontre comme
7.2. L'assertion sur la disjonction des images pour de différentes
$\rho$ se démontre comme suit. Soit $a\in\bbA^\heartsuit(\ovl k)$
un point géométrique dans l'image de $\bbA_H\rta\bbA$. Soit $u$ un
point de l'ouvert $U_a$ où $J_a$ est un tore. La monodromie de
$J_a|_{U_a}$ est donné par un homomorphisme $\pi_1(U_a,u_a)\rta
(W'_\geom)_\kappa$ lequel passe au quotient pour déterminer
$\rho:\pi_1(\ovl X,x)\rta (W'_\geom)_\kappa/W_H$ à conjugaison
près. Ainsi la monodromie du tore $J_a|_{U_a}$ détermine $\rho$ à
conjugaison près. \findem

\begin{theoreme}
Supposons que la caractéristique de $k$ ne divise pas $|W|$.
Supposons que $\mathbbm X^\sv$ n'a pas d'invariants sous
$W'_\geom$.

Soient $\kappa$ un élément d'ordre fini de $\hat T$. Le facteur
direct $^p\rmH^j(f^{\rm ell}_*\QQ_\ell)_{[\kappa]}$ est supporté
par la réunion des images des $\bbA_{H_\rho}$ pour différentes
données endoscopiques géométriques non ramifiées $(\kappa,\rho)$.
En particulier, si la $W$-orbite de $\kappa$ n'est pas stable sous
$\Theta_\geom$, ce facteur direct est nul.
\end{theoreme}

\dem Soit $a$ un point géométrique de $\bbA^{\rm ell}$ tel que la
fibre en $a$ de l'un des faisceau pervers $^p\rmH^j(f^{\rm
ell}_*\QQ_\ell)_{[\kappa]}$ est non nulle. Le point $a$ correspond
à un morphisme $h_a:\ovl X\rta [\car/\GG_m]$. Notons $U_a$ l'image
inverse de $[\car^{\rm ssr}/\GG_m]$ et $u$ un point géométrique de
$U_a$. Soit $\widetilde U_a$ le revêtement étale galoisien de
$U_a$ de groupe de Galois $W'_\geom$ comme dans 6.8. Soit
$\widetilde u$ un point géométrique de $\widetilde U_a$ au-dessus
de $u$. Ce choix donne un homomorphisme $\rho_a:\pi_1(U_a,u_a)\rta
W'_\geom$ d'image $\Sigma'_a$.

D'après le théorème 8.4, pour que la fibre en $a$ de l'un des
faisceau pervers $^p\rmH^j(f^{\rm ell}_*\QQ_\ell)_{\kappa}$ soit
non nulle, il est nécessaire que l'un des conjugués de $\Sigma'_a$
soit contenu dans $(W'_\geom)_\kappa$ le sous-groupe de $W'_\geom$
fixateur de $\kappa$. On peut supposer que $\Sigma'_a\subset
(W'_\geom)_\kappa$ quitte à changer $\widetilde u$.

Rappelons que l'homomorphisme $\pi_1(U_a,u_a)\rta \pi_1(\ovl
X,u_a)$ est surjectif et a un noyau engendré par les sous-groupes
d'inertie $I_x$ pour les points $x\in \ovl X-U_a$. Pour démontrer
que $a$ provient d'un espace de Hitchin $\bbA_{H_\rho}$ défini
précédemment, il suffit de démontrer que l'homomorphisme composé
de $\rho_a:\pi_1(U_a,u)\rta (W'_\geom)_\kappa$ et
$(W'_\geom)_\kappa\rta (W'_\geom)_\kappa/W_H$ se factorise par
$\pi_1(\ovl X,u)$. Il revient au même de démontrer que pour tout
$x\in \ovl X-U_a$, la restriction de cet homomorphisme aux
sous-groupes d'inertie $I_x$ est triviale.

Le schéma en groupes commutatif lisse $J_a$ sur $\ovl X$ admet un
sous-schéma en groupes des composantes neutres $J_a^0$ : pour tout
point géométrique $x\in \ovl X$, la fibre en $x$ de $J_a^0$ est la
composante neutre de la fibre en $x$ de $J_a$. On a une suite
exacte de faisceaux
$$0\rta J_a^0\rta J_a \rta J_a/J_a^0\rta 0$$
dont le conoyau $J_a/J_a^0$ est supporté par $\ovl X-U_a$. En
effet, au-dessus de $U_a$, $J_a$ est un tore et en particulier a
des fibres connexes.

Le faisceau $J_a/J_a^0$ étant supporté par un schéma de dimension
nulle, on en déduit une suite exacte
$$\rmH^0(\ovl X,J_a/J_a^0)\rta \rmH^1(\ovl X,J_a^0)\rta
\rmH^1(\ovl X,J_a)\rta 0.
$$
On en déduit la suite exacte à droite
$$\rmH^0(\ovl X,J_a/J_a^0)\rta \pi_0(\rmH^1(\ovl X,J_a^0))\rta
\pi_0(\rmH^1(\ovl X,J_a))\rta 0.
$$
Puisque $J_a^0$ est un schéma en groupes lisse à fibres connexes
et génériquement torique, on a un isomorphisme
$$\pi_0(\rmH^1(\ovl X,J_a^0)) \isom \bbX^\sv_{\Sigma'_a}$$
de groupe des composante connexes du champ des $J_a^0$-torseurs
est le groupe des coinvariants de $\pi_1(U_a,u_a)$ dans le groupe
des cocaractères $\bbX^\sv$. L'élément $\kappa\in\hat T$ définit
donc un caractère
$$\kappa:\pi_0(\rmH^1(\ovl X,J_a^0)) \rta \QQ_\ell^\times.$$
L'hypothèse que la fibre de l'un des faisceaux pervers
$^p\rmH^j(f^{\rm ell}_*\QQ_\ell)_{\kappa}$ est non nulle, implique
que $\kappa$ se factorise à travers $\pi_0(\rmH^1(\ovl X,J_a))$
autrement dit la restriction de $\kappa$ au groupe fini
$\rmH^0(\ovl X,J_a/J_a^0)$ est nulle.

Soit $x$ un point dans $\ovl X-U_a$. On a un diagramme commutatif
$$
\xymatrix{
  \pi_0((J_a)_x) \ar[d]_{} \ar[r]^{} & \pi_0(J_a(F_x)/J_a^0(\calO_x))
  \ar[d]_{} \ar[r]^{} & \pi_0(J_a(F_x)/J_a(\calO_x)) \ar[d]^{} \\
  \rmH^0(\ovl X,J_a/J_a^0) \ar[r]^{} & \pi_0(\rmH^1(\ovl X,J_a^0)) \ar[r]^{}
  & \pi_0(\rmH^1(\ovl X,J_a))   }
$$
où $\pi_0(J_a(F_x)/J_a^0(\calO_x))=\bbX^\sv_{I_x}$. L'élément
$\kappa$ peut être vu comme un caractère
$$\kappa:\pi_0(J_a(F_x)/J_a^0(\calO_x))\rta \QQ_\ell^\times$$
dont la restriction à $\pi_0((J_a)_x)$ est nulle. Autrement dit
$\kappa$ est dans l'image de l'application duale
$$\pi_0(J_a(F_x)/J_a(\calO_x))^* \rta
\pi_0(J_a(F_x)/J^0_a(\calO_x))^*=\hat T^{I_x}.$$

Il est maintenant nécessaire de prendre un groupe auxiliaire $G_1$
comme dans \cite[7.5]{K-EST}. La construction de loc. cit. fournit
une suite exacte de schémas en groupes réductifs connexes sur $X$
$$1\rta G\rta G_1\rta C\rta 1$$
vérifiant les propriétés suivantes :
\begin{itemize}
\item $C$ est un tore sur $X$,
\item $G_1$ est un groupe réductif de centre connexe,
\item son dual complexe $\hat G_1$ a un groupe dérivé simplement
connexe.
\end{itemize}
On renvoie à loc. cit. pour la construction de $G_1$. Disons
seulement que dans le cas $G$ semi-simple, l'argument qui suit,
avec quelques modifications mineures, marchera aussi  si on prend
$G_1$ le groupe adjoint de $G$.

On peut définir sur $\car_{G_1}=\frakt_1/W$ le centralisateur
régulier $J_1$ de $G_1$. Puisque $G_1$ est un groupe réductif
connexe à centre connexe, $J_1$ est alors un schéma en groupes
commutatif, lisse et {\em à fibres connexes}. On a de plus un
homomorphisme de l'image réciproque de $J$ sur $\car_{G_1}$ dans
$J_1$, qui se déduit de l'homomorphisme $G\rta G_1$. En retirant
sur $X$ par l'homomorphisme $h_a$, on a une suite exacte de
schémas en groupes lisses
$$1\rta J_a \rta J_{1,a}\rta C\rta 1.$$
On en déduit une suite exacte
$$1\rta J_a(F_x)/J_a(\calO_x)\rta J_{1,a}(F_x)/J_{1,a}(\calO_x)
\rta C(F_x)/C(\calO_x) \rta 1.
$$
Puisque $C$ est un tore sur $X$, $C(F_x)/C(\calO_x)\simeq \ZZ^c$
où $c$ est le rang de $C$. En particulier, $C(F_x)/C(\calO_x)$ est
un groupe discret. Il s'ensuit que la composante neutre de
$J_a(F_x)/J_a(\calO_x)$ s'envoie bijectivement sur la composante
neutre de $J_{1,a}(F_x)/J_{1,a}(\calO_x)$. Par conséquent, la
flèche
$$\pi_0(J_a(F_x)/J_a(\calO_x))\rta
\pi_0(J_{1,a}(F_x)/J_{1,a}(\calO_x))
$$
est injective. Après dualisation, on obtient un {\em homomorphisme
surjectif}
$$\hat T_1^{I_x}=\pi_0(J_{1,a}(F_x)/J_{1,a}(\calO_x))^* \rta
\pi_0(J_a(F_x)/J_a(\calO_x))^*$$ où $\hat T_1$ est le tore maximal
de $\hat G_1$. Puisque $\kappa$ appartient à l'image
$$\pi_0(J_a(F_x)/J_a(\calO_x))^* \rta
\pi_0(J_a(F_x)/J_a^0(\calO_x))^*=\hat T^{I_x},
$$
il s'ensuit que $\kappa$ appartient à l'image de l'application
composée
$$\hat T_1^{I_x} \rta \pi_0(J_a(F_x)/J_a(\calO_x))^* \rta \hat T^{I_x}.$$
Il suffit maintenant d'invoquer un lemme de Kottwitz pour terminer
la dé\-mon\-stration du théorème.

\begin{lemme} Si $\kappa$ appartient à l'image de $\hat T_1^{I_x}
\rta \hat T^{I_x}$, alors l'homomorphisme $I_x\rta {\rm
Out}(\mathbbm H)$ est trivial.
\end{lemme}

\dem Pour la commodité du lecteur, on rappellera l'argument de
Kottwitz avec nos notations légèrement différentes des siennes.

Rappelons qu'on a une suite exacte
$$1\rta \hat C\rta \hat G_1\rta
\hat G \rta 1.$$ Notons $\hat T_1$ l'image réciproque de $\hat T$
dans $\hat G_1$. Les plongements $\eta_{\hat T}:\hat T\rta \hat G$
et $\eta_{\hat T_1}:\hat T_1\rta \hat G_1$ sont $I_x$-équivariant
à conjugaison près. Contenu du fait que $I_x$ agit trivialement
sur $\hat G$ et $\hat G_1$ pour tout $\sigma\in I_x$, il existe
$g_1^\sigma\in \hat G_1$ tel que pour tous $t_1\in \hat T_1$, on a
$$\eta_{\hat T_1}(\sigma(t_1))={\rm Int}(g_1^\sigma)
(\eta_{\hat T_1}(t_1)).
$$
Soit $g^\sigma$ l'image de $g_1^\sigma$ dans $\hat G$, on a alors
la relation similaire
$$\eta_{\hat T}(\sigma(t))={\rm Int}(g^\sigma)(\eta_{\hat T}(t))$$
pour tout $t\in \hat T$.

Soit $\kappa_1\in \hat T_1^{I_x}$ ayant $\kappa$ comme image dans
$\hat T^{I_x}$. Soit $\hat H$ la composante neutre du
centralisateur de $\kappa$ dans $\hat G$ et $\hat H_1$ la
composante neutre du centralisateur de $\kappa_1$ dans $\hat G_1$.
Puisque $\hat G_1$ est simplement connexe, $\hat H_1$ est tout le
centralisateur de $\kappa_1$. L'homomorphisme $\hat H_1\rta \hat
H$ est surjectif car il induit un isomorphisme sur les algèbres de
Lie. Nous allons fixer des épinglages compatibles de $\hat H, \hat
H_1$ et de  $\hat G, \hat G_1$ avec tores maximaux $\hat T, \hat
T_1$. D'après Langlands, il existe un unique action de $I_x$ sur
$\hat H_1$ et $\hat H$ fixant ces épinglages tels que les
plongements $\hat T_1\rta \hat H_1 \rta \hat G_1$ et $\hat T\rta
\hat H\rta \hat G$ sont $I_x$-équivariants à conjugaison près. De
plus l'homomorphisme surjectif $\hat H_1\rta \hat H$ est
$I_x$-équivariant. Notons que les faits que $\hat T_1\rta \hat
H_1$ est $I_x$-équivariant par conjugaison et que $\kappa_1\hat
T_1^{I_x}$ est centralisé par $\hat H_1$ impliquent que $\kappa_1$
est aussi fixé par l'action de $I_x$ sur $\hat H_1$.

Il s'agit de démontrer que pour tout $\sigma\in I_x$, l'image de
$\sigma$ dans ${\rm Out}(\hat H)$ est trivial. Il suffit de
démontrer que son image dans ${\rm Out}(\hat H_1)$ est trivial.
Or, puisque l'inclusion $\eta_{\hat H_1}:\hat H_1\rta\hat G_1$ est
$I_x$-équivariant à conjugaison près, il existe $g_1^\sigma\in
\hat G_1$ tel que
$$\eta_{\hat H_1}(\sigma(h_1))={\rm Int}(g_1^\sigma)
(\eta_{\hat H_1}(h_1))
$$
pour tout $h\in \hat H_1$. Appliquer cette relation à $\kappa_1$
qui est fixé par $\sigma$, cette relation implique que $g_1$
centralise $\kappa_1$. Puisque $\hat H_1$ est tout le
centralisateur de $\kappa_1$, ceci implique que l'action de
$\sigma$ sur $\hat H_1$ est intérieur.

Ceci termine la démonstration du lemme 10.5 et donc celle du
théorème 10.4. \findem

\bigskip
Cette description du support de la partie $\kappa$-isotypique joue
le rôle d'un énoncé de pureté, analogue à la conjecture de pureté
de Goresky-Kottwitz-MacPherson \cite{GKM}, importante dans leur
approche du lemme fondamental par la cohomologie équivariante.

Au lieu de démontrer la pureté de la cohomologie des fibres de
Springer conjecturée par Goresky, Kottwitz et MacPherson, on va
démontrer que les images réciproques de $^p\rmH^j(f^{\rm
ell}_*\QQ_\ell)_{[\kappa]}$ aux espaces affines de Hitchin des
groupes endoscopiques $\bbA_{H}^{\rm max}$ sont pures. Dans le cas
unitaire, cet énoncé a été suffisant pour démontrer le lemme
fondamental.

En utilisant un résultat de Faltings \cite{Fa}, on peut vérifier
qu'au-dessus de l'ouvert elliptique $\bbA^{\rm ell}$, le champ
$\calM^{\rm ell}$ est un champ de Deligne-Mumford lisse, et que le
morphisme de Hitchin $f^{\rm ell}:\calM^{\rm ell} \rta \bbA^{\rm
ell}$ est propre et de type fini. Alors, les facteurs directs
$^p\rmH^j(f^{\rm ell}_*\QQ_\ell)_{[\kappa]}$ de l'image directe,
sont purs par un théorème de Deligne (\cf \cite{D}). Maintenant,
les images réciproques aux $\bbA_{H}^{\rm ell}$ sont pures parce
que $\bbA_{H}^{\rm max}$ est fini et étale au-dessus de son image
dans $\bbA^{\rm ell}$ et que le facteur direct $^p\rmH^j(f^{\rm
ell}_*\QQ_\ell)_{[\kappa]}$ est supporté par la réunion disjointe
des adhérences des images des $\bbA_{H_\rho}^{\rm max}$ pour les
différentes classes de de données endoscopiques non ramifiées
$(\kappa,\rho)$.

Notons enfin l'énoncé suivant qui découle immédiatement de la
description de Donagi-Gaitsgory du centralisateur régulier, voir
proposition 3.8.

\begin{proposition}
Soit $p^H_G:\car_H\rta \car_G$ le morphisme canonique de l'espace
des caractéristiques de $H$ dans celui de $G$. On a alors un
homomorphisme canonique
$$(p^H_G)^* J_G \rta J_H$$
entre les centralisateur réguliers qui est un isomorphisme
au-dessus du lieu régulier semi-simple.
\end{proposition}

Ces ingrédients ont été utilisés dans la démonstration du lemme
fondamental pour le groupe unitaire \cite{L-N} et sont donc
établis ici en toute généralité. Il y a d'autres ingrédients de
loc. cit. qui ne peuvent pas se généraliser. En particulier, dans
loc. cit. on a utilisé un plongement $H\rta G$ qui localement pour
la topologie étale de $X$ fait de $H$ un sous-groupe de Levi de
$G$. Ce plongement n'existe pas en général et il faudra suppléer
d'autres arguments.

\section{Le cas $\SL(2)$}

On discute dans cette dernière section du cas simple du groupe
$\SL(2)$ et redémontre le théorème 10.2 dans ce cas-ci. Pour cela,
on donnera une description de la flèche $\rmH^0(\ovl
X,J_a/J_a^0)\rta \pi_0(P'_a)$ en utilisant les modèles de Néron
plutôt que la contourner par le lemme 10.3.

Pour $G=\SL(2)$, l'espace affine de Hitchin est l'espace
$$\bbA(\ovl k)=\rmH^0(\ovl X,\calO_X(2D)).$$ Pour tout $a\in
\rmH^0(\ovl X,\calO_X(2D))$, on trace une courbe $Y_a$ d'équation
$t^2-a=0$ sur l'espace total $V_D$ du fibré en droites
$\calO_X(D)$ au-dessus de $\ovl X$. Cette courbe vient avec un
morphisme fini $p_a:Y_a\rta \ovl X$ de degré $2$ dont le groupe
des automorphismes est $\{1,\tau\}$. L'involution $\tau$ agit par
$\tau(t)=-t$ et a pour points fixes dans $Y_a$ les points
d'intersection de $Y_a$ avec la section zéro de $V_D$. Notons
$U_a$ l'ouvert maximal de $\ovl X$ au-dessus duquel $Y_a$ est
étale : c'est le complément des points dans $X$ images des points
d'intersection de $Y_a$ avec la section nulle. La caractéristique
$a\in\bbA(\ovl k)$ est elliptique si et seulement si $Y_a$ est
irréductible. Supposons que $a$ est elliptique.

Le schéma en groupes $J_a=h_a^*[J]$ défini dans la section 4
associe à tout $X$-schéma $S$ le groupe
$$J_a(S)=\{s\in\Gamma(S\times_X Y_a,\calO^\times)\mid \tau(s)s=1\}.$$
La restriction de $J_a$ à $U_a$ est donc un tore et ses fibres en
les points $x\in \ovl X-U_a$ ont pour groupe de composantes
connexes $\ZZ/2\ZZ$. Soit $J_a^0$ le sous-schéma en groupes de
$J_a$ de ses composantes neutres. On a alors la suite exacte de
faisceaux
$$0\rta J_a^0 \rta J_a \rta \bigoplus_{x\in\ovl X-U_a}
(\ZZ/2\ZZ)_{x}\rta 0$$ qui induit une suite longue
$$\rmH^0(\ovl X,\bigoplus_{x\in\ovl X-U_a}
(\ZZ/2\ZZ)_{x}) \rta \rmH^1(\ovl X,J_a^0) \rta \rmH^1(\ovl
X,J_a)\rta 0.
$$
Cette suite induit une suite exacte à droite des $\pi_0$
$$\rmH^0(\ovl X,\bigoplus_{x\in\ovl X-U_a}
(\ZZ/2\ZZ)_{x}) \rta \pi_0(P'_a) \rta \pi_0(P_a)\rta 0$$ où $P'_a$
est le champ de Picard des $J_a^0$-torseurs. On sait d'après 6.4
que $\pi_0(P'_a)=\ZZ/2\ZZ$, le groupe des coinvariants de
$\pi_1(U_a,u)$ dans $\bbX^\sv=\ZZ$ de sorte que $\pi_0(P_a)$ est
trivial ou est égal à $\ZZ/2\ZZ$ dépendant de la nullité ou de la
non-nullité de la flèche
$$\beta:\rmH^0(\ovl X,\bigoplus_{x\in\ovl X-U_a}
(\ZZ/2\ZZ)_{x}) \rta \pi_0(P'_a).$$ La flèche $\beta$ est une
somme des flèches de bord
$$\beta_x:\pi_0(J_{a,x})\rta \pi_0(P'_a)=(\bbX^\sv_u)_{\pi_1(U_a,u)}$$
qui peuvent être décrites de façon relativement concrète. Soit
${\rm Ner}(J_a|_{U_a})$ le modèle de Néron du tore $J_a|_{U_a}$.
Par la propriété universelle des modèles de Néron, on a un
homomorphisme
$$J_a\rta {\rm Ner}(J_a|_{U_a})$$
qui induit en particulier des flèches
$$\pi_0(J_{a,x})\rta \pi_0({\rm Ner}(J_a|_{U_a})).$$
D'après Kottwitz et Rapoport (\cf \cite{Ra} 2.2 (iii)), le groupe
des composantes de la fibre en $x$ du modèle de Néron est le
groupe des coinvariants sous l'inertie du groupe des cocaractères
$$\pi_0({\rm Ner}(J_a|_{U_a})_x)=(\bbX^\sv_{u})_{I_x}$$
où $I_x\subset \pi_1(U_a,u)$ est le sous-groupe d'inertie en $x$.
La flèche $\beta_v$ est alors le composé de la flèche purement
locale $\pi_0(J_{a,x})\rta \pi_0({\rm Ner}(J_a|_{U_a}))$ et de la
flèche évidente $(\bbX^\sv_{u})_{I_x} \rta
(\bbX^\sv_u)_{\pi_1(U_a,u)}$.

 \begin{lemme} Soit $a\in\bbA^{\rm ell}(\ovl k)$ une
 caractéristique elliptique pour $\SL(2)$.
Alors $\pi_0(P_a)=0$  dès qu'il existe un point de ramification de
$\pi_a:Y_a\rta X$ qui est unibranche. Dans le cas où tous les
points de ramification de $Y_a\rta X$ ont deux branches, on a
$\pi_0(P_a)=\ZZ/2\ZZ$.
 \end{lemme}

 \dem Il revient au même de démontrer que la flèche de bord
 $$\beta_x:(\ZZ/2\ZZ)_x \rta \pi_0(P'_a)=\ZZ/2\ZZ$$
 est non nulle si et seulement si le point de $Y_a$ au-dessus de
 $x$ est unibranche. Examinons les deux cas.

\medskip
{\em Le cas d'un point de ramification à deux branches} : Soit $x$
un point de branchement de $\ovl X$ au-dessus duquel le point
$\widetilde x$ de $Y_a$ a deux branches. Soit $\calO_x$ l'anneau
local complété de $\ovl X$ en $x$ et $F_x$ son corps des
fonctions. Soit $\calO_{Y_a,\widetilde x}$ le complété de $Y_a$ en
$\widetilde x$ et $E_x$ l'anneau total des fractions de
$\calO_{Y_a,\widetilde x}$. Puisque $Y_a$ est deux branches en
$x$, $E_x$ est un produit de deux corps. Puisque $E_x/F_x$ est de
degré $2$, on a $E_x=F_x\times F_x$. On en déduit que l'inertie
$I_x$ agit trivialement sur $\bbX^\vee$ de sorte que
$\bbX^\vee_{I_x}=\ZZ$.

La flèche $\beta_x:(J_{a,x}/J^0_{a,x}) \rta
\bbX^\vee_{\pi_1(U_a)}$ se factorise à travers $\bbX^\vee_{I_v}$.
Or la flèche $(J_{a,x}/J^0_{a,x}) \rta \bbX^\vee_{I_x}$ est
nécessairement nulle car $(J_{a,x}/J^0_{a,x}) =\ZZ/2\ZZ$ est un
groupe fini et $\bbX^\vee_{I_x}=\ZZ$ n'a pas de torsion. Donc
$\beta_x=0$ dans le cas de deux branches.

\medskip
{\em Le cas d'un point de ramification à une seule branche} :
Mettons-nous maintenant au voisinage d'un point de branchement $x$
tel qu'en le point $\widetilde x$ de $Y_a$ au-dessus de $\ovl X$,
$Y_a$ est unibranche. Gardons les mêmes notations que dans le cas
précédent. : $E_x$ est une extension ramifiée de degré $2$ de
$F_x$ et le groupe d'inertie $I_x$ agit sur $\bbX^\vee=\ZZ$ par
$d\mapsto -d$. La flèche $\bbX^\sv_{I_x}\rta
\bbX^\sv_{\pi_1(U_a,u)}$ est donc un isomorphisme.

Le groupe $J(F_x)$ est l'ensemble des éléments de norme $1$ de
$E^\times_x$. On constate que $J(F_x)$ est contenu dans les
entiers de $E_x$. En particulier, il est borné. Le modèle Néron
${\rm Ner}(J_a|_{U_a})$ est donc de type fini et a pour
$\calO_v$-points les points
$${\rm Ner}(J_a|_{U_a})(\calO_x)=J(F_x).$$
Sa fibre spéciale a deux composantes connexes selon que la section
$s$ prend la valeur $1$ ou $-1$ en le point $\widetilde x$ de
$Y_a$.

On a $J_a(\calO_v)\subset {\rm Ner}(J_a|_{U_a})(\calO_v)$ et
$J_a^0(\calO_v)\subset {\rm Ner}^0(J_a|U_a)(\calO_v)$.
L'homomorphisme canonique
$$J_a(\calO_v)/J^0_a(\calO_v)\rta {\rm Ner}(J_a|_{U_a})(\calO_v)/
{\rm Ner}(J_a|_{U_a})^0(\calO_v)
$$
est un isomorphisme. En effet on a un triangle commutatif
$$\xymatrix{
  J_a(\calO_v)/J^0_a(\calO_v) \ar[rr]^{} \ar[dr]_{}
                &  &   {\rm Ner}(J_a|_{U_a})(\calO_v)/
{\rm Ner}(J_a|_{U_a})^0(\calO_v)  \ar[dl]^{}    \\
                & \{\pm 1\}                }
$$
où les deux flèches en biais associent à la classe d'une section
$s$ de $J_a$ ou de ${\rm Ner}(J_a|_{U_a})$, définie au voisinage
de $\widetilde x$, sa valeur prise en $s$. Celles-ci étant des
isomorphismes, il en est de même de la flèche horizontale. Dans le
cas d'un point de ramification unibranche, la flèche $\beta_v$ est
un donc isomorphisme.\findem

\begin{lemme} Soit $Y_a$ un revêtement double de $\ovl X$ génériquement
étale dont tous les points de ramifications ont deux branches.
Alors, la normalisation $\widetilde Y_a$ de $Y_a$ est un
revêtement partout non ramifié de $X$.
\end{lemme}

\dem Après normalisation, au-dessus d'un point de branchement
$x\in\ovl X$, il y a deux points de $\widetilde Y_a$ correspondant
aux deux branches passant par le point de ramification au-dessus
de $x$. Le revêtement $\widetilde Y_a\rta X$ est alors
nécessairement étale. \findem

\bigskip La donnée d'un revêtement double étale de $\ovl X$
est équivalente à la donnée d'un homomorphisme $\rho:\pi_1(\ovl
X,u)\rta \ZZ/2\ZZ$. Elle est aussi équivalente à la donnée d'un
fibré inversible $\calL_\rho$ sur $\ovl X$ muni d'un isomorphisme
$\calL_\rho^{\otimes 2}\isom \calO_{\ovl X}$. Soit
$V_{\calL_\rho}\rta \bbA^1_X$ le morphisme de l'espace total de
$\calL_\rho$ dans l'espace total de $\calO_{\ovl X}$ donné sur les
sections locales par $s\mapsto s^{\otimes 2}$. Alors, on retrouve
le revêtement double étale de $\ovl X$ en prenant l'image
réciproque de la section $1$ de $\bbA^1_X$.

\begin{lemme}
Soit $a\in\rmH^0(\ovl X,\calO_X(2D))$ tel que la courbe $Y_a$
tracée sur l'espace total $V_D$ du fibré vectoriel $\calO_X(D)$,
n'a que des points de ramification à deux branches. Soit
$\rho:\pi_1(\ovl X,u)\rta\ZZ/2\ZZ$ le quotient du groupe
fondamental correspondant au revêtement non ramifié de $\ovl X$
obtenu en normalisant $Y_a$. Alors le diviseur d'annulation ${\rm
div}(a)$ est de la forme
$${\rm div}(a)=2D'$$
$D'$ étant un diviseur effectif tel que $\calO_{\ovl X}(D'-D)$
soit isomorphe à $\calL_\rho$ où $\calL_\rho$ est le fibré
inversible de carré neutre, associé à $\rho$. Il existe de plus
une section
$$b\in \rmH^0(\ovl X,\calL_\rho\otimes \calO_X(D))$$ bien déterminée à un signe près,
tel que $b^{\otimes 2}=a$.
\end{lemme}

\dem On sait que localement le revêtement $Y_a$ est déterminé par
une équation de la forme $t^2-a=0$ où $t$ est la coordonnée
verticale du fibré en droites et où la section $a$ est vue
localement comme une fonction sur $X$ à l'aide de trivialisation
locale de $\calO_X(2D)$. Les points de branchement $x\in\ovl X$
sont les points où $a$ s'annulent. Le germe de $Y_a$ au-dessus
d'un point de branchement $x$, est unibranche si $a$ s'annule en
$x$ à l'ordre impair, et $Y_a$ a deux branches si $a$ s'annule en
$x$ à l'ordre pair. Ainsi $Y_a$ n'a que des points de ramification
à deux branches si et seulement si le diviseur d'annulation de $a$
est de la forme
$${\rm div}(a)=2D'$$
où $D'$ est un diviseur effectif. Soit $\calL=\calO_{\ovl
X}(D'-D)$. Alors la section $a$ de $\calO_X(2D)$ ayant $2D'$ comme
diviseur d'annulation, définit un isomorphisme $\calL^{\otimes
2}\isom \calO_{\ovl X}$. La section $a$ fournit par ailleurs un
isomorphisme de
$$V_{\calO_X(2D)}\times_{\ovl X}(\ovl X-D')\isom \bbA^1_{\ovl
X-D'}
$$
au-dessus duquel on a un isomorphisme
$$V_{\calO_X(D)}\times_{\ovl X}(\ovl X-D') \isom
V_{\calL}\times_{\ovl X}(\ovl X-D').
$$
Il s'ensuit qu'au-dessus de $\ovl X-D'$, le revêtement double
étale $Y_a$ et celui qui se déduit du fibré de carré neutre
$\calL$ sont isomorphes. L'isomorphisme se prolonge sur tout $\ovl
X$ par normalisation.

L'isomorphisme $\calL\otimes \calO_{\ovl X}(D)\isom \calO_{\ovl
X}(D')$ fournit à $\calL\otimes \calO_{\ovl X}(D)$ une section $b$
dont le diviseur d'annulation est $D'$. La section $b^{\otimes 2}$
de $\calO_X(2D)$ a alors le même diviseur d'annulation que $a$. En
modifiant $b$ par un scalaire, on a $b^{\otimes 2}=a$. \findem

\bigskip
La conjonction des trois lemmes précédents fournit la description
suivante du faisceau $\pi_0(P/\bbA^{\rm ell})$. Soit ${\rm
Pic}_{\ovl X}[2]^*$ l'ensemble des points d'ordre deux non
triviaux de la jacobienne de $\ovl X$. Pour chaque $\calL_\rho\in
{\rm Pic}_{\ovl X}[2]^*$, notons $H_\rho$ le tore endoscopique
elliptique au-dessus de $\ovl X$ défini comme la forme extérieure
de $\GG_m$ associée à la représentation du groupe fondamental
$\rho:\pi_1(\ovl X,u)\rta\ZZ/2\ZZ$. Son espace affine de Hitchin
est
$$\bbA_{H_\rho}=\rmH^0(\ovl X,\calL_\rho\otimes \calO_{\ovl
X}(D)).
$$
Notons ${\mathbb S}_\rho$ l'image de $\bbA_{H_\rho}^{\rm ell}$
dans $\bbA^{\rm ell}$. Alors $\pi_0(P/\bbA^{\rm ell})$ est
supporté par la réunion disjointe
$$\bigsqcup_{\calL_\rho\in
{\rm Pic}_{\ovl X}[2]^*} {\mathbb S}_{\rho}.$$ Sa restriction à
chaque ${\mathbb S}_\rho$ est le faisceau constant $\ZZ/2\ZZ$.

Le morphisme de faisceaux
$$\ZZ/2\ZZ=\pi_0(P'/{\bbA^{\rm ell}}) \rta \pi_0(P/\bbA^{\rm ell})$$
définit une action de $\ZZ/2\ZZ$ sur $^p\rmH^j(f^{\rm
ell}_*\QQ_\ell)$ et décompose ce faisceau pervers en somme de deux
composantes isotypiques
$$^p\rmH^j(f^{\rm ell}_*\QQ_\ell)=\,^p\rmH^j(f^{\rm
ell}_*\QQ_\ell)_+\oplus \,^p\rmH^j(f^{\rm ell}_*\QQ_\ell)_-.
$$
On peut vérifier directement dans le cas $\SL(2)$ que $\calM^{\rm
ell}$ est un champ de Deligne-Mumford lisse et que le morphisme
$f^{\rm ell}:\calM^{\rm ell}\rta \bbA^{\rm ell}$ est un morphisme
propre de type fini. Il s'ensuit que le faisceau pervers
$^p\rmH^j(f^{\rm ell}_*\QQ_\ell)$ est pur, ainsi que ses facteurs
directs $^p\rmH^j(f^{\rm ell}_*\QQ_\ell)_+$ et $^p\rmH^j(f^{\rm
ell}_*\QQ_\ell)_-$. Le morceau $^p\rmH^j(f^{\rm ell}_*\QQ_\ell)_-$
est supporté par $\bigsqcup_{\calL_\rho\in {\rm Pic}_{\ovl
X}[2]^*} {\mathbb S}_{\rho}$ de sorte que son image réciproque à
chaque ${\mathbb S}_\rho$ est pure.

\end{document}